\documentclass[12pt]{article}\usepackage{amsfonts,amssymb}

\usepackage{amsmath,bbm}
\usepackage{theorem}          
\usepackage{latexsym}
\usepackage{float}                   
\usepackage{times}
\usepackage{cite}

\newtheorem{tht}{Theorem}[section]

\newtheorem{thl}[tht]{Lemma}
\newtheorem{thp}[tht]{Proposition}
\newtheorem{thc}[tht]{Corollary}

\newcommand{\anf}{\raisebox{0.2ex}{\scriptsize$\triangleleft$}}
\newcommand{\ang}{\raisebox{0.2ex}{\scriptsize$\triangleright$}}

\newcommand{\mn}{\medskip}    

\newcommand{\rti}{\,{\scriptstyle\rtimes}\,} %

\newcommand{\sumop}{\mathop{\mbox{$\sum$}}}

\newcommand{\hsp}{\hspace{-1pt}} 
\newcommand{\hs}{\hspace{1pt}} 

\newcommand{\cD}{{\mathcal{D}}}

\newcommand{\cG}{{\mathcal{G}}} 
\newcommand{\Hh}{{\mathcal{H}}} 
\newcommand{\cO}{{\mathcal{O}}}
\newcommand{\cK}{{\mathcal{K}}}
\newcommand{\cU}{{\mathcal{U}}}   

\newcommand{\cX}{{\mathcal{X}}}
\newcommand{\cA}{{\mathcal{A}}}

\newcommand{\cC}{{\mathcal{C}}}
\newcommand{\cF}{{\mathcal{F}}}

\newcommand{\cS}{{\mathcal{S}}} 
\newcommand{\cL}{{\mathcal{L}}}  

\newcommand{\cM}{{\mathcal{M}}}
\newcommand{\cY}{{\mathcal{Y}}}

\newcommand{\dZ}{{\mathbb{Z}}}
\newcommand{\dR}{{\mathbb{R}}}
\newcommand{\dC}{{\mathbb{C}}}
\newcommand{\dN}{{\mathbb{N}}}

\newcommand{\fm}{{\mathfrak{m}}}
\newcommand{\fn}{{\mathfrak{n}}}

\newcommand{\Lin}{{\mathrm{Lin}}}

\newcommand{\sign}{{\mathrm{sign}}}

\newcommand{\inv}{{\mathrm{inv}}}
\newcommand{\tr}{{\mathrm{Tr}}}

\newcommand{\SU}{\cO({\mathrm{SU}}_q(2))}

\newcommand{\su}{\cU_q({\mathrm{su}}_2)}

\begin{document}

\date{\small{Fakult\"at f\"ur Mathematik und 
Informatik\\ Universit\"at Leipzig, 
Augustusplatz 10, 04109 Leipzig, Germany\\ 
E-mail: Konrad.Schmuedgen@math.uni-leipzig.de / 
Elmar.Wagner@math.uni-leipzig.de}
}

\title{Representations of cross product algebras 
of Podle\'s quantum spheres }

\author{Konrad Schm\"udgen  and Elmar Wagner}

\maketitle

\renewcommand{\theenumi}{\roman{enumi}}

\begin{abstract}\noindent
Hilbert space representations of the 
cross product $\ast$-al\-ge\-bras
of the Hopf $\ast$-al\-ge\-bra 
$\cU_q(\mathrm{su}_2)$ and 
its module $\ast$-al\-ge\-bras $\cO(\mathrm{S}^2_{qr})$ 
of Podle\'s spheres are investigated and classified by describing the action 
of generators. The representations are analyzed 
within two approaches. 
It is shown that the Hopf $\ast$-al\-ge\-bra 
$\cO (\mathrm{SU}_q(2))$  of the quantum group $\mathrm{SU}_q(2)$ 
decomposes into an orthogonal sum 
of  projective Hopf modules 
corresponding to irreducible integrable $\ast$-re\-pre\-sen\-ta\-tions of 
the cross product algebras 
and that each irreducible integrable $\ast$-re\-pre\-sen\-ta\-tion 
appears with multiplicity one.  
The projections of these projective modules 
are computed. The decompositions of tensor products 
of irreducible integrable $\ast$-re\-pre\-sen\-ta\-tions 
with spin $l$ representations 
of $\cU_q(\mathrm{su}_2)$ are given. The invariant state $h$ on 
$\cO(\mathrm{S}^2_{qr})$ is studied 
in detail. By passing to 
function algebras over the quantum spheres $\mathrm{S}^2_{qr}$,  
we give chart descriptions of quantum line bundles and describe 
the representations from the first approach by means of 
the second approach. 
\end{abstract}
%
%
%
{\small Keywords: Quantum groups, unbounded representations}
%
%
%
%

\setcounter{section}{-1}
\section{Introduction}
Podle\'s quantum spheres \cite{P} (see \cite[Section 4.5]{KS} for 
a short treatment) 
are a one-parameter family $\mathrm{S}^2_{qr}$, $r\in [0,\infty]$, 
of mutually non-isomorphic 
quantum homogeneous spaces of the quantum group $\mathrm{SU}_q(2)$, 
where $0 < q < 1$. 
Each of these spaces can be considered as a quantum analogue of the  classical 
2-sphere. Their coordinate algebras $\cO(\mathrm{S}^2_{qr})$ are right coideal 
$\ast$-sub\-al\-ge\-bras 
of the coordinate Hopf $\ast$-al\-ge\-bra $\cO(\mathrm{SU}_q(2))$ 
of the quantum 
group $\mathrm{SU}_q(2)$ 
and hence left module $\ast$-al\-ge\-bras of 
the Hopf $\ast$-al\-ge\-bra $\cU_q(\mathrm{su}_2)$. 
Therefore, the left cross product $\ast$-al\-ge\-bra 
$\cO(\mathrm{S}^2_{qr}) \rti \cU_q(\mathrm{su}_2)$ 
is defined. The subject of this paper are Hilbert space representations 
of the 
$\ast$-al\-ge\-bra $\cO(\mathrm{S}^2_{qr}) \rti \cU_q(\mathrm{su}_2)$.

There is now an extensive literature on Podle\'s   spheres 
(see e.g.\ \cite{B,BM,DK,H,HM,MS}). 
Let us restate some algebraic results from these papers that are relevant 
for our 
investigations. Let $\cC$ denote the coalgebra 
$\cO(\mathrm{SU}_q(2)) / \cO(\mathrm{S}^2_{qr})^+\cO(\mathrm{SU}_q(2)) $ 
with quotient map $\rho : \cO(\mathrm{SU}_q(2)) \rightarrow \cC$, 
where 
 $\cO(\mathrm{S}^2_{qr})^+ = \{x \in \cO(\mathrm{S}^2_{qr})\,;\,
\varepsilon (x) = 0\}$. 
The algebra 
$\cO(\mathrm{S}^2_{qr})$ is then the algebra of left 
$\cC$-co\-var\-i\-ant elements of 
$\cO(\mathrm{SU}_q(2))$. 
Note that only in the case $r=0$ the coalgebra $\cC$ is a Hopf algebra and 
the quantum 
sphere $\mathrm{S}^2_{qr}$ is the quotient by a quantum subgroup. 
M. Dijkhuizen and T. Koornwinder 
\cite{DK} have found a skew-primitive element 
$X_r \in \cU_q(\mathrm{su}_2)$ such that 
$\cO(\mathrm{S}^2_{qr})$ 
is the subalgebra of right $X_r$-in\-var\-iant elements of 
$\cO(\mathrm{SU}_q(2))$. 
A major 
step have been the results 
of E. F. M\"uller and H.-J. Schneider \cite{MS}. 
They showed that $\cO(\mathrm{SU}_q(2))$ is faithfully flat as 
a left (and right) 
$\cO(\mathrm{S}^2_{qr})$-mo\-dule and that $\cC$ is spanned by group-like 
elements. As a consequence, 
$\cC$ is the direct sum of simple subalgebras $\cC_j$. 
Then
$$
M_j = \{x \in \cO(\mathrm{SU}_q(2))\,;\, 
\rho (x_{(1)}) \otimes x_{(2)} \in \cC_j \otimes \cO(\mathrm{SU}_q(2)) \}
$$
is a finitely generated projective relative 
$(\cO(\mathrm{SU}_q(2)),\ \cO(\mathrm{S}^2_{qr}))$-Hopf module 
and $\cO(\mathrm{SU}_q(2))$ is the direct sum of these 
Hopf modules $M_j$ \cite[p.\ 186]{MS}. In the 
subgroup case $r=0$, the corresponding projections and 
their Chern numbers have been 
computed in \cite{HM} and \cite{H}, respectively. A family of group-like 
elements spanning 
the coalgebra $\cC$ was determined in \cite{BM}. 

In the present paper, we reconsider and extend these algebraic results 
in the Hilbert 
space setting. We prove that $\cO(\mathrm{SU}_q(2))$ is 
the orthogonal direct sum of the Hopf modules 
$M_j$ and  we give an 
explicit description of this decomposition. 
Here the skew-primitive element $X_r$ plays a crucial role. 
Moreover, we determine the 
projections of the projective modules $M_j$. 
All this is carried out in Section \ref{sec-decomp}. 
Each Hopf module $M_j$ corresponds to 
an irreducible $\ast$-re\-pre\-sen\-ta\-tion $\pi_j$ 
of the cross product $\ast$-al\-ge\-bra 
$\cO(\mathrm{S}^2_{qr}) \rti \cU_q(\mathrm{su}_2)$ such that the 
restriction of $\pi_j$ to $\cU_q(\mathrm{su}_2)$ is a direct sum of spin $l$ 
representations $T_l$, 
$l\in \frac{1}{2}\dN_0$. Let us call a $\ast$-re\-pre\-sen\-ta\-tion of 
$\cO(\mathrm{S}^2_{qr}) \rti \cU_q(\mathrm{su}_2)$ 
integrable if it has the latter property. 
In Section \ref{S5},  
we classify integrable 
$\ast$-re\-pre\-sen\-ta\-tions of 
$\cO(\mathrm{S}^2_{qr}) \rti \cU_q(\mathrm{su}_2)$ and prove that each 
irreducible integrable $\ast$-re\-pre\-sen\-ta\-tion of 
$\cO(\mathrm{S}^2_{qr}) \rti \cU_q(\mathrm{su}_2)$ 
is unitarily equivalent to one of the representations $\pi_j$. 
In the course of this classification,  we describe the structure  of these 
$\ast$-re\-pre\-sen\-ta\-tions $\pi_j$ by explicit 
formulas for the actions of generators of $\cO(\mathrm{S}^2_{qr})$ on 
an orthonormal basis of 
weight vectors for the representation $T_l$ of $\cU_q(\mathrm{su}_2)$. 
In the terminology of 
our previous paper \cite{SW}, this is the first approach 
to representations of the 
cross product algebra $\cO(\mathrm{S}^2_{qr}) \rti \cU_q(\mathrm{su}_2)$.
We also derive a formula for the decomposition of 
the tensor product representation 
$\pi_j \otimes T_l$ 
into a direct sum of representations $\pi_j$.

The corresponding second approach \cite{SW} is 
developed in Section \ref{sec-2app}. 
Here we begin with a  
representation  of the $\ast$-al\-ge\-bra $\cO(\mathrm{S}^2_{qr})$ given 
in a canonical form and  
we extend it to a representation of 
$\cO(\mathrm{S}^2_{qr}) \rti \cU_q(\mathrm{su}_2)$. The main 
technical tool  for this is to ``decouple'' the cross relations 
of the cross product algebra by finding   
an auxiliary $\ast$-sub\-al\-ge\-bra $\cY_r$ of 
$\cO(\mathrm{S}^2_{qr}) \rti \cU_q(\mathrm{su}_2)$ which 
commutes 
with the $\ast$-al\-ge\-bra $\cO(\mathrm{S}^2_{qr})$. 
Such decouplings have been found 
and studied in \cite{F}.

Section \ref{***} starts by defining algebras of functions which extend the 
coordinate algebras $\cO(\mathrm{S}^2_{qr})$ and by describing 
invariant functionals on such function algebras. The algebras of functions 
together with the invariant functionals will be used to give 
another description of irreducible integrable $\ast$-re\-pre\-sen\-ta\-tions. 
It should be emphasized that though all irreducible 
integrable representations $\pi_j$ of 
$\cO(\mathrm{S}^2_{qr}) \rti \cU_q(\mathrm{su}_2)$ involve unbounded
operators, their restrictions to the $\ast$-sub\-al\-ge\-bra 
$\cO(\mathrm{S}^2_{qr})\otimes\cY_r$ are given by 
bounded operators only.
In Subsection \ref{sec-repcpa}, we show that the restriction of an 
irreducible integrable $\ast$-re\-pre\-sen\-ta\-tion to 
$\cO(\mathrm{S}^2_{qr})\otimes\cY_r$ decomposes into a direct sum of 
two representations which can be realized on algebras of functions 
with support in the positive and the negative spectrum of a certain 
self-adjoint operator. 
This self-adjoint operator represents a coordinate function of the 
quantum sphere and the  two algebras of functions can be 
considered as ``charts'' of the projective module $M_j$.  
The representation of 
$\cO(\mathrm{S}^2_{qr})\otimes\cY_r$ 
on each chart leads again to an irreducible $\ast$-re\-pre\-sen\-ta\-tion of 
the cross product algebra which is not integrable and can be described 
by the formulas from the second approach. 
To round off this circle of investigations, we recover in Section \ref{last}  
the irreducible integrable representation
by taking the direct sum of both charts and passing to another domain.

In Section \ref{S1}, we briefly mention the correspondence 
between relative Hopf modules 
and modules of cross product algebras and we characterize 
direct sums of Heisenberg 
representations for  cross product algebras 
$\cA \rti \cU$ of Hopf $\ast$-al\-ge\-bras 
$\cA$ of compact quantum groups. Section \ref{S2} 
collects a number of definitions and 
basic facts on Podle\'s spheres and on the cross product algebras 
$\cO(\mathrm{S}^2_{qr}) \rti \cU_q(\mathrm{su}_2)$ 
which are needed in what follows.

All facts and notions on quantum groups used in this paper can be found, 
for instance, in \cite{KS}. The algebra $\cU_q(\mathrm{sl}_2)$ is due to
P.~P. Kulish and N.~Y. Reshetikhin \cite{KR}. The quantum group 
$\mathrm{SU}_q(2)$ was
invented in \cite{W} and \cite{VS} and the quantum spheres $\mathrm{S}^2_{qr}$
were discovered in \cite{P}.

Let us introduce some notation. 
Throughout this paper, $q$ is a real number of the open interval 
$(0,1)$. We abbreviate 
$$\lambda:=q\!-\!q^{-1},\quad
\lambda_n:=(1\!-\!q^{2n})^{1/2},\quad 
[n]:= (q^n \!-\! q^{-n})/(q \!-\! q^{-1}),
$$ 
where $n\in\dN_0$. 
For an algebra $\cX$, the notation 
$\cX^n$ stands for the 
direct sum of $n$ copies of $\cX$ 
and $\mathrm{M}_{n} (\cX)$ denotes the set of 
$n\times n$-matrices with entries 
from $\cX$. 
Let $I$ be an at most countable index set, $V$ a linear space and
$\cD=\mathop{\oplus}_{i\in I}V_i$, where $V_i=V$ for all 
$i\in I$. We denote by $\eta_i$ the vector of $\cD$ which has the 
element $\eta\in V$ as its $i$-th component and zero otherwise.
It is understood that $\eta_i=0$ whenever $i\notin I$. 
If $V$ is a Hilbert space, then $\mathop{\bar \oplus}_{i\in I}V_i$ 
denotes the closed linear span of $\{\eta_i\,;\,\eta\in V,\; i\in I\}$. 
Given a dense linear subspace $\cD$ of a Hilbert space, 
$$
   \cL^+(\cD):=\{\,y\in {\rm End}(\cD)\,;\,\cD\subset \cD(y^*),\ 
y^*\cD\subset \cD\,\} 
$$
is a unital *-al\-ge\-bra of closeable operators with involution 
$y\mapsto y^*\lceil \cD$. 
As customary, $\cD(x)$ denotes the domain of an 
operator $x$.  
By a $\ast$-re\-pre\-sen\-ta\-tion of a unital 
$\ast$-al\-ge\-bra $\cX$ on the domain $\cD$, we mean a unit preserving 
$\ast$-homomorphism $\pi$ from  $\cX$ into $\cL^+(\cD)$
(see e.g.\ \cite{S1}). 
When no confusion can arise, we omit the letter which denotes 
the representation and write $x$ instead of $\pi(x)$. 
%
\section{Relative Hopf modules and representations of the 
cross product algebras}                                         \label{S1}
%
Let $\cU$ be a Hopf $\ast$-al\-ge\-bra and let $\cX$ be a left $\cU$-mo\-dule 
$\ast$-al\-ge\-bra, that is,  $\cX$  is a unital $\ast$-al\-ge\-bra with left 
$\cU$-action $\ang$ satisfying 
\begin{equation}                                                 \label{E1}
f\ang xy = (f_{(1)}\ang x) (f_{(2)}\ang y),\quad 
f\ang 1  = \varepsilon(f) 1,\quad 
(f\ang x)^\ast = S(f)^\ast \ang x^\ast 
\end{equation}
for $x$, $y \in \cX$ and $f \in \cU$. 
Here $\Delta (f)=f_{(1)} \otimes f_{(2)}$ is the Sweedler notation for 
the comultiplication $\Delta(f)$ of $f\in \cU$. 
Then the {\it left cross product $\ast$-al\-ge\-bra} $\cX\rti\cU$ is 
the $\ast$-al\-ge\-bra generated by the two $\ast$-sub\-al\-ge\-bras $\cX$ 
and $\cU$ with respect to the cross commutation relations 
\begin{equation}\label{fxrel}
fx = (f_{(1)}\ang x) f_{(2)}, \quad x\in\cX,\ \,f\in\cU.
\end{equation}
Let $\cA$ be a Hopf algebra and $\cX$ a right $\cA$-coideal subalgebra 
of $\cA$. A relative Hopf module in ${_\cX}\cM^\cA$ (see e.g. \cite{M})
is a right $\cA$-co\-mo\-dule $M$ which is also a left $\cX$-mo\-dule 
such that 
the right coaction of $\cA$ is $\cX$-linear, that is,
\begin{equation}\label{linear}
(xm)_{(1)}\otimes (xm)_{(2)}=x_{(1)} m_{(1)} \otimes x_{(2)} m_{(2)},\quad
x\in \cX,\ \, m\in M.
\end{equation}
Here we write simply $xm$ for the left module action of $x$ at $m$ 
and we use the Sweedler notation for the coaction.

Let $\langle\cdot,\cdot\rangle$ be a dual pairing of Hopf 
algebras $\cU$ and $\cA$. Then any right $\cA$-co\-mo\-dule $M$ 
determines a left $\cU$-mo\-dule by
\begin{equation}\label{lmod1}
f\ang m := \langle f, m_{(2)}\rangle m_{(1)},\quad f\in \cU,\ \, m\in M,
\end{equation}
and the right $\cA$-co\-mo\-dule algebra $\cX$ 
becomes a left $\cU$-mo\-dule algebra. 
Hence the cross product algebra $\cX\rti\cU$ is defined. The 
cross relations (\ref{fxrel}) are given by 
\begin{equation}                                     \label{crossrel}
 fx=\langle f_{(1)}, x_{(2)}\rangle x_{(1)}f_{(2)},\quad 
x\in \cX,\ \, f\in\cU.
\end{equation}

The following simple well-known lemma is crucial in what follows.

\begin{thl}                                   \label{Lem1}                  
Let $M$ be a left $\cX$-mo\-dule and a right $\cA$-co\-mo\-dule.
\begin{enumerate}
\item[(i)] If $M\in{_\cX}\cM^\cA$, then $M$ is a left $\cX\rti\cU$-mo\-dule.
\item[(ii)] Suppose that $\cU$ separates the points of $\cA$, that is, 
$\langle f,a\rangle=0$ for all $f\in\cU$ implies $a=0$. If $M$ is a 
left $\cX\rti\cU$-mo\-dule, then $M\in{_\cX}\cM^\cA$.
\end{enumerate}
\end{thl}
{\bf Proof.} As $M$ is a right $\cA$-co\-mo\-dule, it is a left $\cU$-mo\-dule.
From (\ref{lmod1}), we obtain
\begin{align}\label{lmod2}
&f\ang (xm)=\langle f, (xm)_{(2)}\rangle (xm)_{(1)},\\
\label{lmod3}
&\langle f_{(1)}, x_{(2)}\rangle x_{(1)} (f_{(2)}\ang m)= 
\langle f, x_{(2)} m_{(2)}\rangle x_{(1)} m_{(1)}.
\end{align}
If $M\in {_\cX}\cM^\cA$, then (\ref{linear}) holds and hence the right hand 
sides of (\ref{lmod2}) and (\ref{lmod3}) coincide. That is, 
the cross relations (\ref{crossrel}) of $\cX\rti\cU$ are 
satisfied, so we have a 
well defined left $\cX\rti\cU$-mo\-dule. Conversely,  if the right hand 
sides of (\ref{lmod2}) and (\ref{lmod3}) are equal,  then (\ref{linear}) 
follows by using the assumption that $\cU$ separates $\cA$.\hfill$\Box$
\mn

Suppose now that $\langle\cdot,\cdot\rangle$ is a dual pairing of 
Hopf $\ast$-al\-ge\-bras $\cU$ and $\cA$ and that $\cX$ 
is a $\ast$-in\-var\-iant 
right $\cA$-coideal. Suppose that there exists a positive linear 
functional $h$ on $\cX$ (i.e.\  $h(x^\ast x)\ge 0$ for $x\in\cX$) which 
is\ \,$\cU$-in\-var\-iant 
(i.e.\ $h(f\ang x)=\varepsilon (f) h(x)$ for $x\in\cX$ 
and $f\in\cU$). Then there is a unique $\ast$-re\-pre\-sen\-ta\-tion $\pi_h$, 
called the {\it Heisenberg representation} 
of the $\ast$-al\-ge\-bra $\cX\rti\cU$, with domain $\cD_h$ 
such that $\pi_h\lceil\cX$ is the 
GNS-re\-pre\-sen\-ta\-tion of $\cA$ with cyclic vector 
$\varphi_h=\pi_h(1)$, $\cD_h=\pi(\cX)\varphi_h$  and 
$\pi(f)\varphi_h=\varepsilon (f)\varphi_h$ for $f\in\cU$. 
(Note that we have taken  in \cite{SW} 
the closure of $\pi_h$ as the Heisenberg representation obtained from $h$.)

Recall that a Hopf $\ast$-al\-ge\-bra $\cA$ is called a {\it CQG-al\-ge\-bra} 
(compact quantum group algebra) if $\cA$ is the linear span of matrix 
elements of finite dimensional unitary corepresentations of 
$\cA$ \cite{DK}, \cite[Subsection 11.3.1]{KS}. 
A  CQG-al\-ge\-bra has a unique Haar 
state $h$.

The next proposition is a Hilbert space version of the well-known algebraic 
fact that Hopf modules in $_\cA{\cM^\cA}$ are trivial.
\begin{thp}                                        \label{1.2}
Suppose that $\langle\cdot,\cdot\rangle$ is a dual pairing of a  Hopf 
$\ast$-al\-ge\-bra $\cU$ and a CQG Hopf $\ast$-al\-ge\-bra $\cA$ such that 
$\cU$ separates the points of $\cA$. Let $\pi$ be a 
$\ast$-re\-pre\-sen\-ta\-tion 
of the cross product $\ast$-al\-ge\-bra $\cA\rti\cU$ on a domain $\cD$. 
Then $\pi$ is unitarily equivalent to a direct sum of Heisenberg 
representations $\pi_h$ of $\cA\rti\cU$
if and only if $\cD$ is a right $\cA$-co\-mo\-dule such that 
$\pi(f)\varphi = f\ang \varphi$ 
for $f\in \cU$ and $\varphi \in\cD$.
\end{thp}
{\bf Proof.} 
To prove the necessity, it suffices to check the 
above condition for the Heisenberg 
representation $\pi_h$. Because the Haar state of $\cA$ is faithful, 
there is a 
well defined right coaction $\phi$ of $\cA$ on 
$\cD$ given by 
$\phi (\varphi) = \pi_h (a_{(1)}) \varphi_h \otimes a_{(2)}$  for  
$\varphi = \pi_h (a) \varphi_h$, $a\in\cA$. Let $f\in \cU$. 
Using the cross relation 
(\ref{fxrel}) and $\cU$-invariance of the cyclic vector $\varphi_h$,  
we obtain
\begin{align*}
&\pi_h (f) \varphi = \pi_h (fa) \varphi_h 
= \langle f_{(1)}, a_{(2)}\rangle \pi_h  
(a_{(1)} f_{(2)}) \varphi_h
=\langle f, a_{(2)}\rangle \pi_h (a_{(1)}) \varphi_h = f\ang \varphi.
\end{align*}

Now we prove the sufficiency. Suppose $\cD$ is a right $\cA$-co\-mo\-dule 
such that 
$\pi (f) \varphi = f\ang \varphi \equiv \langle f, 
\varphi_{(2)}\rangle \varphi_{(1)}$ 
for $\varphi \in \cD$ and $f\in \cU$. Let $\cD_{\inv}$ 
be the subspace of vectors 
$\omega(\varphi) :=\pi(S^{-1}(\varphi_{(2)}))\varphi_{(1)}$,  
$\varphi\in\cD$. Let $f\in\cU$ and $\varphi\in\cD$. Using the fact that 
$\pi$ is a representation and relation (\ref{fxrel}), we compute
\begin{align*}                                                 \label{inv}
\pi(f)\omega(\varphi)&=\pi(fS^{-1}(\varphi_{(2)}))\varphi_{(1)}
=\langle f_{(1)}, S^{-1} (\varphi_{(2)})\rangle 
\pi(S^{-1} (\varphi_{(3)}))\pi(f_{(2)}) \varphi_{(1)}\\
&=\langle f_{(1)}, S^{-1} (\varphi_{(3)})
\rangle \pi(S^{-1} (\varphi_{(4)}))\langle f_{(2)}, 
\varphi_{(2)}\rangle \varphi_{(1)}\\
&=\langle f, S^{-1} (\varphi_{(3)})\varphi_{(2)}\rangle 
\pi(S^{-1} (\varphi_{(4)})) \varphi_{(1)}=\varepsilon(f)\omega( \varphi).
\end{align*}
Hence the functional $\langle \pi(\cdot) \omega(\varphi),
\omega(\varphi)\rangle$ on $\cA$ is $\cU$-in\-var\-iant. 
Since $\cU$ separates the points of $\cA$, it is also $\cA$-in\-var\-iant 
and so a multiple of the Haar state $h$. Thus, we have
\begin{equation}                                      \label{haar}
\langle \pi(x)\omega(\varphi), 
\omega(\varphi)\rangle = \|\omega(\varphi)\|^2 h(x),\quad  x\in\cA.
\end{equation}
Next we show that $\omega(\varphi)\bot \omega(\varphi^\prime)$ 
for $\varphi,\varphi^\prime\in\cD$ implies 
\begin{equation}\label{ortho}
\pi(\cA\rti\cU)\omega(\varphi)\,\bot\, \pi(\cA\rti\cU)\omega(\varphi^\prime).
\end{equation}
Let $u\in\dC$, $|u|=1$, and let $x\in \cA$ be hermitian.  
Clearly, $\omega(\varphi)+u\omega(\varphi^\prime)\in\cD_{\inv}$. 
By (\ref{haar}), we obtain
\begin{align*}
&(\|\omega(\varphi)\|^2 + \|\omega(\varphi^\prime)\|^2 ) h(x)
=\langle \pi(x) (\omega(\varphi) + u \omega (\varphi^\prime)), 
\omega (\varphi) + u \omega (\varphi^\prime)\rangle\\
&\qquad=\langle \pi(x) \omega(\varphi),\omega (\varphi)\rangle + \langle  
\pi(x) \omega (\varphi^\prime),
\omega(\varphi^\prime) \rangle 
+ 2\mathrm{Re}\,u\, \langle \pi(x) \omega (\varphi), 
\omega(\varphi^\prime)\rangle\\
&\qquad=(\|\omega(\varphi)\|^2 
+ \| \omega(\varphi^\prime)\|^2 ) h(x) + 2 \mathrm{Re}\,u\, 
\langle \pi (x) \omega (\varphi), \omega (\varphi^\prime)\rangle.
\end{align*}
Therefore, $\mathrm{Re}\,u\,\langle \pi(x) \omega (\varphi), 
\omega (\varphi^\prime)\rangle=0$ for all $u\in\dC$ with $|u|=1$
which implies 
$\langle \pi(x)\omega(\varphi), \omega(\varphi^\prime)\rangle=0$. 
Since $\cA$ is spanned by hermitian elements and 
$\pi$ is a $\ast$-re\-pre\-sen\-ta\-tion, 
$\pi(\cA)\omega(\varphi)\bot \pi(\cA)\omega(\varphi^\prime)$. 
From (\ref{fxrel}), we obtain 
$\pi(f) \pi(\cA) \omega (\psi)\subseteq \pi(\cA)\omega(\psi)$ 
for all $\psi\in\cD$, so (\ref{ortho}) follows.

Using Zorn's lemma, we choose a maximal set 
$\{\omega(\varphi_i)\,;\,i \in I\}$ of orthonormal vectors of the vector 
space $\cD_{\inv}$. Let $\cD_i = \pi (\cA\rti\cU) \omega (\varphi_i)$ 
and let $\pi_i$ be the restriction of the representation $\pi$ to $\cD_i$. 
Since $\langle \pi_i(x)\omega(\varphi_i), \omega(\varphi_i)\rangle= h(x)$ 
for $x\in \cA$ by (\ref{haar}), $\pi_i$ is unitarily equivalent to the 
Heisenberg representation $\pi$ of $\cA\rti\cU$. By (\ref{ortho}), 
$\cD_i\bot\cD_j$ for $i\ne j$. Since 
$\varphi=\pi(\varphi_{(2)})\omega(\varphi_{(1)})$ for $\varphi\in\cD$, 
the domain $\cD$ is the linear span of subspaces $\cD_i$, $i\in I$. 
Putting the preceding together, we have shown that $\pi$ is 
unitarily equivalent to a direct sum of Heisenberg representations.
                                                       \hfill $\Box$
\mn 

Let $\cA$ be the $CQG$-al\-ge\-bra $\cO(\mathrm{SU}_q(2))$ and let  
$\cU = \cU_q(\mathrm{su}_2)$. 
Then the hypothesis in Proposition \ref{1.2} is satisfied 
if and only if $\pi$ is integrable, that is, its restriction 
to $\cU_q(\mathrm{su}_2)$ is a direct sum of spin $l$ 
representations $T_{l}$, $l\in \frac{1}{2}\dN_0$. 
Therefore, by Proposition \ref{1.2}, 
a $\ast$-re\-pre\-sen\-ta\-tion $\pi$ of the $\ast$-al\-ge\-bra 
$\cO(\mathrm{SU}_q(2))\rti \cU_q(\mathrm{su}_2)$ is 
a direct sum of Heisenberg representations $
\pi_h$ if and only $\pi$ is integrable. 
This was proved in \cite{SW} by another method. 
A similar result holds for compact 
forms of standard quantum groups.

%
\section{The cross product algebra 
{\mathversion{bold}$\cO({\mathrm{S}}^2_{qr})\rti \cU_q({\mathrm{su}}_2)$} }  
                                                                  \label{S2}
%
The Hopf $\ast$-al\-ge\-bra $\cU_q(\mathrm{su}_2)$ is generated by elements 
$E,F,K,K^{-1}$ with relations
\begin{equation}\label{urel}
KK^{-1} \!=\! K^{-1} K\!=\!1,\ KE\!=\!qEK,\, FK\!=\!qKF,\
EF- FE\!=\!\lambda^{-1}(K^2-K^{-2}),
\end{equation}
involution $E^\ast=F$, $K^\ast=K$, comultiplication 
$$
\Delta (E)\!=\!E\otimes K+K^{-1} \otimes E,\ \ 
\Delta (F)\!=\!F\otimes K + K^{-1}\otimes F,\ \  \Delta (K)\!=\!K\otimes K,
$$
counit $\varepsilon(E)=\varepsilon(F)=\varepsilon(1-K)=0$ 
and antipode $S(K)=K^{-1}$, $S(E)=-qE$, $S(F)=-q^{-1}F$.
There is a dual pairing $\langle \cdot,\cdot\rangle$ of the Hopf 
$\ast$-al\-ge\-bras $\su$ and $\SU$ given on generators by
$$
\langle K^{\pm 1},d\rangle = \langle K^{\mp 1},a\rangle = q^{\pm 1/2}, \quad
\langle E,c\rangle = \langle F,b\rangle=1
$$
and zero otherwise, where $a,b,c,d$ are the usual generators 
of $\SU$ (see e.g.\ \cite[Chapter 4]{KS}). 

We shall use the definition of the coordinate algebras 
$\cO(\mathrm{S}^2_{qr})$, $r\in [0,\infty]$, 
of Podle\'s spheres as given in \cite{P}. 
For $r\in[0,\infty)$, $\cO(\mathrm{S}^2_{qr})$ is the $\ast$-al\-ge\-bra 
with generators $A\!=\!A^\ast, B,B^\ast$ and defining relations
\begin{equation}\label{podrel}
AB\!=\!q^{-2} BA,\ AB^\ast\!=\!q^2 B^\ast A,\ B^\ast B\!=\!A-A^2+r,\ 
BB^\ast \!=\!q^2 A-q^4 A^2+r.
\end{equation}
For $r=\infty$, the defining relations of 
$\cO(\mathrm{S}^2_{q ,\infty})$ are
\begin{equation}\label{prodrel1}
AB=q^{-2} BA,\  AB^\ast = q^2B^\ast A,\ B^\ast B\!=\! -A^2+1,\ 
BB^\ast\!=\!-q^4 A^2+1.
\end{equation}
Let $r<\infty$. 
As shown in \cite{P}, 
$\cO(\mathrm{S}^2_{qr})$ is a right $\SU$-co\-mo\-dule 
$\ast$-al\-ge\-bra such that 
\begin{equation}\label{xgener}
x_{-1}\!:=\! q^{-1}(1+q^{2})^{1/2}B,\ \, x_1 :=-(1+q^{2})^{1/2}B^\ast,\ \,
x_0 := 1-(1+q^2)A.
\end{equation}
transform by the spin $1$ matrix corepresentation 
$(t^1_{ij})$ of $\mathrm{SU}_q(2)$. Hence $\cO(\mathrm{S}^2_{qr})$ 
is a left $\su$-mo\-dule $\ast$-al\-ge\-bra with left action given by 
$f\ang x_j=\sum_i x_i\langle f,t^1_{ij}\rangle$ 
for $f\in \su$, $j=-1,0,1$. Inserting the form of the 
matrix $(t^1_{ij})$ (see \cite{P} or \cite[Subsection 4.5.1]{KS}) and 
the Hopf algebra pairing $\langle\cdot,\cdot\rangle$ into (\ref{crossrel}),   
we derive the following cross relations for the cross product algebra
$\cO(\mathrm{S}^2_{qr})\rti \su$:
\begin{align*}
&KA \!=\! AK,\ \ EA \!=\! AE + q^{-1/2} B^\ast K,\ \ 
FA \!=\!AF - q^{-3/2} BK,\\
&KB \!=\! q^{-1} BK,\ \ EB \!=\! q BE -q^{1/2} (1+q^2) AK + q^{1/2} K, \  \ 
FB \!=\! q BF,\\
&KB^\ast \!=\! q B^\ast K,\, EB^\ast \!=\! q^{-1} B^\ast E, \,
FB^\ast \!=\! q^{-1} B^\ast F + q^{-1/2} (1+q^2) AK - q^{-1/2} K.
\end{align*}
That is, $\cO(\mathrm{S}^2_{qr})\rti \cU_q(\mathrm{su}_2)$ is the 
$\ast$-al\-ge\-bra 
with generators $A$, $B$, $B^\ast$, $E$, $F$, $K$, $K^{-1}$, with defining 
relations (\ref{urel}), (\ref{podrel}) and 
the preceding set of cross relations.

For $r=\infty$, we set 
\begin{equation}\label{xgener1}
x_{-1}:= q^{-1}(1+q^{2})^{1/2}B,\ \,x_1 :=-(1+q^{2})^{1/2}B^\ast,\ \,
x_0 := -(1+q^2)A.
\end{equation}
Then the cross relations for 
$\cO(\mathrm{S}^2_{q\infty})\rti \su$ can be written as 
\begin{align*}
&KA = AK,\ \ EA = AE + q^{-1/2} B^\ast K,\ \ 
FA =AF - q^{-3/2} BK,\\
&KB = q^{-1} BK,\ \  EB = q BE -q^{1/2} (1+q^2) AK, \ \   
FB = q BF,\\
&KB^\ast = q B^\ast K,\ \ EB^\ast = q^{-1} B^\ast E, \  \ 
FB^\ast = q^{-1} B^\ast F + q^{-1/2} (1+q^2) AK.
\end{align*}

There is an infinitesimal description of Podle\'s quantum spheres which was
discovered in  \cite{DK}. 
For $r\in[0,\infty]$, define an element $X_r\in\cU_q(\mathrm{su}_2)$ by
\begin{align}
\begin{split}                                            \label{Xr}
  X_r &= q^{1/2}(q^{-1}-q)^{-1}r^{-1/2}(1-K^2)+EK+qFK,\quad r\in(0,\infty),\\
  X_0 &= 1-K^2,\quad   X_\infty =EK + qFK.
\end{split}
\end{align}
Then $\Delta(X_r)=1\otimes X_r+X_r\otimes K^2$. 
As shown in \cite{DK}, the right coideal $\ast$-sub\-al\-ge\-bra
$$
\cX := \{x\in\cO(\mathrm{SU}_q(2))\,;\,\langle X_r, x_{(1)}\rangle x_{(2)}=0\}
$$
of infinitesimal invariants with respect to $X_r$ can be identified with 
the coordinate $\ast$-al\-ge\-bra $\cO(\mathrm{S}^2_{qr})$ of 
Podle\'s quantum 2-spheres.
The generators $x_{-1}$, $x_0$, $x_1$ of $\cO(\mathrm{S}^2_{qr})$ are then 
identified with the elements 
\begin{align*}
x_{-1}&=(1+q^{-2})^{1/2}(r^{1/2}a^2+ac-qr^{1/2}c^2),\\
x_0&=(1+q^{-2})r^{1/2}ab+1+(q+q^{-1})bc-(1+q^{2})r^{1/2}dc,\\
x_1&=(1+q^{-2})^{1/2}(r^{1/2}b^2+bd-qr^{1/2}d^2)
\quad \mbox{for}\ \,r<\infty;\\[6pt]
x_{-1}&=(1+q^{-2})^{1/2}(a^2-qc^2),   \\
x_0&=(1+q^{-2})(ab-q^2dc),\\
x_1&=(1+q^{-2})^{1/2}(b^2-qd^2)
\quad \mbox{for}\,\ r=\infty.
\end{align*}
%
%
\section{Decomposition of {\mathversion{bold}$\cO({\mathrm{SU}}_q(2))$} } 
                                                   \label{sec-decomp}
Throughout this section, we denote by $\cX$ the coordinate algebras 
$\cO(\mathrm{S}^2_{qr})$ and by $\cA$ and $\cU$ the  Hopf $\ast$-al\-ge\-bras  
$\cO (\mathrm{SU}_q(2))$ and $\cU_q(\mathrm{su}_2)$, 
respectively. 

First we recall a few crucial algebraic results from \cite{MS}
needed in what follows.
Let $\cC$ denote the coalgebra $\cA/\cX^+\cA$ with quotient map  
$\rho:\cA\rightarrow\cC$, where $\cX^+=\{x\in\cX\,;\,\varepsilon (x)=0\}$. 
By \cite{MS}, $\cC$ is spanned by group-like elements and, as a consequence,  
it is the direct sum of simple 
subcoalgebras $\cC_j$, $j\in\frac{1}{2}\dZ$. 
Set 
$$
M_j=\{x\in \cA\,;\, \rho(x_{(1)}) \otimes x_{(2)}\in \cC_j\otimes\cA\}.
$$
Then $M_j$ is a finitely generated projective Hopf module in ${_\cX}\cM^\cA$
and $\cA$ is the direct sum of all $M_j$ 
(cf.\ \cite[p.\ 163]{MS}). 
Since the dual pairing of $\cU$ and $\cA$ is non-degenerate 
\cite[Section 2.4]{KS}, it follows from Lemma \ref{Lem1} that 
Hopf modules in ${_\cX}\cM^\cA$
are in one-to-one correspondence with left $\cX\rti\cU$-mo\-dules.

In this section, we reconsider these algebraic results in the 
Hilbert space setting. 
Since the Haar state on $\cA$ is faithful, we can equip $\cA$ with 
the inner product 
\begin{equation}                                            \label{ip}
\langle a,b\rangle := h(b^\ast a),\quad a,b\in\cA.
\end{equation}
The Heisenberg representation of $\cA$ is just the left $\cA \rti
\cU$-mo\-dule $\cA$ endowed with the inner product (\ref{ip}). 
Hence each left 
$\cX\rti\cU$-mo\-dule $M_j$ corresponds to a $\ast$-re\-pre\-sen\-ta\-tion,
denoted by $\hat{\pi}_j$, of the cross product
$\ast$-al\-ge\-bra $\cO(\mathrm{S}^2_{qr})\rti\cU_q(\mathrm{su}_2)$. We
will return to this representation in Proposition \ref{cor} below.

In Theorem \ref{T1}, we give a direct proof 
that $\cA=\cO(\mathrm{SU}_q(2))$ is the orthogonal 
direct sum of Hopf modules $M_j$ and describe 
this decomposition explicitely. Moreover, the projections 
of the projective $\cX$-mo\-dules $M_j$ are computed. 

Let us first introduce some more notation. 
We abbreviate 
\begin{align*}
&\lambda_\pm := 1/2 \pm (r+ 1/4)^{1/2} \ \,\mbox{for}\ \,r<\infty,\quad 
\lambda_\pm := \pm1 \ \,\mbox{for}\ \,r=\infty, \\
&s:=0\ \,  \mbox{for}\ \,r=0,\quad 
   s:= -r^{-1/2} \lambda_-  
\ \,  \mbox{for}\ \,r\in (0,\infty), \quad   
   s:= 1 \ \, \mbox{for}\ \,r=\infty.
\end{align*}
Note that $s= r^{1/2} \lambda_+^{-1}$ when $r\in [0,\infty)$. 
For $j\in\frac{1}{2}\dN$, define
\begin{align}                                         
u_j &= (d + q^{-1}s b)(d+q^{-2} s b)\dots (d+q^{-2j}s b),\label{uj}\\
w_j &= (a - q sc)(a-q^2 s c) \dots (a-q^{2j} s c), \label{wj}\\        
u_{-j} &= E^{2j}\ang w_j,  \label{uj-}
\end{align}
and set $u_0=w_0=1$. From \cite{BM}, 
we conclude that $\rho_L(u_j)$, $\rho_L(w_j)$, $j\in\frac{1}{2}\dN_0$, 
are group like elements that span the coalgebra $\cA/\cX^+\cA$, where 
$\rho_L :\cA\rightarrow \cA/\cX^+ \cA$ is the canonical mapping. 
(In order to apply the results in \cite{BM}, 
one has to interchange the right $\cA$-co\-mo\-dule  
algebra $\cX$ with the left $\cA$-co\-mo\-dule 
algebra $\theta(\cX)$ using the 
$\ast$-al\-ge\-bra automorphism and coalgebra 
anti-homomorphism $\theta:\cA\rightarrow\cA$ determined by
$\theta (a) = a$, $\theta (d)=d$, $\theta(b) =-qc$,  
$\theta(c)=-q^{-1} b$.) 
In particular, the simple coalgebras $\cC_j$ are given 
by $\cC_j=\dC\hspace{1pt}\rho_L(u_j)$ 
and $\cC_{-j}=\dC\hspace{1pt}\rho_L(w_j)$, 
$j\in\frac{1}{2}\dN_0$. As a consequence, 
$u_j\in M_j$ for $j\in\frac{1}{2}\dZ$.
A crucial role 
will play the elements  
\begin{equation}                                            \label{vlkj}
v^l_{kj} := N^l_{kj}\, F^{l-k} \ang (x^{l-|j|}_1 u_j), 
\quad l\in\mbox{$\frac{1}{2}\dN_0$}, \ \, 
j,k=-l,-l+1,{\dots}\,,l,
\end{equation}
of $\cA$, where $N^l_{kj} = \| F^{l-k}\ang (x^{l-|j|}_1 u_j)\|^{-1}$.
To describe the projection of the projective module $M_j$, 
we define a $(2|j|{+}1){\times}(2|j|{+}1)$-matrix $P_j$ 
with entries from $\cA$ by
$$
P_j = \Big(\,q^{-(n+m)} [2|j|\!+\!1]^{-1} v^{|j|}_{nj} v^{|j|\ast}_{mj}\,
\Big)^{|j|}_{n,m=-|j|}.
$$
\begin{tht}                                                 \label{T1}
The decomposition $\cA=\oplus_{j\in\frac{1}{2}\dZ} M_j$ is an orthogonal 
direct sum with respect to the inner product given by (\ref{ip}).
The set 
$\{v^l_{kj}\,;\, l\in\frac{1}{2}\dN_0,\ k,\,j=-l$, $-l+1,{\dots}\,,l\}$ 
is an orthonormal basis of $\cA$ and 
$$
M_j = \Lin \{v^{|j|+n}_{kj}\,;\, n\in\dN_0,\ 
k=-(|j|+n),-(|j|+n)+1,{\dots}\,,|j|+n\}.
$$
As a left $\cX$-mo\-dule, $M_j$ is generated by 
$\{v^{|j|}_{kj}\,;\, k=-|j|,-|j|+1,{\dots}\,, |j|\}$, the matrices 
$P_j$ are orthogonal projections in $\mathrm{M}_{2|j|+1} (\cX)$, 
and  $M_j$ is isomorphic to $\cX^{2|j|+1} P_j$.
Each vector $v^{|j|}_{kj}$ is cyclic for the 
$\ast$-re\-pre\-sen\-ta\-tion $\hat\pi_j$ 
of $\cX\rti\cU$ on $M_j$.
\end{tht}
{\bf Proof.} Obviously, $v^l_{lj}$ is a highest weight 
vector of weight $l$ and the linear space 
$
V^l_j \hsp:=\hsp\Lin  \{v^l_{kj}\,;\, k \hsp=\hsp -l,-l+1,{\dots}\,,l\}
$
is an irreducible $\cU_q(\mathrm{su}_2)$-mo\-dule of spin $l$. 
Hence $v^l_{kj}$ and 
$v^{l^\prime}_{k^\prime j^\prime}$ are orthogonal whenever $l\ne l^\prime$ 
or  $k\ne k^\prime$ because then they belong to representations of 
different spin or are vectors of different weights. 

It remains to prove that $v^l_{kj}$ and $v^l_{kj^\prime}$ are orthogonal 
when $j\ne j^\prime$. The idea of the proof \cite{MS} 
is to show that $v^l_{kj}$ 
and $v^l_{kj^\prime}$ are eigenvectors of different eigenvalues of a 
hermitian operator $\hat X_r$ acting on $\cA$.
Let $\hat{X}_r$ be defined by
$$
\hat{X_r} (a) := a\anf X_r=\langle X_r,a_{(1)}\rangle a_{(2)},\quad  a\in\cA, 
$$
where $X_r\in\cU_q(\mathrm{su}_2)$ is given by (\ref{Xr}).
The relation $S(X_r)^\ast = S(X_r)$ implies that the operator  $\hat{X}_r$ is
hermitian. Indeed, for $a,b \in \cA$, we have 
\begin{align*}
\langle a,b\anf X_r\rangle &= h ((b\anf X_r)^\ast a) 
= h((b^\ast \anf S(X_r)^\ast)a)=h((b^\ast \anf S(X_r)) a) \\
&= \varepsilon (X_{r(2)}) h((b^\ast \anf S(X_{r(1)}))a)
=h((b^\ast \anf S(X_{r(1)}) X_{r(2)}) (a\anf X_{r(3)})) \\
&= h(b^\ast (a\anf X_r))
=\langle a\anf X_r, b\rangle.
\end{align*}
For $j\in \frac{1}{2}\dZ$, set
\begin{align*}
&\mu_j := 1-q^{2j}\quad \mbox{for}~~ r=0,\\
&\mu_j := q^{1/2} (q^{-1} -q)^{-1} r^{-1/2} (1-q^{-2j} \lambda_- 
- q^{2j} \lambda_+)\quad \mbox{for}~~r\in (0,\infty),\\
&\mu_j := q^{1/2} (q^{-1} -q)^{-1} (q^{-2j} - q^{2j})\quad 
\mbox{for}~~r=\infty.
\end{align*}
We claim that $u_j \anf X_r = \mu_j u_j$, where the vectors 
$u_j$, $j\in \frac{1}{2}\dZ$, 
are defined by (\ref{uj})--(\ref{uj-}).
Let $r\in (0,\infty)$. For $j=0$, we have 
$u_0 \anf X_r = \varepsilon (X_r) = 0$.
Assume that the assertion holds for $j\in \frac{1}{2}\dN_0$.
Using $\Delta (X_r)= 1 \otimes X_r + X_r \otimes K^2$ and 
$s=-r^{-1/2}\lambda_- = r^{-1/2} \lambda_+^{-1}$, we compute
\begin{align*}
&u_{j+1/2} \anf X_r 
= u_j \big((d+q^{-(2j+1)} s b) \anf X_r \big)
+ u_j \anf X_r \big((d+q^{-(2j+1)}s b)\anf K^2\big)\\
&\quad= u_j \big[ ( q^{1/2} r^{-1/2} 
(q^{-1} - q)^{-1} (1-q) - q^{1/2} r^{-1/2} q^{-2j} \lambda_- + q \mu_j) d\\
&\quad\quad +(q^{1/2} r^{-1/2} (q^{-1}- q)^{-1} (1-q^{-1}) 
+ q^{1/2} r^{-1/2} q^{2j} \lambda_+ 
+ q^{-1} \mu_j)q^{-(2j+1)}s b\big]\\
&\quad=\mu_{j+1/2} u_j (d+q^{-(2j+1)} s b) 
= \mu_{j+1/2} u_{j+1/2}.
\end{align*}
By induction, the claim follows for $j\in\frac{1}{2}\dN_0$.
Similarly, one proves that 
$w_j \anf X_r = \mu_{-j} w_j$ for $j\in \frac{1}{2}\dN$. 
Since $(E^{2j} \ang w_j)\anf X_r = E^{2j} \ang (w_j \anf X_r)$, 
we obtain $u_j \anf X_r = \mu_j u_j$ for all $j\in\frac{1}{2}\dZ$. 
Analogous, but simpler, computations show that the claim also 
holds for $r=0$ and $r=\infty$.
Using $x_1\anf X_r = 0$, we obtain 
\begin{align*}
\big(N^l_{kj} F^{l-k} \ang (x^{l-|j|}_1 u_j)\big) \anf X_r 
&= N^l_{kj} F^{l-k} \ang \big(x^{l-|j|}_1 (u_j \anf X_r) 
+ (x^{l-|j|}_1 \anf X_r)(u_j\anf K^2)\big)\\  
&= \mu_j N^l_{kj} F^{l-k} \ang (x^{l-|j|}_1 u_j),
\end{align*}
so that $\hat{X_r} (v^l_{kj})= \mu_j v^l_{kj}$. 
Since $\hat{X_r}$ is hermitian and $\mu_{j}\neq  \mu_{j^\prime}$, 
we conclude that 
$v^l_{kj}$ and $v^l_{kj^\prime}$ are orthogonal 
whenever $j\ne j^\prime$.

The decomposition of $\cA$ into an orthogonal direct sum is 
a consequence of above results. 
Let $t^l_{kj}\in\cA$, $l\in\frac{1}{2}\dN_0$, $k,j=-l,-l+1,{\dots}\,,l$, 
denote the matrix elements from the Peter-Weyl decomposition of $\cA$ 
\cite[Section 4.2]{KS}. As $v^l_{kj}$ is a weight 
vector with weight $k$ of a spin $l$ representation of 
$\cU_q(\mathrm{su}_2)$, we know that 
$v^l_{kj}\! \in\!\Lin \{t^l_{ki}\,;\, {i=-l,-l+1,{\dots}\,, l}\}$. 
A simple dimension argument shows that 
$\Lin\{v^l_{ki}\,; \,i=-l,-l+1,\dots\,,l\}
=\Lin \{t^l_{ki}\,;$ $i=-l,-l+1,{\dots}\,,l\}$. 
Since the elements $t^l_{kj}$ 
span $\cA$, we conclude that 
$\{v^l_{kj}\,;\, l\in\frac{1}{2}\dN_0,\ k,j=-l,-l+1,{\dots}\,, l\}$ 
is an orthonormal basis of $\cA$.
Recall that $v^{|j|}_{|j|,j}=N^{|j|}_{|j|,j} u_j \in M_j$.     
As $M_j\in\,\!_\cX \cM^\cA$, it follows by the definition of $v^l_{kj}$ 
that $v^l_{kj}\in M_j$ for all $l=|j|,|j|+1,\dots $ 
and $k=-l,-l+1,{\dots}\,,l$.
Since $\cA$ is the direct sum of $M_j$, 
we conclude that 
$M_j = \Lin \{v^l_{k,j}\,;\, l=|j|,|j|+1,\dots ,\ 
k=-l,-l+1,{\dots}\,,l\}$ 
and the decomposition $\cA=\oplus_{j\in\frac{1}{2}\dZ}~M_j$ 
is an orthogonal sum.

Writing
\begin{align*}
v^l_{kj} 
&=N^l_{kj} \big((F^{l-k})_{(1)} \ang x^{l-|j|}_1\big) 
(F^{l-k})_{(2)} \ang u_j\\
&=N^l_{kj} (N^{|j|}_{|j|,j})^{-1} 
\big((F^{l-k})_{(1)} \ang x^{l-|j|}_1\big) (F^{l-k})_{(2)} \ang v^{|j|}_{|j|,j}
\end{align*}
and keeping in mind that
$(F^{l-k})_{(2)} \ang v^{|j|}_{|j|,j}\in \Lin \{v^{|j|}_{kj}\,;\, 
k=-|j|, -|j|+1,{\dots}\,,|j|\}$, it is clear that $M_j$ is generated by 
$\{v^{|j|}_{kj}\,;\,k=-|j|,-|j|+1,{\dots}\,,|j|\}$ as a left $\cX$-mo\-dule
and that $v^{|j|}_{kj}$ is cyclic for the 
$\ast$-re\-pre\-sen\-ta\-tion $\hat\pi_j$ 
of $\cX\rti\cU$ on $M_j$.

We turn now to the projections of the projective 
modules $M_j$. 
Defining
\begin{equation}\label{vj}
v_j := [2|j| \!+\!1]^{-1/2} 
(q^{|j|}v^{|j|}_{-|j|,j},{\dots}\,,q^{-|j|} v^{|j|}_{|j|,j})^{\mathrm{t}},
\end{equation}
we can write $P_j = v_jv^\ast_j$. This immediately implies that 
$P^\ast_j = P_j$. 
In order to prove $P^2_j=P_j$, it is  sufficient to show that 
$v^\ast_j v_j = 1$. 
Recall that $K\ang v^l_{kj} = q^k v^l_{kj}$ and 
$F \ang v^l_{kj} = [l\!+\!k]^{1/2} [l\!-\!k+1]^{1/2} v^l_{k-1,j}$
(see also Equation (\ref{kefop}) below).
By the third equation in (\ref{E1}), 
$K\ang v^{l\ast}_{kj} =q^{-k}v^{l\ast}_{kj}$ 
and 
$F \ang v^{l\ast}_{kj} = 
-q^{-1} [l\!-\!k]^{1/2} [l\!+\!k\!+\!1]^{1/2} v^{l\ast}_{k+1,j}$. 
From this, we conclude that 
$v^\ast_j v_j 
= [2|j|\!+\!1]^{-1} \big(q^{2|j|} v^{|j|\ast}_{-|j|,j} v^{|j|}_{-|j|,j}
+{\dots}+ q^{-2|j|}v^{|j|\ast}_{|j|,j} v^{|j|}_{|j|,j}\big)$ 
is a linear combination of vectors of weight $0$.
We  show that $v^\ast_j v_j$ belongs to a spin $0$ representation. 
Since $K \ang v_j^\ast v_j = v_j^\ast v_j$, 
this is equivalent to $F\ang (v^\ast_j v_j)=0$.
Inserting the expressions for $v_j$ gives  
\begin{align*}
F \ang (v^\ast_j v_j) 
&= \sumop^{|j|}_{k=-|j|} q^{-2k} [2|j|\!+\!1]^{-1} 
\Big((F \ang v^{|j|\ast}_{kj} )(K \ang v^{|j|}_{kj})
+(K^{-1} \ang v^{|j|\ast}_{kj}) (F\ang v^{|j|}_{kj})\Big)\\
&=[2|j|\!+\!1]^{-1} \sumop^{|j|}_{k=-|j|}\Big( 
-q^{-k-1}[|j|\!-\!k]^{1/2} [|j|\!+\!k\!+\!1]^{1/2} 
v^{|j|\ast}_{k+1,j} v^{|j|}_{kj}\\
&\qquad 
+q^{-k} [|j|\!+\!k]^{1/2}[|j|\!-\!k\!+\!1]^{1/2} 
v^{|j|\ast}_{kj} v^{|j|}_{k-1,j}\Big)
\end{align*}
which telescopes to zero. Since $v^\ast_j v_j$ belongs to 
a spin $0$ representation and 
$h(v^{|j|\ast}_{kj} v^{|j|}_{k,j})\hsp=\hsp\|v^{|j|}_{k,j}\|^2=1$, 
we have
$
v^\ast_j v_j
\hsp =\hsp h(v^\ast_jv_j)\hsp= \hsp\sum^{|j|}_{k=-|j|} 
q^{-2k} [2|j|{+}1]^{-1}\hsp=\hsp 1
$
by (\ref{vj}) as desired. Hence  $ P_j$ is an orthogonal projection.

Next we verify that $v^{|j|}_{nj}v^{|j|\ast}_{mj}$ 
belongs to $\cX$. In order to do so,  
we use the fact that $\cX$ is the set of elements 
$x\in \cA$ such that $x\anf X_r=0$. 
Since $v^{|j|}_{nj}\anf X_r=\hat{X_r}(v^{|j|}_{nj})=\mu_j v^{|j|}_{nj}$ 
and $v^{|j|\ast}_{mj}= N^{|j|}_{mj}\, S(F^{|j|-m})^\ast\ang u^\ast_j$, 
we get
\begin{align*}
(v^{|j|}_{nj}v^{|j|\ast}_{mj})\anf X_r
&\!=\!N^{|j|}_{mj}\, [v^{|j|}_{nj}\big(S(F^{|j|-m})^\ast \ang u^\ast_j\big) 
\anf X_r+(v^{|j|}_{nj}\anf X_r)
\big(S(F^{|j|-m})^\ast \ang u^\ast_j\big)\anf K^2]\\
&\!=\!N^{|j|}_{mj}\, v^{|j|}_{nj} S(F^{|j|-m})^\ast \ang [u^\ast_j \anf X_r 
+ \mu_j u^\ast_j \anf K^2].
\end{align*}
Hence it suffices to show that 
$u^\ast_j \anf X_r + \mu_j u^\ast_j\anf K^2=0$ for all 
$j\in \frac{1}{2}\dZ$. This can be done by induction.
Let $r\in (0,\infty)$. 
For $j=0$, the assertion is true since $u_0^\ast\anf X_r=0$ and $\mu_0=0$. 
Assume that $u^\ast_j\anf X_r + \mu_j u_j^\ast \anf K^2=0$ holds for 
$j\in \frac{1}{2}\dN_0$. Then 
$u^\ast_j \anf X_r = -\mu_j (u^\ast_j \anf K^2)$. 
As $u^\ast_{j+1/2} = (a-q^{-2j} s c)u^\ast_j$, 
we compute
\begin{align*}
&u^\ast_{j+1/2} \anf X_r 
+ \mu_{j+1/2} (u^\ast_{j+1/2}  \anf K^2)
= (a-q^{-2j} s c)(u^\ast_j \anf X_r)\\
&\qquad +(a-q^{-2j} s c) \anf X_r (u^\ast_j \anf K^2) 
+ \mu_{j+1/2}  (a-q^{-2j} s c)\anf K^2 
(u^\ast_j \anf K^2)\\
&=\big[-\mu_j(a-q^{-2j} s c)
+\big(q^{1/2} r^{-1/2}
 (q^{-1} - q)^{-1} (1-q^{-1}) - q^{-1/2}q^{-2j} s \big) a\\
&\qquad - \big( q^{1/2} r^{-1/2} (q^{-1}-q)^{-1}(1-q) 
- q^{2j+3/2} s^{-1}\big) q^{-2j} s c\\
&\qquad \qquad + \mu_{j+1/2} (q^{-1} a
- q^{-2j+1} s c)\big] (u^\ast_j \anf K^2).
\end{align*}
Inserting the expressions for $s$, $\mu_j$ and $\mu_{j+1/2}$, 
the preceding equation yields zero. 
This proves that $u^\ast_j \anf X_r + \mu_j u^\ast_j \anf K^2=0$ and 
hence $(v^{|j|}_{nj} v^{|j|\ast}_{mj})\anf X_r=0$ 
for $j\in \frac{1}{2}\dN_0$.
In the same way, one can show that  
$w^\ast_j\anf  X_r+\mu_{-j} w^\ast_j\anf K^2=0$ for $j\in\frac{1}{2}\dN_0$.
Since  
$u^\ast_{-j}\anf X_r + \mu_{-j} u^\ast_{-j} \anf K^2 
=  S(E^{2j})^\ast \ang [w^\ast_j \anf X_r + \mu_{-j} w^\ast_j\anf K^2]$, 
we conclude that $u^\ast_j\anf X_r + \mu_{j} u^\ast_{j} \anf K^2 = 0$ 
holds for all $j\in\frac{1}{2}\dZ$.
For $r=0$ and $r=\infty$, the proof is similar.

Finally we prove that $M_j$ is isomorphic 
to $\cX^{2|j|+1} P_j$ as a left  $\cX$-mo\-dule. 
Define a mapping $\Psi_j:\cX^{2|j|+1} P_j\rightarrow M_j$ by
\begin{equation}                                            \label{psi}
\Psi_j ((y_{-|j|},{\dots}\,, y_{|j|})P_j) :=  
[2|j|\!+\!1]^{-1/2}\! \sumop^{|j|}_{k=-|j|} q^{-k} y_k v^{|j|}_{kj}.
\end{equation}
Recall that $P_j=v_jv^\ast_j$ and $v^\ast_j v_j=1$. 
Suppose we are given  $y_{-|j|},{\dots}\,,y_{|j|}\in\cX$ such that 
$(y_{-|j|},{\dots}\,,y_{|j|})P_j=0$. 
Multiplying by $v_j$ from the right yields
$$
0\!\hspace{1pt}=\!\hspace{1pt}(y_{-|j|},{\dots}\,, y_{|j|})P_jv_j
=(y_{-|j|},{\dots}\,, y_{|j|}) v_j (v^\ast_j v_j)
=[2|j|\!+\!1]^{-1/2}\!\! \sumop^{|j|}_{k=-|j|}\! q^{-k} y_k v^{|j|}_{kj}.
$$
Hence $\Psi_j$ is well defined. 
Furthermore,
$[2|j|\!+\!1]^{-1/2} \sum^{|j|}_{k=-|j|} 
q^{-k} y_k v^{|j|}_{kj} =0$ implies
$$0=\big([2|j|\!+\!1]^{-1/2}  
\!\sumop^{|j|}_{k=-|j|} q^{-k} y_k v^{|j|}_{kj}\big)v^\ast_j 
= (y_{-|j|},{\dots}\,, y_{|j|}) P_j,
$$
so $\Psi_j$ is injective. 
Since $M_j$ is generated by $\{v^{|j|}_{kj}\,;\, k=-|j|,{\dots}\,,|j|\}$ 
as a left $\cX$-mo\-dule, $\Psi_j$ is also surjective. 
Whence $\Psi_j$ realizes 
the desired isomorphism. This completes the proof.\hfill $\Box$
%
%
\section{Integrable
{\mathversion{bold}$\ast$}-re\-pre\-sen\-ta\-tions of the cross product algebras }
                                                                 \label{S5}
%
%
%
\subsection{Classification of integrable
{\mathversion{bold}$\ast$}-re\-pre\-sen\-ta\-tions of the cross product algebra 
{\mathversion{bold}$\cO({\mathrm{S}}^2_{qr})\rti\cU_q({\mathrm{su}}_2)$} }
                                                     \label{sec-S5}
For $l\in\frac{1}{2}\dN_0$, let $T_l$ denote the type 1 spin $l$ 
representations of $\cU_q(\mathrm{su}_2)$. Recall that $T_l$ is 
an irreducible $\ast$-rep\-re\-sen\-ta\-tion of the 
$\ast$-al\-ge\-bra $\cU_q(\mathrm{su}_2)$ 
acting on a $(2l+1)$-dimensional Hilbert space with orthonormal 
basis $\{v^l_j\,;\,j=-l,-l+1,{\dots}\,,l\}$ by the formulas 
(see, for instance, \cite[Subsection 3.2.1]{KS})
\begin{equation}                                        \label{kefop}
K v^l_j = q^j v^l_j,\ E v^l_j 
= [l\!-\!j]^{1/2} [l\!+\!j\!+\!1]^{1/2} v^l_{j+1},\
F v^l_j = [l\!-\!j\!+\! 1]^{1/2} [l\!+\!j]^{1/2} v^l_{j-1}.
\end{equation}
As mentioned in the introduction, we call a
$\ast$-rep\-re\-sen\-ta\-tion of 
$\cO(\mathrm{S}^2_{qr})\rti \cU_q(\mathrm{su}_2)$
{\it integrable} if it has the 
following property:

\medskip
{\it The restriction to $\cU_q(\mathrm{su}_2)$ is the direct 
sum of representations $T_l$, $l\in \frac{1}{2}\dN_0$.} 
\medskip\\
The aim of this subsection 
is classify all integrable $\ast$-rep\-re\-sen\-ta\-tions of the cross product 
$\ast$-al\-ge\-bra  
$\cO(\mathrm{S}^2_{qr})\rti \cU_q(\mathrm{su}_2)$.
We shall use the generators $x_{j}$ defined by (\ref{xgener}) 
(resp.\ (\ref{xgener1})) rather than $A,B,B^\ast$.

Suppose we have such a representation acting on the domain $\Hh$. 
Then $\Hh$ can be written as 
$\Hh=\oplus_{l\in\frac{1}{2}\dN_0}V_l$ 
such that $V_l=\oplus_{j=-l}^{l}V^l_j$, where each $V^l_j$ is the same 
Hilbert space $V^l_l$, say, and the 
generators of\ \,$\cU_q(\mathrm{su}_2)$ 
act on $V_l$ by (\ref{kefop}). 
We claim that there exist operators 
$\alpha^\pm (l,j)\hsp:\hsp 
V^l_j\rightarrow V^{l\pm 1}_{j+1}$,\ \,$\alpha^0 (l,j)\hsp:\hsp
V^l_j\rightarrow V^{l}_{j+1}$,\ \,$\beta^+ (l,j)\hsp:\hsp
V^l_j\rightarrow V^{l+1}_{j}$, and 
self-adjoint operators 
$\beta^0 (l,j)\hsp:\hsp V^l_j\rightarrow V^{l}_{j}$ 
such that 
\begin{align}\label{xrel1}
x_1 v^l_j &= \alpha^+ (l,j) v^{l}_{j} + \alpha^0 (l,j) v^l_{j} 
+ \alpha^- (l,j) v^{l}_{j},\\
\label{xrel2}
x_0 v^l_j &= \beta^+ (l,j) v^{l}_j + \beta^0 (l,j) v^l_j + 
\beta^+ (l\!-\!1,j)^\ast v^{l}_j,\\
\label{xrel3}
x_{-1} v^l_j &= -q^{-1} 
\big(\,\alpha^-(l\!+\!1, j\!-\!1)^\ast v^{l}_{j} + 
\alpha^0(l,j\!-\!1)^\ast v^l_{j} +
\alpha^+(l\!-\!1,j\!-\!1)^\ast v^{l}_{j}\big) 
\end{align}
for $v^l_j\in V^l_j$. Indeed, let $v^l_j\in V^l_j$.  
Since $Kx_1=q x_1 K$, $x_1v^l_j$ is a weight vector with weight $j{+}1$,
and since 
$E^{l-j+1} x_1 v^l_j=q^{-l+j-1} x_1 E^{l-j+1} v^l_j=0$, $x_1 v^l_j$ is in 
the linear span of vectors  $w^r_{j+1}\in V^r_{j+1}$, 
where $r\le l{+}1$. Similarly, 
replacing $x_1$ by $x_{-1}$ and $E$ by $F$, we conclude 
that $x_{-1} v^l_j$ 
belongs to the span of vectors  $w^r_{j-1}\in V^r_{j-1}$, $r\le l{+}1$. 
Therefore, since $x_{-1}=-q^{-1} x^\ast_1$, we have
$x_{\pm 1} v^l_j\in 
V^{l-1}_{j\pm 1}\oplus V^l_{j\pm 1}\oplus V^{l+1}_{j\pm 1}$, 
so $x_{\pm 1} v^l_j$ is of the form (\ref{xrel1}) (resp.\ (\ref{xrel3})). 
From the last two relations of (\ref{podrel}) (resp.\ (\ref{prodrel1})) 
and from $x_0^\ast=x_0$, it follows that $x_0v^l_j$ 
is of the form (\ref{xrel2}). 
Note that all operators $\alpha^\pm (l,j)$, $\alpha^0(l,j)$,  
$\beta^{+} (l,j)$, 
$\beta^{0} (l,j)$  
are bounded because the operators $x_{-1}$, $x_0$, $x_1$ 
are bounded for any $\ast$-re\-pre\-sen\-ta\-tion 
of the $\ast$-al\-ge\-bra
$\cO(\mathrm{S}^2_{qr})$.

Inserting (\ref{kefop}) and (\ref{xrel1})
into the equation  $E x_1 v^l_j = q^{-1} x_1 E v^l_j$, we get
\begin{align*}
[l\!-\!j]^{1/2} [l\!+\!j\!+\!3]^{1/2} \alpha^+(l,j)
&=q^{-1} [l\!-\!j]^{1/2} [l\!+\!j\!+\!1]^{1/2} \alpha^+(l,j\!+\!1),\\
[l\!-\!j\!-\!1]^{1/2} [l\!+\!j\!+\!2]^{1/2} \alpha^0(l,j)
&=q^{-1} [l\!-\!j]^{1/2} [l\!+\!j\!+\!1]^{1/2} \alpha^0(l,j\!+\!1),\\
[l\!-\!j\!-\!2]^{1/2} [l\!+\!j\!+\!1]^{1/2} \alpha^-(l,j)
&=q^{-1} [l\!-\!j]^{1/2} [l\!+\!j\!+\!1]^{1/2} \alpha^-(l,j\!+\!1).
\end{align*}
The solutions of these recurrence relations are given by 
\begin{align*}
\alpha^+ (l,j) &= q^{-l+j} [l\!+\!j\!+\!1]^{1/2} 
[l\!+\!j\!+\!2]^{1/2} [2l\!+\!1]^{-1/2} [2l\!+\!2]^{-1/2} \alpha^+(l,l),\\
\alpha^0 (l,j) &= q^{-l+j+1} [l\!-\!j]^{1/2} [l\!+\!j\!+\!1]^{1/2} [2l]^{-1/2} 
\alpha^0(l,l\!-\!1),\\
\alpha^- (l,j) &= q^{-l+j+2} [l\!-\!j\!-\!1]^{1/2} [l\!-\!j]^{1/2} [2]^{-1/2} 
\alpha^-(l,l\!-\!2).
\end{align*}
Similarly, from the equation $E x_0 v^l_j=(x_0E+[2]^{1/2} x_1 K) v^l_j$,  
we obtain
\begin{align*}
&[l{-}j{+}1]^{1/2} [l{+}j{+}2]^{1/2} \beta^+(l,j)
\!=\![l{-}j]^{1/2} [l{+}j{+}1]^{1/2} 
\beta^+(l,j{+}1)\!+\![2]^{1/2} q^j \alpha^+(l,j),\\
&[l{-}j]^{1/2} [l{+}j{+}1]^{1/2} \beta^0(l,j)=[l{-}j]^{1/2} [l{+}j{+}1]^{1/2} 
\beta^0(l,j\!+\!1)+[2]^{1/2} q^j \alpha^0(l,j).
\end{align*}
The equation $E x_0 v^l_l =(x_0 E+[2]^{1/2} x_1 K)v^l_l$ yields in addition 
\begin{equation}                                      \label{baplus}
\beta^+ (l,l)= q^l[2]^{1/2} [2l\!+\!2]^{-1/2} \alpha^+(l,l).
\end{equation}
Further, the equation 
$0=\langle v^l_l, (x_1F+q[2]^{1/2} x_0 K-qF x_1)v^l_l\rangle$ implies that 
\begin{equation}                                        \label{banull}
\alpha^0 (l,l\!-\!1)= -[2]^{1/2} [2l]^{-1/2} q^{l+1} \beta^0(l,l).
\end{equation}
As a consequence, $\alpha^0 (l,l\!-\!1)$ is self-adjoint. 
Using (\ref{baplus}) and (\ref{banull}), it follows that the above 
recurrence relations for $\beta^+(l,j)$ and $\beta^0(l,j)$ have the 
following solutions:
\begin{align}    \label{bplus}
\beta^+ (l,j) &= 
q^j[l\!-\!j\!+\!1]^{1/2} 
[l\!+\!j\!+\!1]^{1/2}[2]^{1/2}[2l\!+\!1]^{-1/2} [2l\!+\!2]^{-1/2} 
\alpha^+(l,l),\\
\beta^0 (l,j)
&=(1-q^{l+j+1}[l\!-\!j][2][2l]^{-1})\beta^0(l,l).\label{bnull}&
\end{align}
From $0=\langle v^l_{l-1}, 
(E x_{-1}-qx_{-1} E-[2]^{1/2} x_0 K)v^{l-1}_{l-1}\rangle$, 
we derive  $\beta^+(l\!-\!1,l\!-\!1)=
-q^{-l}[2l\!-\!1]^{1/2} \alpha^-(l,l\!-\!2)^\ast$. 
Combining this equation  with (\ref{baplus}) gives
\begin{equation}\label{aplusminus}
\alpha^- (l,l\!-\!2)= 
-q^{2l-1} [2]^{1/2}[2l\!-\!1]^{-1/2} [2l]^{-1/2}
\alpha^+(l\!-\!1,l\!-\!1)^\ast.
\end{equation}
From (\ref{baplus})--(\ref{aplusminus}), it follows that the representation 
is completely described if the operators $\alpha^+(l,l)$ and $\beta^0(l,l)$ 
for $l\in\frac{1}{2}\dN_0$ are known. Our next aim is to 
determine these operators. This will be done by induction on $l$, 
where consecutive steps increase $l$ by 1.

We begin by analyzing  the relations between $\beta^0(l,l)$ and 
$\alpha^+(l,l)$. Let $r<\infty$ and abbreviate $\rho := 1+[2]^2 r$. 
On $V^l_l$, the two equations $BB^\ast = r + q^2 A - q^4 A^2$ and  
$B^\ast B = r+A-A^2$ lead to the operator relations
\begin{align*}
&|\alpha^+ (l\!-\!1,l\!-\!1)^\ast|^2 
+ \alpha^0 (l,l\!-\!1)^2 
+|\alpha^- (l\!+\!1,l\!-\!1)^\ast|^2 \\
&\qquad \qquad =(1+q^2)^{-1} \big(q^2 \rho + (1-q^2) \beta^0 (l,l) 
- (\beta^0 (l,l)^2 + |\beta^+ (l,l)|^2)\big),\\
&|\alpha^+(l,l)|^2 
= (1\!+\!q^2)^{-1} \big(q^2 \rho - (1\!-\!q^2) q^2 \beta^0(l,l) 
- q^4 (\beta^0 (l,l)^2 + |\beta^+ (l,l)|^2)\big).
\end{align*}
Here it is understood that $\alpha^0 (0,-1)=0$ and 
$\alpha^+ (l\!-\!1, l\!-\!1)=0$ for $l=0,1/2$. 
Inserting (\ref{baplus}), 
(\ref{banull}) and (\ref{aplusminus}), these two relations give 
after some calculations
\begin{align}                                              
&|\alpha^+ (l\!-\!1, l\!-\!1)^\ast|^2 
+ q^{2l+1}[2l\!+\!1]^{-1} [2l\!+\!2]^{-1}[2l\!+\!3] |\alpha^+(l,l)|^2 
\nonumber\\
&\qquad 
=[2]^{-1}\big(q\rho + (q^{-1} - q) \beta^0 (l,l) - q^2 [2l]^{-1} ([2l\!+\!3] 
+ q^{2l}) \beta^0 (l,l)^2\big),                           \label{alpharel-}\\
&[2l\!+\!2]^{-1} [2l\!+\!3]|\alpha^+ (l,l)|^2        \label{alpharel}
=[2]^{-1}\big(\rho - (1-q^2)\beta^0(l,l) - q^2 \beta^0 (l,l)^2\big).
\end{align}
Eliminating $\alpha^+(l,l)$ from these equations yields 
\begin{equation}\label{beta0}
[2][2l\!+\!1] |\alpha^+ (l\!-\!1, l\!-\!1)^\ast|^2 
{=}[2l]\rho \!+\! (1\!-\!q^2)[2l\!+\!2] \beta^0 (l,l) 
\!-\! [2l]^{-1} [2l\!+\!2]^2 q^2 \beta^0(l,l)^2.
\end{equation}

Now we start with the induction procedure. Our first aim is to show that
$\beta^0 (0,0)= 0$. 
A computation similar to the above shows that 
\begin{align*}
&q[2]^{-1} [3] |\alpha^+ (0,0)|^2 
= (1+ q^2)^{-1} \big(q^2 \rho - (1-q^2)q^2 \beta^0 (0,0) 
- q^4 \beta^0 (0,0)^2\big),\\
&q[2]^{-1} [3] |\alpha^+ (0,0)|^2 
= (1+ q^2)^{-1} \big(q^2 \rho + (1-q^2) \beta^0 (0,0) -  \beta^0 (0,0)^2\big).
\end{align*}
Eliminating $|\alpha^+ (0,0)|^2$ gives $0=\beta^0 (0,0)-\beta^0 (0,0)^2$. 
Since $\beta^0 (0,0)$ is self-adjoint, it is an orthogonal projection.  
Assume to the contrary that $\beta^0 (0,0)\ne 0$. 
Then there is $v^0_0 \in V^0_0$ such that 
$\beta^0 (0,0) v^0_0 = v^0_0$ and $\|v^0_0\|=1$. 
Note that $h(\cdot) := \langle (\cdot) v^0_0, v^0_0\rangle$ 
is the unique $\cU_q(\mathrm{su}_2)$-in\-var\-iant state on 
$\cO(\mathrm{S}^2_{qr})$.
From this, we conclude 
$1=\langle \beta^0 (0,0)v^0_0, v^0_0\rangle 
= \langle x_0v^0_0, v^0_0\rangle = h(x_0)=0$,
which is a contradiction. Thus $\beta^0 (0,0)=0$.

For $l=1/2$, Equation (\ref{beta0}) becomes 
$$
0=(q[3] \beta^0(1/2, 1/2))^2 
- (1-q^2)[3] \beta^0(1/2, 1/2)-\rho.
$$
This operator identity can only be satisfied if the self-adjoint 
operator $\beta^0(1/2, 1/2)$ has purely discrete 
spectrum with eigenvalues
$$
\beta^0(1/2), 1/2)^\pm_{1/2} 
:= [3]^{-1} (q^{-2} \lambda_\pm - \lambda_\mp).
$$
Clearly, $\beta^0(1/2, 1/2)^+_{1/2}\ne \beta^0(1/2, 1/2)^-_{1/2}$. 
Denoting the corresponding 
eigen\-spaces by $\cK_{-1/2}$ and $\cK_{1/2}$, 
we can write 
$V^{1/2}_{1/2} = \cK_{-1/2} \oplus \cK_{1/2}$, 
and $\beta^0(1/2, 1/2)$ acts on $V^{1/2}_{1/2}$ by
\begin{align*}
&\beta^0(1/2, 1/2) w_{-1/2 } 
= \beta^0(1/2, 1/2)^-_{1/2} w_{-1/2}, \quad 
w_{-1/2} \in \cK_{1/2},\\
&\beta^0(1/2, 1/2) w_{1/2 } 
= \beta^0(1/2, 1/2)^+_{1/2} w_{1/2},  \quad 
w_{1/2} \in \cK_{1/2}.
\end{align*}
Here we do not exclude the cases 
$\cK_{-1/2}=\{0\}$ and $\cK_{1/2}=\{0\}$.

Next let $l\in\frac{1}{2}\dN_0$. Assume that there exist 
(possibly zero) Hilbert spaces $\cK_{-l}, \cK_{-l+1},{\dots}\,,\cK_l$ 
such that 
$V^l_l=\cK_{-l}\oplus{\dots}\oplus\cK_l$. For $j\in\frac{1}{2}\dN_0$, 
$j\le l$, set 
\begin{align}                                        \label{beta+-}
\beta^0(l,l)_j^\pm
&= [2l\!+\!2]^{-1}\big([2j](q^{-2} \lambda_\pm - \lambda_\mp) 
- (1-q^{-2})[l\!-\!j][l\!+\!j\!+\!1]\big),\\          
\alpha^+(l,l)^\pm_j &= [2]^{1/2}[2l\!+\!3]^{-1/2}[2l\!+\!2]^{-1/2} 
\big( [2l\!+\!2]^2 (c\!+\!1/4)                                 \nonumber\\
&\quad -\{(q^{-1}\!-\!q) [l\!-\!j\!+\!1][l\!+\!j\!+\!1]/2
\pm [2j](c\!+\!1/4)^{1/2} \}^2\big)^{1/2}.              \label{alpha+-}
\end{align}
where we already inserted (\ref{beta+-}) into (\ref{alpharel}). 
Assume that $\beta^0(l,l)$ acts on $V^l_l$ by 
\begin{equation}                                        \label{K+}
\beta^0(l,l)w_{-j} = \beta^0(l,l)^-_j w_{-j},\ \, w_{-j} \in \cK_{-j},\quad 
\beta^0(l,l) w_j = \beta^0(l,l)^+_j w_j,\ \, w_j\in \cK_j.
\end{equation}
We show that there exist Hilbert spaces $\cK_{-(l+1)}$ and 
$\cK_{l+1}$ such that, up to unitary equivalence, 
$V^{l+1}_{l+1} 
= \cK_{-(l+1)} \oplus \cK_{-l} \oplus{\dots}\oplus\cK_l\oplus\cK_{l+1}$, 
the operator $\beta^0(l+1,l+1)$ acts on $V^{l+1}_{l+1}$ 
by the formulas (\ref{K+}), and
$\alpha^+(l,l):V^l_l\rightarrow V^{l+1}_{l+1}$ is given by
\begin{equation}                                        \label{alphaact}
\alpha^+(l,l) w_{-j} =\alpha^+(l,l)^-_j w_{-j},\ \, w_{-j} \in \cK_{-j}, 
\ \ \alpha^+(l,l)w_j = \alpha^+(l,l)^+_j w_j,\ \,w_j \in \cK_j, 
\end{equation}
where $j=-l,-l+1,{\dots}\,,l$.

Observe that $\ker B =\{0\}$ (cf.\ Subsection \ref{sec-4} below). 
Hence $\ker \alpha^+(l,l)=\{0\}$. Let $\alpha^+(l,l)=U|\alpha^+(l,l)|$ 
denote the polar decomposition of $\alpha^+(l,l)$. 
Then $U$ is an isometry from $V^l_l$ onto $\overline{\alpha^+(l,l)V^l_l}$. 
On the other hand, we have the decomposition 
$V^{l+1}_{l+1} 
= \overline{\alpha^+(l,l)V^l_l}\oplus \ker \alpha^+(l,l)^\ast$. 
After applying a unitary transformation, we can assume that
$V^{l+1}_{l+1} = V^l_l \oplus \ker \alpha^+(l,l)^\ast 
= \cK_{-l} \oplus {\dots} \oplus \cK_l \oplus \ker \alpha^+(l,l)^\ast$ 
and $\alpha^+(l,l)v^l_l = |\alpha^+(l,l)|v^l_l$ for all $v^l_l \in V^l_l$.
From (\ref{alpharel}) and (\ref{K+}), it follows that the action 
of $\alpha^+(l,l)$ on $V^l_l$ is determined by Equation (\ref{alphaact}) 
as asserted. We proceed by describing the action of $\beta^0(l\!+\!1,l\!+\!1)$.
By (\ref{beta0}), $|\alpha^+(l,l)^\ast|^2$ commutes 
with $\beta^0(l\!+\!1,l\!+\!1)$. Hence $\beta^0(l\!+\!1,l\!+\!1)$ leaves the 
subspaces $V^l_l$ and $\ker \alpha^+(l,l)^\ast$ of $V^{l+1}_{l+1}$ 
invariant. Let $\tilde{\beta^0}(l\!+\!1,l\!+\!1)$ and 
$\hat{\beta^0}(l\!+\!1,l\!+\!1)$ denote the restrictions 
of $\beta^0(l\!+\!1,l\!+\!1)$ to $V^l_l$ and $\ker \alpha^+(l,l)^\ast$, 
respectively. 
Evaluating  the relation 
$0=\langle (x_0x_1-q^2 x_1 x_0 - (1-q^2)x_1)v^l_l, 
v^{l+1}_{l+1}\rangle$ for  $v^l_l \in V^l_l$ and 
$v^{l+1}_{l+1}\in V^{l+1}_{l+1}$ 
yields
\begin{multline*}
\beta^0(l\!+\!1,l\!+\!1) \alpha^+(l,l)
-q^2 \alpha^0(l\!+\!1,l\!+\!1) \beta^+(l,l) \\
- q^2 \alpha^+(l,l) \beta^0(l,l) - (1-q^2) \alpha^+(l,l)=0.
\end{multline*} 
Inserting (\ref{baplus}) and observing that  
$\alpha^+(l,l)   \beta^0(l,l)= \beta^0(l,l) \alpha^+(l,l)$ 
by (\ref{K+}) and (\ref{alphaact}), 
we deduce the following operator equation on $V^l_l$
$$\{ (1+q^{2l+4} [2][2l\!+\!2]^{-1}) \tilde{\beta}^0 (l\!+\!1,l\!+\!1) 
-q^2 \beta^0 (l,l) - (1-q^2)\} \alpha^+ (l,l)= 0.
$$ 
Since $V^l_l = \overline{\alpha^+(l,l)V^l_l}$, the operator
in braces must be zero.
This implies that 
$$
\tilde{\beta^0} (l\!+\!1,l\!+\!1)
=[2l\!+\!4]^{-1} [2l\!+\!2] (\beta^0(l,l)+q^{-2}-1)
$$
Hence the operator 
$\beta^0(l\!+\!1,l\!+\!1)$ acts on 
$V^l_l=\cK_{-l} \oplus{\dots}\oplus\cK_l$ by
\begin{multline*}
\beta^0(l\!+\!1,l\!+\!1)w_j 
=[2l\!+\!4]^{-1} [2l\!+\!2](\beta^0(l,l)^\epsilon_j + q^{-2} -1) w_j\\
=[2l\!+\!4]^{-1} 
\left( [2j] 
(q^{-2} \lambda_\epsilon - \lambda_{-\epsilon}
+(q^{-2} -1)[l\!-\!j\!+\!1][l\!+\!j\!+\!2]
\right) w_j,
\end{multline*}
where $\epsilon=\sign(j)$. 
The last equation is obtained by straightforward computations.
On $\ker \alpha^+(l,l)^\ast$, Equation (\ref{beta0}) reads  
$$
0\!=\!\rho+(1\!-\!q^2)[2l\!+\!2]^{-1}[2l\!+\!4]
\hat{\beta^0}(l\!+\!1,l\!+\!1)
-([2l\!+\!2]^{-1} [2l\!+\!4] q\hat{\beta^0}(l\!+\!1,l\!+\!1))^2.
$$ 
Since the solution of the quadratic equation $0=-\rho-(q^{-1}-q)t+t^2$ 
is given by $ t^\pm=q^{-1}\lambda_\pm-q\lambda_\mp$, 
we conclude that $\hat{\beta}(l+1,l+1)$ has purely 
discrete spectrum consisting of the two distinct eigenvalues
$$
\beta^0(l\!+\!1,l\!+\!1)^\pm_{l+1} 
=[2l\!+\!4]^{-1} [2l\!+\!2] (q^{-2} \lambda_\pm - \lambda_\mp),
$$
and $\ker \alpha^+(l,l)^\ast$ splits into 
the direct sum $\ker \alpha^+(l,l)^\ast = \cK_{-(l+1)} \oplus \cK_{l+1}$, 
where $\cK_{-(l+1)}$ and  $\cK_{l+1}$ denote the eigenspaces of 
$\hat{\beta^0}(l\!+\!1,l\!+\!1)$ 
corresponding to the eigenvalues 
$\beta^0(l\!+\!1,l\!+\!1)^-_{l+1}$ and 
$\beta^0(l\!+\!1,l\!+\!1)^+_{l+1}$, respectively.  Accordingly, the operator 
$\beta^0(l\!+\!1,l\!+\!1)$ acts on $\ker\alpha^+(l,l)^\ast$ by
\begin{align*}
\beta^0(l\!+\!1,l\!+\!1)w_{-(l+1)}
&=\beta^0(l\!+\!1,l\!+\!1)^-_{l+1}w_{-(l+1)},
\quad w_{-(l+1)} \in \cK_{-(l+1)},\\ 
\beta^0(l\!+\!1,l\!+\!1)w_{l+1}
&=\beta^0(l\!+\!1,l\!+\!1)^+_{l+1}w_{l+1}, 
\quad w_{l+1}\in \cK_{l+1}.
\end{align*}
This shows that $\beta^0(l\!+\!1,l\!+\!1)$ is of the 
same form as  $\beta^0(l,l)$. 
By induction, the operators $\beta^0(l,l)$ and $\alpha^+(l,l)$ 
are now completely determined.

In the case $r=\infty$, there are only minor changes in the 
preceding reasoning. 
Set $\rho=(q+q^{-1})^2$. Then,  
Equations (\ref{alpharel})--(\ref{beta0}) remain 
valid if one omits the first order term of $\beta^0(l,l)$. 
Equations  (\ref{beta+-}) and (\ref{alpha+-}) become
\begin{align}                 \label{betainfty}
\beta^0(l,l)^\pm_j &= \pm q^{-1}[2][2j][2l+2]^{-1},\\
\alpha^+(l,l)^\pm_j                                             \label{alphainfty}
&=[2]^{1/2}[2(l\!+\!j\!+\!1)]^{1/2}[2(l\!-\!j\!+\!1)]^{1/2}[2l\!+\!2]^{-1/2}[2l\!+\!3]^{-1/2}, 
\end{align}
and $\alpha^+(l,l)$ is again given by (\ref{alphaact}).

Note that whenever a Hilbert space $\cK_j$ is non-zero, 
it appears as a direct summand in each $V^l_k$, 
where $l=|j|,|j|+1,{\dots}$ and $k=-l,-l+1,{\dots}\,,l$.
Moreover, the generators of $\cO(\mathrm{S}^2_{qr})\rti \cU_q(\mathrm{su}_2)$ 
leave this decomposition invariant. Hence the representation 
of $\cO(\mathrm{S}^2_{qr})\rti\cU_q(\mathrm{su}_2)$ on 
$\oplus_{l\in\frac{1}{2}\dN_0}\! \oplus^l_{k=-l} V^l_k$ 
splits into a direct sum of representations on 
$\oplus_{n\in\dN_0}\! \oplus^{|j|+n}_{k=-(|j|+n)} \cK^{|j|+n}_k$, 
where each $\cK^{|j|+n}_k$ is the same Hilbert space $\cK_j$ 
and is considered as a direct summand of $V^{|j|+n}_k$.

It still remains to prove that Equations (\ref{kefop})--(\ref{xrel3}) 
define a representation of 
$\cO(\mathrm{S}^2_{qr})\rti\cU_q(\mathrm{su}_2)$ when 
we insert the expressions for the operators obtained in the previous 
discussion. This can be done by showing that the defining 
relations of $\cO(\mathrm{S}^2_{qr})$ and the cross relations of 
$\cO(\mathrm{S}^2_{qr})\rti \cU_q(\mathrm{su}_2)$ are satisfied. 
We have checked this; the details of these lengthy 
and tedious computations are omitted.

We summarize the outcome of above 
considerations in the next theorem. 
%
\begin{tht}                                                 \label{rep}
Each integrable $\ast$-re\-pre\-sen\-ta\-tion of  
the cross product $\ast$-al\-ge\-bra  
$\cO(\mathrm{S}^2_{qr})\rti \cU_q(\mathrm{su}_2)$  
is, up to unitary equivalence, a direct sum of representations $\pi_j^\pm$, 
$j\!\in\!\frac{1}{2}\dN_0$, of the following form: 
\mn

\noindent
The domain 
is the direct sum $\oplus_{l\in\dN_0}\!\oplus^{j+l}_{k=-(j+l)} \cK^{j+l}_k$, 
where each $\cK^{j+l}_k$ is the same Hilbert space $\cK$. 
The generators $E$, $F$, $K$ of $\cU_q(\mathrm{su}_2)$ act on 
$\oplus^{j+l}_{k=-(j+l)} \cK^{j+l}_k$ by (\ref{kefop}). 
The actions of the generators $x_1$, $x_0$, $x_{-1}$ of 
$\cO(\mathrm{S}^2_{qr})$ are  
determined by   
\begin{align*}
x_1 v^l_k =\  &q^{-l+k} [l\!+\!k\!+\!1]^{1/2} [l\!+\!k\!+\!2]^{1/2} 
[2l\!+\!1]^{-1/2} [2l\!+\!2]^{-1/2} \alpha^+ (l,l)^\pm_j v^{l+1}_{k+1}\\
        &-q^{k+2} [l\!-\!k]^{1/2} [l\!+\!k\!+\!1]^{1/2} [2]^{1/2} [2l]^{-1} 
\beta^0 (l,l)^\pm_j v^l_{k+1}\\
           &-q^{l+k+1} [l\!-\!k\!-\!1]^{1/2} [l\!-\!k]^{1/2} 
[2l\!-\!1]^{-1/2} [2l]^{-1/2} \alpha^+ (l\!-\!1,l\!-\!1)^\pm_j v^{l-1}_{k+1},
\\[-24pt]
\end{align*}
\begin{align*}
x_0 v^l_k =\ &q^k [l\!-\!k\!+\!1]^{1/2} [l\!+\!k\!+\!1]^{1/2} [2]^{1/2} 
[2l\!+\!1]^{-1/2} [2l\!+\!2]^{-1/2} \alpha^+ (l,l)^\pm_j v^{l+1}_k\\
           &+\big(1-q^{l+k+1} [l\!-\!k][2][2l]^{-1}\big)
            \beta^0 (l,l)^\pm_j v^l_k\\
           &+q^k [l\!-\!k]^{1/2} [l\!+\!k]^{1/2} [2]^{1/2} 
[2l\!-\!1]^{-1/2} [2l]^{-1/2} \alpha^+ (l\!-\!1, l\!-\!1)^\pm_j v^{l-1}_k,
\\[-24pt]
\end{align*}
\begin{align*}
x_{-1} v^l_k =\ &q^{l+k} [l\!-\!k\!+\!1]^{1/2} [l\!-\!k\!+\!2]^{1/2} 
[2l\!+\!1]^{-1/2} [2l\!+\!2]^{-1/2} \alpha^+ (l,l)^\pm_j v^{l+1}_{k-1}\\
               &+q^k [l\!-\!k+1]^{1/2} [l\!+\!k]^{1/2} [2]^{1/2} 
[2l]^{-1} \beta^0 (l,l)^\pm_j v^l_{k-1}\\
               &-q^{-l+k-1} [l\!+\!k\!-\!1]^{1/2} [l\!+\!k]^{1/2} 
[2l\!-\!1]^{-1/2} [2l]^{-1/2} \alpha^+ (l\!-\!1,l\!-\!1)^\pm_j v^{l-1}_{k-1},
\end{align*}
where, for $r<\infty$, 
the real numbers $\beta^0(l,l)^\pm_j$ and $\alpha^+(l,l)^\pm_j$ 
are given by (\ref{beta+-}) and (\ref{alpha+-}), respectively, 
and, for $r=\infty$, 
by (\ref{betainfty}) and (\ref{alphainfty}), respectively. 
 
Representations corresponding to different pairs of labels $(j,\pm)$ 
are not unitarily equivalent 
(with only one obvious exception:
$\pi_0^+=\pi_0^-$). 
A representation of this list is irreducible if and only if $\cK=\dC$.
\end{tht}
An immediate consequence of Theorem \ref{rep} is the following
\begin{thc} \label{sumirr}                 
(i) Each integrable $\ast$-re\-pre\-sen\-ta\-tion of 
$\cO(\mathrm{S}^2_{qr})\rti
\cU_q(\mathrm{su}_2)$ is a direct sum of
  integrable irreducible $\ast$-re\-pre\-sen\-ta\-tions. \\
(ii) Each integrable irreducible  $\ast$-re\-pre\-sen\-ta\-tion of 
$\cO(\mathrm{S}^2_{qr})\rti\cU_q(\mathrm{su}_2)$ is unitarily equivalent 
to a $\ast$-re\-pre\-sen\-ta\-tions $\pi^{\pm}_j$, 
$j\in \frac{1}{2}\dN_0$, with $\cK=\dC$.
\end{thc}

It will be convenient to introduce the following notation. If 
$\pi^\pm_j$, $j\in\frac{1}{2}\dN_0$, denotes an {\it irreducible} 
integrable $\ast$-re\-pre\-sen\-ta\-tion of 
$\cO(\mathrm{S}^2_{qr})\rti\cU_q(\mathrm{su}_2)$ 
(that is, $\pi^{\pm}_j$ is
the representation from Theorem \ref{rep} with $\cK=\dC$), then we set
$$
\pi_j:= \pi^-_{-j}\ \,\mbox{for}\ \,j<0,\quad  
\pi_j:=\pi^+_j\ \,\mbox{for}\ \ j\geq 0,\quad j\in 
\mbox{$\frac{1}{2}$} \dZ.
$$ 

\subsection{Decomposition of tensor products of irreducible 
integrable representations with spin {\mathversion{bold}$l$} representations }
                                                     \label{sec-S6}
%
\begin{thl} \label{tensorpro}
Suppose that $\cX$ is a left module $\ast$-al\-ge\-bra of 
a Hopf $\ast$-al\-ge\-bra  $\cU$. 
Let $\pi$ and $T$ be $\ast$-re\-pre\-sen\-ta\-tions 
of the  $\ast$-al\-ge\-bras $\cX\rti \cU$ and $\cU$ 
on domains $\cD$ and $V$, respectively. Then
there is a $\ast$-re\-pre\-sen\-ta\-tion, denoted by $\pi \otimes T$, of the 
$\ast$-al\-ge\-bra $\cX \rti \cU$ on the domain $\cD \otimes V$ such that
$(\pi \otimes T)(x)= \pi(x)\otimes T(1)$ for 
$x\in \cX$ and $(\pi \otimes T)(f)= \pi(f_{(1)})\otimes T(f_{(2)})$ for
$f \in \cU$. 
\end{thl}
{\bf Proof.} It suffices to check that $\pi \otimes T$ respects
the cross relation (\ref{fxrel}), that is, 
$(\pi\hsp \otimes\hsp T)(f)\hs(\pi\hsp \otimes\hsp T)(x)
= (\pi\hsp\otimes\hsp T) (f_{(1)}\ang x)\hs(\pi\hsp\otimes\hsp T)(f_{(2)})$  
for $x \in \cX$ and $f \in \cU$. 
The details of this easy verification are left to the reader. 
                                                        \hfill$\Box$
\mn

Clearly, the tensor product $\pi\otimes T_l$ of an integrable
$\ast$-re\-pre\-sen\-ta\-tion $\pi$ of the cross product algebra 
$\cO(\mathrm{S}^2_{qr})\rti\cU_q(\mathrm{su}_2)$ 
and a spin $l$ representation $T_l$ 
of $\cU_q(\mathrm{su}_2)$ is again integrable,  
so Corollary \ref{sumirr} applies. 
The decomposition of $\pi_j\otimes T_l$ 
into a direct sum of irreducible integrable representations of 
$\cO(\mathrm{S}^2_{qr})\rti\cU_q(\mathrm{su}_2)$ is
described in the next proposition. 
\begin{thp}                                                 \label{deco}
For $j\in\frac{1}{2}\dZ$ and $l\in\frac{1}{2}\dN_0$, let 
$\pi_j$ be an irreducible integrable $\ast$-re\-pre\-sen\-ta\-tion 
of $\cO(\mathrm{S}^2_{qr})\rti\cU_q(\mathrm{su}_2)$ and $T_l$ a spin $l$ 
representation of $\cU_q(\mathrm{su}_2)$. 
Then, up to unitary equivalence, 
$\pi_j\otimes T_l=\pi_{j-l}\oplus \pi_{j-l+1}\oplus \ldots \oplus\pi_{j+l}$.
\end{thp}
{\bf Proof.} 
For $l=0$, there is nothing to prove. 
Let $l=1/2$. From Theorem \ref{rep}, it follows that the restriction of 
$\pi_j$ to $\cU_q(\mathrm{su}_2)$ 
is of the form $\oplus_{n\in\dN_0}T_{|j|+n}$. 
Assume first that $j\neq 0$. By the Clebsch--Gordon decomposition, we have
\begin{equation}                                             \label{T}
    \pi_j\otimes T_{1/2}=T_{|j|-1/2}\oplus(\oplus_{n\in\dN_0}2T_{|j|+n+1/2})
\end{equation}
as representations of $\cU_q(\mathrm{su}_2)$.
By Corollary \ref{sumirr}, $\pi_j\otimes T_{1/2}$ decomposes into 
a direct sum of irreducible integrable $\ast$-re\-pre\-sen\-ta\-tions 
described 
in Theorem \ref{rep}. 
In (\ref{T}), there occur no spin $k$ representations for 
$k<|j|-1/2$, exactly one spin $|j|-1/2$ representation, and two 
spin $k$ representations for $k=|j|+1/2,|j|+3/2,\ldots$\ \,Thus 
we must have  $\pi_j\otimes T_{1/2}
=\pi_{|j|-1/2}^{\epsilon_1}\oplus\pi_{|j|+1/2}^{\epsilon_2}$, 
where $\epsilon_1,\epsilon_2\in\{+,-\}$ are to be determined. 
Moreover, $\pi_{|j|-1/2}^{\epsilon_1}$ and $\pi_{|j|+1/2}^{\epsilon_2}$ 
are irreducible. Again by the Clebsch--Gordon decomposition, 
\begin{equation}             \label{cgd}
v^{|j|-1/2}_{|j|-1/2}= 
[2|j|\!+\!1]^{-1/2} \big(q^{1/2} [2|j|]^{1/2} 
v^{|j|}_{|j|}\otimes v^{1/2}_{-1/2} - q^{-|j|} 
v^{|j|}_{|j|-1}\otimes v^{1/2}_{1/2}\big)
\end{equation}
is the (unique) highest weight vector of the spin $|j|-1/2$ 
representation in (\ref{T}). If $|j|=1/2$, then $v^0_0$ belongs to 
a spin $0$ representation and
$\pi_{|j|-1/2}^{\epsilon_1}=\pi_0$ is the Heisenberg representation. 
Let  $|j|>1/2$ and $\epsilon =\sign(j)$. 
Then a straightforward computation gives 
$$
\langle v^{|j|-1/2}_{|j|-1/2},\pi_j\otimes T_{1/2}(x_0)
 v^{|j|-1/2}_{|j|-1/2}\rangle=
[2|j|\!+\!2][2|j|\!-\!1][2|j|\!+\!1]^{-1}[2|j|]^{-1}
\beta^0(|j|,|j|)^{\epsilon}_{|j|}.
$$
Hence $\epsilon_1\!\hspace{1pt}=\!\hspace{1pt}
\sign(\langle v^{|j|-1/2}_{|j|-1/2},\pi_j\otimes T_{1/2}(x_0)
 v^{|j|-1/2}_{|j|-1/2}\rangle)=
\sign(\beta^0(|j|,|j|)^{\epsilon}_{|j|})=\sign(j)$. 
The linear space of highest weight vectors belonging to  
spin $|j|\!+\!1/2$ representations in (\ref{T}) is 
2-dimensional and spanned by 
the orthonormal vectors 
\begin{align*}
u^{|j|+1/2}_{|j|+1/2}&= 
[2|j|\!+\!3]^{-1/2} \big(q^{1/2} [2|j|\!+\!2]^{1/2} 
v^{|j|+1}_{|j|+1}\otimes v^{1/2}_{-1/2} - q^{-|j|-1} 
v^{|j|+1}_{|j|}\otimes v^{1/2}_{1/2}\big),\\
w^{|j|+1/2}_{|j|+1/2}&=v^{|j|}_{|j|}\otimes v^{1/2}_{1/2}. 
\end{align*}
The vector 
\begin{multline*}
v:=
[2|j|\!+\!3]^{1/2}[2|j|\!+\!2]^{-1/2}
\alpha^+ (|j|,|j|)^\epsilon_{|j|}w^{|j|+1/2}_{|j|+1/2}\\
-q[2]^{1/2}[2|j|]^{-1/2}\beta^0(|j|,|j|)^{\epsilon}_{|j|}u^{|j|+1/2}_{|j|+1/2}
\end{multline*}
is orthogonal to $x_1v^{|j|-1/2}_{|j|-1/2}$. 
Hence $\|v\|^{-1}v$ is the (unique) highest weight vector 
of the spin $|j|\!+\!1/2$ representation 
belonging to the decomposition of $\pi_{|j|+1/2}^{\epsilon_2}$. 
Our goal is to determine 
$\epsilon_2=\sign(\langle \|v\|^{-1}v,x_0\|v\|^{-1}v\rangle)
=\sign(\langle v,x_0v\rangle)$. 
Using (\ref{beta+-}), (\ref{alpha+-}) and the formulas in 
Theorem \ref{rep}, we obtain 
\begin{multline*}
\langle v,x_0v\rangle=\big\{
[2|j|\!+\!3][2|j|\!-\!1][2|j|]^{-1}[2|j|\!+\!1]^{-1}
(\alpha^+ (|j|,|j|)^{\epsilon}_{|j|})^2+\\
[2][2|j|]^{-1}[2|j|\!+\!2]^{-1}
\big\}\beta^0(|j|,|j|)^{\epsilon}_{|j|}.
\end{multline*}
As the expression in braces is positive, we deduce  
$\epsilon_2=\sign(\beta^0(|j|,|j|)^{\epsilon}_{|j|})=\sign(j)$. 

Next,  assume that $j=0$. 
By similar arguments as above, we conclude that $\pi_0\otimes T_{1/2}
=\pi_{1/2}^{\epsilon_1}\oplus\pi_{1/2}^{\epsilon_2}$.
Let $u_1$ and $u_2$ denote the highest weight vectors of the 
spin $1/2$ representations belonging to 
$\pi_{1/2}^{\epsilon_1}$ and $\pi_{1/2}^{\epsilon_2}$, 
respectively. 
Then $\langle u_1,x_0u_1\rangle=\beta^0(1/2,1/2)^{\epsilon_1}_{1/2}$ 
and $\langle u_2,x_0u_2\rangle=\beta^0(1/2,1/2)^{\epsilon_2}_{1/2}$. 
The vector  
$w^{1/2}_{1/2}=v^{0}_{0}\otimes v^{1/2}_{1/2}$ belongs to the span 
of the orthonormal vectors $u_1$ and $u_2$. 
Since $\langle w^{1/2}_{1/2},x_0w^{1/2}_{1/2}\rangle=0$, 
we get $\sign(\beta^0(1/2,1/2)^{\epsilon_1}_{1/2})\neq
\sign(\beta^0(1/2,1/2)^{\epsilon_2}_{1/2})$, 
whence $\epsilon_1\neq\epsilon_2$. 
Summarizing the preceding results, we conclude that, 
for all $j\in\frac{1}{2}\dZ$, the representation 
$\pi_j\otimes T_{1/2}$ decomposes into $\pi_j\otimes T_{1/2}
=\pi_{j-1/2}\oplus\pi_{j+1/2}$.  

The proposition can now be proved by induction. 
Let $k\in\frac{1}{2}\dN$. 
Assume that Proposition \ref{deco} holds for 
all $l=0,1/2,\ldots, k$. By Corollary \ref{sumirr}, 
$\pi_j\otimes T_k\otimes T_{1/2}=
(\pi_j\otimes T_k)\otimes T_{1/2}
=\pi_j\otimes (T_k\otimes T_{1/2})$ decomposes into irreducible 
integrable $\ast$-re\-pre\-sen\-ta\-tions. 
By our induction hypothesis and (\ref{T}), we have 
$$
(\pi_j\otimes T_k)\otimes T_{1/2}=
\pi_{j-k-1/2}\oplus(\oplus_{n=0}^{2k-1}2\pi_{j-k+n+1/2})\oplus\pi_{j+k+1/2}.
$$
On the other hand, 
$$
\pi_j\otimes (T_k\otimes T_{1/2})
=(\oplus_{n=0}^{2k-1}\pi_{j-k+n+1/2})\oplus (\pi_j\otimes T_{k+1/2}). 
$$
Comparing both results shows that 
$$
\phantom{\Box}\hspace{67pt}
\pi_j\otimes T_{k+1/2}=
\pi_{j-k-1/2}\oplus\pi_{j-k+1/2}
\oplus\ldots\oplus\pi_{j+k+1/2}.\hspace{67pt}\Box
$$
%
\subsection{Realization of irreducible integrable 
{\mathversion{bold}$\ast$}-re\-pre\-sen\-ta\-tions of 
{\mathversion{bold}$\cO({\mathrm{S}}^2_{qr})\rti\cU_q({\mathrm{su}}_2)$} 
on {\mathversion{bold}$\cO({\mathrm{SU}}_q(2))$} } 
                                                              \label{real}
%
In this subsection, we relate the irreducible representation $\pi_j$
from Theorem \ref{rep} to 
the representation $\hat{\pi}_j$
corresponding to the 
Hopf module $M_j$ from Section \ref{sec-decomp}.  
%
\begin{thl}                                                   \label{vl}
Let $l\in\frac{1}{2}\dN$. Define $e^+_{1/2} = d+q^{-2l-1}s b$ and 
$e^-_{-1/2} = c-q^{2l+1} s a$. Set $e^+_{-1/2} = F \ang
 e^+_{1/2}$ and  $e^-_{1/2} = E \ang e^-_{-1/2}$. 
Then 
\begin{align}                                              \label{vl+}
q^{1/2} [2l]^{1/2} v^l_{l l} e^+_{-1/2} 
{-} q^{-l} v^l_{l-1,l} e^+_{1/2}= 
q^{-1/2} [2l]^{1/2} v^l_{-l,-l} e^-_{1/2} 
{-} q^{l} v^l_{-l+1,-l} e^-_{-1/2} =0. 
\end{align}
\end{thl}
{\bf Proof.}
Recall that 
$v^l_{ll} = N^l_{ll} u_l = N^l_{ll} (d + q^{-1}  s  b)
\cdot{\dots}\cdot (d+q^{-2l}  s  b)$ 
and that  
$v^l_{l-1,l}=[2l]^{-1/2} F\ang v^l_{ll}=[2l]^{-1/2} N^l_{ll} (F\ang u_l)$. 
A straightforward induction argument shows that 
$F\ang u_l = q^{l-1/2} [2l] u_{l-1/2} (c+q^{-2l}  s  a)$. 
Now a direct calculation gives
\begin{align*}
&u_l e^+_{-1/2} - q^{-l-1/2} [2l]^{-1} (F\ang u_l) e^+_{1/2}\\
& =u_{l-1/2} (d+q^{-2l} s  b) (c+q^{-2l-1}s a) 
- q^{-1} u_{l-1/2} (c + q^{-2l}  s  a)(d+q^{-2l-1}s b)=0
\end{align*}
which implies the first equality in (\ref{vl+}). The second equality is proved 
in the same way by using  $v^{l}_{-l,-l}=\|w_l\|^{-1} w_l$.
\hfill$\Box$
\begin{thp}                                    \label{cor}
Let  $j\!\in\!\frac{1}{2}\dZ$. 
The irreducible integrable
$\ast$-re\-pre\-sen\-ta\-tion $\pi_j$ of  
$\cO(\mathrm{S}^2_{qr})\rti \cU_q(\mathrm{su}_2)$ 
from Theorem \ref{rep} is unitarily equivalent 
to the $\ast$-re\-pre\-sen\-ta\-tion $\hat{\pi}_j$ on 
the left $\cO(\mathrm{S}^2_{qr})\rti \cU_q(\mathrm{su}_2)$-mo\-dule $M_j$ 
from Theorem \ref{T1}. 
\end{thp}
{\bf Proof.} By Theorem \ref{T1}, the restriction of the 
representation $\hat\pi_j$ to $\cU_q(\mathrm{su}_2)$ is of the form
$\oplus_{n\in\dN_0}T_{|j|+n}$. Likewise, by Theorem \ref{rep}, 
the restriction of the representation $\pi_j$ to $\cU_q(\mathrm{su}_2)$ is 
$\oplus_{n\in\dN_0}T_{|j|+n}$. Therefore, by Corollary \ref{sumirr}(i), 
$\hat\pi_j$ is unitarily equivalent either to 
$\pi_{|j|}^-$ or to $\pi_{|j|}^+$. 
That is, it only remains to specify the label $+$ or $-$. Recall from Theorem
\ref{T1} that $M_j$ is the linear span of vectors $v^{|j|+n}_{kj}$, where
$n\in \dN_0$ and $k=-(|j|+n),\dots,|j|+n$.

For $j=0$, there is only the unique Heisenberg representation. 
A direct calculation yields 
$\langle v^{1/2}_{1/2, -1/2}, 
x_0 v^{1/2}_{1/2, -1/2}\rangle 
= \beta^0 (1/2, 1/2)^-_{1/2}$ and 
$\langle v^{1/2}_{1/2, 1/2}, 
x_0 v^{1/2}_{1/2, 1/2}\rangle 
= \beta^0 (1/2, 1/2)^+_{1/2}$
so that $\hat\pi_{-1/2}=\pi_{-1/2}$ and $\hat\pi_{1/2}=\pi_{1/2}$.

We proceed by induction. 
Let $l\in\frac{1}{2}\dN$.  
Assume that $\pi_l$ is unitarily equivalent to $\hat{\pi_l}$. 
Recall that $T_{1/2}$ is the spin $1/2$ representation 
on $V_{1/2}=
\Lin \{v_{-1/2}^{1/2}, v_{1/2}^{1/2}\}$, where $v_{\pm 1/2}^{1/2}$ are the
weight vectors.
By Proposition \ref{deco}, the tensor product representation 
$\hat{\pi}_l\otimes T_{1/2} \cong \pi_l \otimes T_{1/2}$ 
decomposes into the direct 
sum $\pi_{l-1/2}\oplus\pi_{l+1/2}$ on 
$M_l\otimes V_{1/2}= W_{l-1/2}\oplus W_{l+1/2}$ . 
Set $e_{1/2}:=\|d+q^{-2l-1}sb\|^{-1}(d+q^{-2l-1}sb)$ 
and $e_{-1/2}:=F\ang e_{1/2}$.  
Consider the linear mapping 
$\Phi : M_l\otimes V_{1/2}\rightarrow \cO(\mathrm{SU}_q(2))$ 
defined by  $\Phi(v\otimes v_{\pm 1/2}^{1/2})=ve_{\pm 1/2}$. 
Clearly, $\Phi$ intertwines the actions of 
$\cO(\mathrm{S}^2_{qr})\rti \cU_q(\mathrm{su}_2)$ 
on $M_l\otimes V_{1/2}$ and 
$\cO(\mathrm{SU}_q(2))$. From (\ref{cgd}) and the first equality in
(\ref{vl+}), we conclude that $\Phi(v^{l-1/2}_{l-1/2,l-1/2})=0$. Since the
highest weight vector $v^{l-1/2}_{l-1/2,l-1/2}$ is cyclic for the
representation $\pi_{l-1/2}$, the latter implies that 
$\Phi(W_{l-1/2})=\{0\}$. On the other hand, it follows from 
(\ref{uj}) and (\ref{vlkj}) that there is a non-zero constant $\gamma_l$ 
such that 
$\Phi(v^l_{ll}\otimes v_{1/2}^{1/2})=\gamma_l v^{l+1/2}_{l+1/2,l+1/2}$. 
Therefore, $\Phi$ maps $W_{l+1/2}$ into $M_{l+1/2}$, so $\Phi$ is a non-trivial
intertwiner of the irreducible integrable representations $\pi_{l+1/2}$ on 
$W_{l+1/2}$ and $\hat{\pi}_{l+1/2}$ on $M_{l+1/2}$. By a standard
application of Schur's lemma, $\pi_{l+1/2}$ and  $\hat{\pi}_{l+1/2}$ are
unitarily equivalent. This proves the assertion for $l{+}1/2$. 

The proof for negative integers is similar. One replaces 
$v^l_{ll}$ by $v^l_{-l,-l}$ and uses 
$e_{-1/2}:=\|a-q^{2l+1}sc\|^{-1}(a-q^{2l+1}sc)$ 
and the second equality in (\ref{vl+}).   \hfill$\Box$
\mn 

By Proposition \ref{cor} and Theorem \ref{T1}, 
each irreducible integrable 
$\ast$-re\-pre\-sen\-ta\-tion of the cross product algebra 
$\cO(\mathrm{S}^2_{qr})\rti \cU_q(\mathrm{su}_2)$ 
can be realized on $\cO(\mathrm{SU}_q(2))$ and 
$\cO(\mathrm{SU}_q(2))$ decomposes into the orthogonal 
direct sum of all irreducible integrable 
$\ast$-re\-pre\-sen\-ta\-tion of
$\cO(\mathrm{S}^2_{qr})\rti \cU_q(\mathrm{su}_2)$, 
each with multiplicity one. 
%
%
\section{Representations of the cross product 
{\mathversion{bold}$\ast$}-al\-ge\-bras: 
Second approach}
%
                                               \label{sec-2app}
%
\subsection{``Decoupling'' of the cross product algebras} 

                                                        \label{S3}
%
%
Let us first suppose that $r\in[0,\infty)$.
From the relations $AB=q^{-2} BA$ and $AB^\ast = q^2 B^\ast A$, it is 
clear that $\cS=\{A^n\,;\, n\in\dN_0\}$ is a left and right Ore set 
of the algebra $\cO(\mathrm{S}^2_{qr})$. 
Moreover the algebra $\cO(\mathrm{S}^2_{qr})$ has 
no zero divisors. Hence the localization algebra, 
denoted by $\hat{\cO}(\mathrm{S}^2_{qr})$,  
of $\cO(\mathrm{S}^2_{qr})$ at $\cS$ exists. 
The $\ast$-al\-ge\-bra $\cO(\mathrm{S}^2_{qr})$ is then 
a $\ast$-sub\-al\-ge\-bra 
of $\hat{\cO}(\mathrm{S}^2_{qr})$ and all elements $A^n$, $n\!\in\!\dN_0$, 
are invertible in $\hat{\cO}(\mathrm{S}^2_{qr})$. 
From \cite[Theorem 3.4.1]{LR}, 
we conclude that $\hat{\cO} (\mathrm{S}^2_{qr})$ is a left 
$\su$-mo\-dule $\ast$-al\-ge\-bra such that $\cO(\mathrm{S}^2_{qr})$ 
is a left $\su$-mo\-dule $\ast$-sub\-al\-ge\-bra. The left action of the 
generators $E,F,K$ on $A^{-1}$ is given by
$$
E\ang A^{-1} = -q^{-5/2} B^\ast A^{-2},\ \, F\ang A^{-1} = q^{1/2} BA^{-2}, 
\ \,K\ang A^{-1}=A^{-1}.
$$
Hence the left cross product algebra 
$\hat{\cO}(\mathrm{S}^2_{qr})\rti \su$ is a 
well defined $\ast$-al\-ge\-bra containing $\cO(\mathrm{S}^2_{qr})\rti\su$
as a $\ast$-sub\-al\-ge\-bra.

Let $\cY_r$ denote the 
$\ast$-sub\-al\-ge\-bra of $\cO(\mathrm{S}^2_{qr})\rti\su$ 
generated by
\begin{align}\label{xydef}
&X := q^{3/2} \lambda FK^{-1}A + q B = q^{3/2} \lambda AFK^{-1}
 + q^{-1} B, \\
&X^\ast := q^{3/2} \lambda AK^{-1}E + q B^\ast = 
q^{3/2} \lambda  K^{-1}EA+ q^{-1} B^\ast, \\ 
&Y := qK^{-2} A \label{xydef1}
\end{align}
and let $\hat{\cY_r}$ be the subalgebra of 
$\hat{\cO}(\mathrm{S}^2_{qr})$ generated by $\cY_r$ and 
$Y^{-1}$.
Note that $Y=Y^\ast$. It is straightforward to check that 
the elements $X$, $X^\ast$, $Y$, $Y^{-1}$ commute with the generators 
$A$, $B$, $B^\ast$ of $\cO(\mathrm{S}^2_{qr})$. 
Hence the algebras $\hat{\cY_r}$ and 
$\hat{\cO}(\mathrm{S}^2_{qr})$ {\it commute} inside the cross product algebra 
$\hat{\cO} (\mathrm{S}^2_{qr})\rti \su$. Moreover, the generators of 
$\cY_r$ satisfy the commutation relations 
\begin{equation}                                                 \label{xyrel}
YX = q^2 XY,\ \, YX^\ast = q^{-2} X^\ast Y,\ \,
X^\ast X - q^2 XX^\ast = (1-q^2) (Y^2+r).
\end{equation}

We denote by $\hat{\cU}_q(\mathrm{su}_2)$ the Hopf $\ast$-sub\-al\-ge\-bra 
of $\cU_q(\mathrm{su}_2)$
generated by 
$e:= EK$, $f:= K^{-1} F$  and $k:=K^2.$ As an algebra, 
$\hat{\cU}_q(\mathrm{su}_2)$ has generators $e$, $f$, $k$, $k^{-1}$ 
with defining 
relations 
\begin{equation}\label{efrel}
kk^{-1} = k^{-1}k=1,\ \, ke = q^2 ek,\ \,kf = q^{-2} fk,\ \,
ef \!-\! fe = \lambda^{-1} (k\!-\!k^{-1}).
\end{equation}
From (\ref{xydef})--(\ref{xydef1}), we obtain
\begin{equation}\label{efdef}
f = q^{-1/2} \lambda^{-1} (X\!-\!qB) A^{-1},\ \, 
e = q^{1/2} \lambda^{-1} (X^\ast\!-\! q^{-1} B^\ast) Y^{-1},\ \, 
k=qY^{-1} A.
\end{equation}
Hence the two commuting algebras $\hat{\cY_r}$ 
and $\hat{\cO}(\mathrm{S}^2_{qr})$ 
generate the $\ast$-sub\-al\-ge\-bra 
$\hat{\cO}(\mathrm{S}^2_{qr})\rti \hat{\cU}_q (\mathrm{su}_2)$ 
of $\hat{\cO}(\mathrm{S}^2_{qr})\rti \cU_q(\mathrm{su}_2)$. 
In this sense, $\hat{\cY_r}$ and $\hat{\cO}(\mathrm{S}^2_{qr})$ 
 ``decouple'' the cross product algebra 
$\hat{\cO}(\mathrm{S}^2_{qr})\rti \hat{\cU}_q (\mathrm{su}_2)$. 
By choosing a PBW-basis of $\hat{\cO}(\mathrm{S}^2_{qr})$ and $\hat{\cY_r}$, 
it is easy to show that the subalgebra generated 
by  $\hat{\cO}(\mathrm{S}^2_{qr})$
and $\hat{\cY_r}$ is isomorphic to the tensor product algebra 
$\hat{\cO}(\mathrm{S}^2_{qr})\!\otimes\! \hat{\cY_r}$. 
Therefore we can consider 
$\hat{\cO}(\mathrm{S}^2_{qr})\!\otimes\! \hat{\cY_r}$ as 
a $\ast$-sub\-al\-ge\-bra
of $\hat{\cO}(\mathrm{S}^2_{qr})\rti \hat{\cU}_q (\mathrm{su}_2)$ 
by identifying $x\!\otimes\! y$ with $xy$, 
$x\!\in\!\hat{\cO}(\mathrm{S}^2_{qr})$, 
$y\!\in\! \hat{\cY_r}$. Similarly, 
$\cO(\mathrm{S}^2_{qr})\!\otimes\! \cY_r$ becomes a $\ast$-sub\-al\-ge\-bra
of $\cO(\mathrm{S}^2_{qr})\rti \cU_q (\mathrm{su}_2)$.

There is an alternative way to define the $\ast$-al\-ge\-bra 
$\hat{\cO}(\mathrm{S}^2_{qr})\rti \hat{\cU}_q (\mathrm{su}_2)$ 
by taking the two sets $A\!=\!A^\ast$, $A^{-1}$, $B$, $B^\ast$ and 
$X$, $X^\ast$, $Y\!=\!Y^\ast$, $Y^{-1}$ of 
pairwise commuting generators 
with defining relations (\ref{podrel}), (\ref{xyrel}) 
and the obvious relations 
\begin{equation}\label{obrel}
AA^{-1} = A^{-1} A=1, \quad  YY^{-1} = Y^{-1} Y=1.
\end{equation}
Indeed, if we define $e,f,k$ by (\ref{efdef}), then the 
relations (\ref{efrel}) of $\cU_q(\mathrm{su}_2)$ and the cross relations 
of $\hat{\cO}(\mathrm{S}^2_{qr})\rti \hat{\cU}_q(\mathrm{su}_2)$ 
can be derived 
from this set of defining relations. 

The larger cross product algebra 
$\hat{\cO}(\mathrm{S}^2_{qr})\rti \cU_q(\mathrm{su}_2)$ can be redefined in 
a similar manner 
if we replace the generator $Y$ by $K$. That is, 
$\hat{\cO}(\mathrm{S}^2_{qr})\rti \cU_q(\mathrm{su}_2)$ 
is the $\ast$-al\-ge\-bra 
with the generators  $A\!=\!A^\ast$, $A^{-1}$, 
$B$, $B^\ast$, $X$, $X^\ast$, $K\!=\!K^\ast$, 
$K^{-1}$ 
satisfying the defining relations (\ref{podrel}), (\ref{obrel}), 
\begin{equation}
\label{xkrel}
KA\!=\!AK,\, BK\!=\!qKB ,\, KB^\ast \!=\! qB^\ast K,\, 
X^\ast X\!-\!q^2 XX^\ast \!=\! (1\!-\!q^2)(q^2 K^{-4} A^2\!+\!r),
\end{equation}
and $A$, $B$, $B^\ast$ commute with $X$ and $X^\ast$.
The generators $F$ and $E$ are then given by 
\begin{equation}\label{ferel}
F=q^{-3/2} \lambda^{-1} (X-qB) KA^{-1},\quad 
E=q^{-3/2} \lambda^{-1} A^{-1} K(X^\ast-q B^\ast).
\end{equation}

The preceding considerations and facts carry over almost verbatim to 
the case $r=\infty$. The only difference is that in the case $r=\infty$ 
one has to set $r=1$ in   
the third equations of (\ref{xyrel}) and (\ref{xkrel}). 
%
\subsection{Operator representations of 
the {\mathversion{bold}$\ast$}-al\-ge\-bra {\mathversion{bold}$\cY_r$} } 
                                                                   \label{S4}
%
For the study of representations of the $\ast$-al\-ge\-bra $\cY_r$, 
we need two 
auxiliary lemmas. The first one restates the {\it Wold decomposition} of 
an  isometry (see \cite[Theorem 1.1]{SF}), while the second 
is Lemma 4.2(ii) in \cite{SW}.

\begin{thl}                                                     \label{L1}
Each isometry $v$ on a Hilbert space $\cK$ is up to unitary equivalence 
of the following form: There exist Hilbert subspaces $\cK^u$ and $\cK^s_0$ 
of $\cK$ and a unitary operator $v_u$ on $\cK^u$ such that $v=v_u\oplus v_s$ 
on $\cK = \cK^u \oplus \cK^s$ and $v_s$ acts on 
$\cK^s = \bar\oplus^\infty_{n=0} \cK^s_n$ by 
$v_s \zeta_n = \zeta_{n+1}$, where each $\cK^s_n$ is $\cK^s_0$ and 
$\zeta_n\in\cK^s_n$. 
Moreover, $\cK^u = \cap^\infty_{n=0} v^n \cK$.
\end{thl}

\begin{thl}                                                     \label{L2}
Let $v$ be the operator $v_s$ on $\cK^s = \bar\oplus^\infty_{n=0} \cK^s_n$, 
$\cK^s_n= \cK^s_0$, from  Lemma \ref{L1} 
and let $Y$ be a self-adjoint operator 
on $\cK^s$ such that $q^2 v Y \subseteq Y v$. Then there is a self-adjoint 
operator $Y_0$ on the Hilbert space  $\cK^s_0$ such that 
$Y\zeta_n=q^{2n} Y_0 \zeta_n$ and $v\zeta_n=\zeta_{n+1}$ for 
$\zeta_n\in \cK^s_n$ and $n\in\dN_0$.
\end{thl}

It suffices to treat the case $r\in [0,\infty)$ because the  
algebra $\cY_\infty$ is isomorphic to $\cY_1$. 
Suppose that we have a representation of relations (\ref{xyrel}) by 
closed operators $X$, $X^\ast$ and a self-adjoint operator $Y$ 
acting on a Hilbert space $\cK$. We assume that $Y$ has trivial 
kernel.

Let $X= v|X|$ be the polar decomposition of the operator $X$. 
Since $q\in (0,1)$, $r\ge 0$ and $\ker Y=\{0\}$, we have $\ker X= \{0\}$ by 
the third equation of (\ref{xyrel}). Hence $v$ is an isometry on $\cK$, 
that is, $v^\ast v=I$.

By (\ref{xyrel}), $YX^\ast X= X^\ast XY$. We assume that the 
self-adjoint operators $Y$ and $X^\ast X$ strongly commute. Then $Y$ 
and $|X|=(X^\ast X)^{1/2}$ also strongly commute. 
Again by (\ref{xyrel}), 
$Y v|X|=YX = q^2 XY=q^2 v|X|Y =q^2 v Y|X|$. Since $\ker|X|=\ker X=\{0\}$, 
we conclude that
\begin{equation}\label{vyrel}
q^2 v Y=Yv.
\end{equation}
Since $X^\ast = |X| v^\ast$, the third equation of (\ref{xyrel}) rewrites as 
\begin{equation}\label{xvrel}
|X|^2 = q^2 v|X|^2 v^\ast + (1\!-\!q^2)(Y^2 +r).
\end{equation}
Multiplying by $v$ gives 
$|X|^2 v= v(q^2 |X|^2 + (1\!-\!q^2)(q^4 Y^2 + r))$ since 
$v^\ast v=I$ and, by (\ref{vyrel}),
$Y^2v=q^4v Y^2$.
Proceeding by induction, 
we derive
\begin{equation}\label{xvnrel}
|X|^2 v^n = v^n \big(q^{2n} |X|^2 + (1\!-\!q^{2n})(q^{2n+2} Y^2 + r)\big).
\end{equation}
Now we use the Wold decomposition $v=v_u\oplus v_s$ on 
$\cK=\cK^u\oplus \cK^s$ of the isometry $v$ by Lemma \ref{L1}. 
Since $\cK^u=\cap^\infty_{n=0} v^n \cK$ and $Yv=q^2 vY$, 
the Hilbert subspace $\cK^u$ reduces the self-adjoint operator $Y$. 
From (\ref{xvnrel}), it follows that $|X|^2$ leaves a dense subspace 
of the space $\cK^u$ invariant. We assume that $\cK^u$ reduces $|X|^2$. 
Hence we can consider all operators occurring in  relation (\ref{xvnrel}) 
on the subspace $\cK^u$. For the unitary part $v_u$ of $v$, we have 
$v_uv^\ast_u=I$ on $\cK^u$. Multiplying (\ref{xvnrel}) by 
$(v^\ast_u)^n$ and using again Equation (\ref{vyrel}), we derive on $\cK^u$
the relation 
$$
0\le q^{2n} v^n_u |X|^2 (v^\ast_u)^n 
= |X|^2 - (1\!-\!q^{2n})r + q^2 Y^2 - q^{-2n+2} Y^2
$$
for all $n\in \dN$. Letting $n\rightarrow\infty$ and remembering that  
$\ker Y=\{0\}$ and $q\in (0,1)$, we conclude that the latter is only possible 
when $\cK^u = \{0\}$. That is, the isometry $v$ 
is (unitarily equivalent to ) the unilateral shift operator $v_s$.

Since $v=v_s$, Lemma \ref{L2} applies to the relation $q^2 v Y=Yv$. 
Hence there exists a self-adjoint operator $Y_0$ on $\cK^s_0$ 
such that $Y\zeta_n = q^{2n} Y_0 \zeta_n$ 
on $\cK^s=\bar\oplus^\infty_{n=0} \cK^s_n$, $\cK^s_n = \cK^s_0$. 
Since $v^\ast\zeta_0=0$, (\ref{xvrel}) yields 
$|X|^2\zeta_0= (1-q^2)(Y^2_0 +r)\zeta_0$. 
Using this equation and (\ref{xvrel}), we compute 
\begin{align}\label{x2rel}
|X|^2\zeta_n&=|X|^2 v^n\zeta_0 
= v^n \big(q^{2n} |X|^2 + (1\!-\!q^{2n})(q^{2n+2} Y^2+r)\big)\zeta_0\nonumber\\
&= v^n \big(q^{2n} (1\!-\!q^2)(Y^2_0+r)\zeta_0 
+ (1\!-\!q^{2n})(q^{2n+2} Y^2_0 +r)\zeta_0\big)\nonumber\\
&=\lambda_{n+1}^2 (q^{2n} Y^2_0 +r)\zeta_n
\end{align}
for $\zeta_n\in\cK^s_n$. We assumed above that the pointwise commuting 
self-adjoint operator $|X|^2$ and $Y$ strongly commute. Hence their 
reduced parts on each subspace $\cK^s_n$ strongly commute. 
Therefore, (\ref{x2rel}) implies 
$|X| \zeta_n=\lambda_{n+1} (q^{2n} Y^2_0 + r)^{1/2}\zeta_n$. 
Recall that $X=v |X|$ and $X^\ast = |X| v^\ast$.  
Renaming $\cK^s_n$ by $\cK_n$ and summarizing the preceding, 
we obtain the following form of the operators $X$, $X^\ast$ and $Y$:
\begin{align*}
&X \zeta_n = \lambda_{n+1} (q^{2n} Y^2_0 + r)^{1/2} \zeta_{n+1},\ \;
X^\ast \zeta_n = \lambda_n (q^{2n-2} Y^2_0 + r)^{1/2} \zeta_{n-1},\\
&Y \zeta_n = q^{2n} Y_0 \zeta_n\quad \mbox{on}\ \,
\cK = \bar\oplus^\infty_{n=0} \cK_n,\ \, \cK_n=\cK_0,
\end{align*}
where $Y_0$ is self-adjoint operator with trivial kernel 
on the Hilbert space $\cK_0$. Conversely, it is easy to check 
that these operators $X,X^\ast, Y$ satisfy the relations (\ref{xyrel}), 
so they define indeed a $\ast$-re\-pre\-sen\-ta\-tion 
of the $\ast$-al\-ge\-bra $\cY_r$. 
The representation is irreducible if and only if $\cK_0=\dC$. 
In this case, $Y_0$ is a non-zero real number.  
 
%
\subsection{Representations of the cross product 
{\mathversion{bold}$\ast$}-al\-ge\-bras}
\label{sec-4}
First let us review the representations of the $\ast$-al\-ge\-bra 
$\cO (\mathrm{S}^2_{qr})$ from \cite{P}. 
Recall the definition of $\lambda_\pm$ from Section \ref{sec-decomp}. 
For $r<\infty$, we set 
$$
c_\pm (n) := (r+ \lambda_\pm q^{2n} - (\lambda_\pm q^{2n})^2)^{1/2}.
$$
Let $\Hh^0,\Hh^+_0$ and $\Hh^-_0$ be Hilbert spaces and let $u$ be 
a unitary operator on 
$\Hh^0$. Let $\Hh^\pm = \bar \oplus^\infty_{n=0}\Hh^\pm_n$, 
where $\Hh^\pm_n = \Hh^\pm_0$, and let 
$\Hh = \Hh^0\oplus\Hh^+ \oplus\Hh^-$ for $r\!\in\! (0,\infty]$ and  
$\Hh = \Hh^0 \oplus \Hh^+$ for $r=0$. 
The generators of $\cO(\mathrm{S}^2_{qr})$ act 
on the Hilbert space $\Hh$ by the following formulas:
\begin{align*}
&r \in [0,\infty):\quad 
A=0,\  B=r^{1/2} u,\  B^\ast = r^{1/2}\  u^\ast \ \ \mbox{ on}~ \Hh^0,\\
&A \eta_n = \lambda_\pm q^{2n} \eta_n, \ B \eta_n = c_\pm (n) \eta_{n-1},\  
B^\ast \eta_n = c_\pm (n\!+\!1) \eta_{n+1} \ \ \mbox{on}~ \Hh^\pm.\\[6pt]
&r = \infty: \quad 
 A = 0, \ B=u, \   B^\ast = u^\ast \ \ \mbox{on}~ \Hh_0,\\
&A \eta_n = \pm q^{2n} \eta_n, \ B \eta_n = (1\!-\!q^{4n})^{1/2} \eta_{n-1},\  
B^\ast \eta_n = (1\!-\!q^{4(n+1)})^{1/2} \eta_{n+1}\ \  \mbox{on}~ \Hh^\pm.
\end{align*}

Recall that, for $r=0$, there is no Hilbert space $\Hh^-$.
From  \cite[Proposition 4]{P}, it follows that, up to unitary equivalence, 
each $\ast$-re\-pre\-sen\-ta\-tion of $\cO(\mathrm{S}^2_{qr})$ 
is of the above form. 
Note that all operators are bounded and $\Hh^0=\ker A$, $A>0$ on $\Hh^+$ 
and $A<0$ on $\Hh^-$.

Next we show that, for each $\ast$-rep\-re\-sen\-ta\-tion 
of the algebras $\cO(\mathrm{S}^2_{qr}) \rti \cU_q (\mathrm{su}_2)$ and       
$\cO(\mathrm{S}^2_{qr})\rti \hat{\cU}_q(\mathrm{su}_2)$,  
the space $\Hh^0$ is $\{0\}$, 
so the operator $A$ is invertible. We carry out the reasoning 
for $\cO(\mathrm{S}^2_{qr})\rti \cU_q(\mathrm{su}_2)$. 
Assuming that the commuting 
self-adjoint operators $A$ and $K$ strongly commute, it follows 
that $K$ leaves $\Hh^0 = \ker A$ invariant. 
By (\ref{podrel}) (resp.\ (\ref{prodrel1})), 
$B\Hh^0\subseteq  \Hh^0$. Let $\xi \in \Hh^0$. Using the relation 
$FA = AF \!-\!q^{-3/2} BK$, we see that 
$$
q^{-3} \| BK \xi\|^2= \|AF \xi\|^2 = \langle q^{-3/2} BK \xi, 
AF\xi\rangle = \langle q^{-3/2} ABK \xi, F\xi\rangle =0.
$$
That is, $BK \xi = 0$ and $AF \xi =0$, so that $F\xi \in\Hh^0$. 
For $r\in (0,\infty]$, $BK$ is invertible on $\Hh^0$ and hence $\xi=0$. 
For $r=0$, we have $B^\ast \xi=B^\ast F\xi =0$. 
From the cross relation 
$FB^\ast = q^{-1} B^\ast F +q^{-1/2} (1\!+\!q^2)AK -q^{-1/2} K$,  
we get $K \xi=0$ and so $\xi=0$. Thus, $\Hh^0=\{0\}$.

Consider now a $\ast$-rep\-re\-sen\-ta\-tion of the $\ast$-al\-ge\-bra 
$\hat{\cO}(\mathrm{S}^2_{qr})\rti 
\hat{\cU}_q(\mathrm{su}_2)$.  Its restriction to 
$\cO(\mathrm{S}^2_{qr})$ is a 
$\ast$-rep\-re\-sen\-ta\-tion of the form described above with  
$\Hh^0 = \{0\}$. 
All operators of the $\ast$-sub\-al\-ge\-bra $\cY_r$ commute with $A$ and $B$. 
Let us assume that the spectral projections of the self-adjoint 
operator $A$ commute also with all operators of $\cY_r$. 
(Note that, for a $\ast$-rep\-re\-sen\-ta\-tion by bounded 
operators on a Hilbert space,  
the latter fact follows. For unbounded $\ast$-rep\-re\-sen\-ta\-tion,  
it does not 
and we restrict ourselves to the class of well behaved representations 
which satisfy this assumption.) Since $\Hh^\pm_n$ is the eigenspace 
of $A$ with eigenvalue $\lambda_\pm q^{2n}$ and these eigenvalues 
are pairwise distinct,  the operators of $\cY_r$ leave $\Hh^\pm_n$ 
invariant. That is, 
we have a $\ast$-rep\-re\-sen\-ta\-tion of $\cY_r$ on $\Hh^\pm_n$. 
But $\cY_r$ commutes also with $B$. Since $B$ is a weighted shift 
operators with weights $c_\pm (n)\ne 0$ for $n\in\dN$, it follows 
that the representations of $\cY_r$ on $\Hh^\pm_n \!=\! \Hh^\pm_0$ are 
the same for all $n\in\dN_0$. Using the structure of the representation 
of the $\ast$-al\-ge\-bra $\cY_r$ on the Hilbert space $\cK := \Hh^\pm_0$ 
derived in Subsection \ref{S4}
and inserting the formulas 
for $X,X^\ast,Y$ and $A,B,B^\ast$ into (\ref{efdef}), 
one obtains the action of the generators 
$e,f,k,k^{-1}$ of $\hat{\cU}_q(\mathrm{su}_2)$. 
We do not list these formulas, but we will do so below for the generators 
of $\cU_q(\mathrm{su}_2)$. 

Let us turn to a $\ast$-rep\-re\-sen\-ta\-tion of the larger 
$\ast$-al\-ge\-bra 
$\hat{\cO}(\mathrm{S}^2_{qr})\rti \cU_q(\mathrm{su}_2)$. Then, in addition 
to the considerations of the preceding paragraph, we have to deal 
with the generator $K$. We assumed above that $K$ and $A$ are strongly 
commuting self-adjoint operators. Therefore, $K$ commutes with the spectral 
projections of $A$. Hence each space $\Hh^\pm_n$ is reducing for $K$. 
The  relation $BK=q KB$ implies that there is an invertible self-adjoint 
operator $K_0$ on $\Hh^\pm_0$ such that  $K \eta_n = q^n K_0\eta_n$ 
for $\eta_n\in \Hh^\pm_n$. Recall that we have $XK= q KX$ and $YK=KY$ in 
the algebra $\hat{\cO}(\mathrm{S}^2_{qr})\rti \cU_q(\mathrm{su}_2)$. 
Inserting for $X$ and $Y$ 
the corresponding operators on 
$\cK\equiv \Hh^\pm_0 = \oplus^\infty_{m=0} \cK_m$ from Subsection \ref{S4},  
we conclude that there exists an invertible self-adjoint 
operator $H$ on $\cK_0$ such that $K_0\zeta_m =q^{-m} H\zeta_m$ for 
$\zeta_m\in\cK_m$. Further, since $Y=q K^{-2} A$, 
we have $Y \zeta_0=Y_0\zeta_0=q H^{-2} \lambda_\pm \zeta_0$ 
for $\zeta_0\in\cK_0$.
Inserting the preceding facts and the results from Subsection \ref{S4} 
into (\ref{ferel})
and renaming $\cK_0$ by $\cG$, 
we obtain the following $\ast$-re\-pre\-sen\-ta\-tions (satisfying the
assumptions made above) 
of the cross product 
$\ast$-al\-ge\-bra $\hat{\cO}(\mathrm{S}^2_{qr})\rti \cU_q(\mathrm{su}_2)$ 
for $r\in(0,\infty)$:
\begin{align*}
(I)_{\pm,H}:\quad &A\eta_{nm} \hsp=\hsp  \lambda_\pm q^{2n} \eta_{nm}, 
\hsp \ \,
B\eta_{nm} \hsp=\hsp  c_\pm (n) \eta_{n-1,m}, \hsp \ \,
B^\ast \eta_{nm} \hsp=\hsp  c_\pm (n\!+\!1) \eta_{n+1,m},\\
&E \eta_{nm} \hsp=\hsp  q^{-1/2} \lambda^{-1} [q^{-n} \lambda_m 
(\lambda_\pm^{-2} q^{-2m} r + H^{-4})^{1/2} H \eta_{n,m-1}\\
&\qquad \quad -q^{-m} (\lambda^{-2}_\pm q^{-2n-2} r 
+ \lambda^{-1}_\pm -q^{2n+2})^{1/2} H \eta_{n+1,m}],\\
&F \eta_{nm} \hsp=\hsp  q^{-1/2} \lambda^{-1} [q^{-n} \lambda_{m+1} 
(\lambda^{-2}_\pm q^{-2m-2} r + H^{-4})^{1/2} H \eta_{n,m+1}\\
&\qquad \quad -q^{-m} (\lambda^{-2}_\pm q^{-2n} r 
+ \lambda^{-1}_\pm - q^{2n})^{1/2} H \eta_{n-1,m}],\\
&K \eta_{nm} \hsp=\hsp q^{n-m} H \eta_{nm},
\end{align*}
where $H$ is an invertible self-adjoint operator on a Hilbert space $\cG$. 
The domain is the direct 
sum $\Hh=\oplus^\infty_{n,m=0} \Hh_{nm}$, where $\Hh_{nm} =\cG$.
In the case $r=0$, there is only the representation $(I)_{+,H}$. 
The case $r=\infty$ has already been treated in \cite{SW}.
%
%
\section{Algebras of functions and representations 
of the cross product algebras}
                                                             \label{***}
%
%
\subsection{Algebras of functions and the invariant state 
on {\mathversion{bold}${\mathrm{S}}^2_{qr}$} }
                                                      \label{sec-algfun}
Let $\cF(\sigma(A))$ be the $\ast$-al\-ge\-bra of all 
complex Borel functions 
on the set 
$\sigma(A):=\{q^{2n} \lambda_+, q^{2n} \lambda_-\,;\, n\in\dN_0\}$, 
let $\cF_\mathrm{b}(\sigma(A))$ be the $\ast$-al\-ge\-bra of all 
bounded complex Borel functions on $\sigma(A)$, and let 
$\cF^\infty(\sigma(A))$ be the set of all $f\in \cF(\sigma(A))$ 
for which there exist 
an $\varepsilon > 0$ and a function 
$\tilde{f}\in C^\infty (-\varepsilon,\varepsilon)$ 
such that $f=\tilde{f}$ on $\sigma(A)\cap (-\varepsilon,\varepsilon)$.   
For $f\in \cF^\infty(\sigma(A))$, we can assign unambiguously 
the value $f(0)$ at $0$ by taking the limit 
$f(0)=\lim_{t\rightarrow 0}f(t)$. 
In order to be in accordance with our previous notation, we write A 
for the function $\mathrm{id}(t)= t$ and also for the argument of functions 
$f\in \cF(\sigma(A))$. We denote by 
$\cF(\mathrm{S}^2_{qr})$, $r\in[0,\infty]$, 
the unital $\ast$-al\-ge\-bra generated by the $\ast$-al\-ge\-bra 
$\cF(\sigma(A))$ and two generators $B$, $B^\ast$ with defining relations
\begin{align}\label{bfrel}
&B f(A) = f(q^2 A) B,\quad f(A) B^\ast = B^\ast f(q^2 A),\quad 
f\in \cF (\sigma(A));\\
\label{barel1}
&B^\ast B = A- A^2 + r,\quad BB^\ast = q^2 A - q^4 A^2 + r,\quad
\mbox{for}~~ r\in[0,\infty),\\
\label{barel2}
&B^\ast B = -A^2 + 1,\quad BB^\ast = -q^4 A^2+ 1 \quad\mbox{for}~~ r=\infty.
\end{align}
For $k\!\in\! \dZ$, we set 
$B^{\# k}\!=\!B^k$ if $k\!\ge\! 0$ 
and $B^{\# k} \!=\! B^{\ast k}$ if $k\!<\!0$. 
As a vector space, $\cF(\mathrm{S}^2_{qr})$ is spanned by 
$\{B^n f_1(A),\, f_2 (A) B^{\ast k}\,;\, f_1, f_2 \!\in\! \cF(\sigma(A)),\  
n,k \!\in\! \dN_0\}
=\{f(A) B^{\# k}\,;\, f\in \cF(\sigma(A)),\ k \in \dZ\}$.
We denote by  $\cF_{\mathrm{b}}(\mathrm{S}^2_{qr})$ and 
$\cF^\infty(\mathrm{S}^2_{qr})$ the $\ast$-subalgebras 
of $\cF(\mathrm{S}^2_{qr})$ generated by 
$\cF_{\mathrm{b}}(\sigma(A))$ and
$\cF^\infty(\sigma(A))$, respectively, and $B$ and $B^\ast$. 
We introduce a left action~$\ang$ of the Hopf algebra 
$\cU_q(\mathrm{su}_2)$ on $\cF^\infty(\mathrm{S}^2_{qr})$
by setting 
\begin{align*}
&E\ang p(B) f(A) 
=q^{1/2} \left[\frac{p(q^{-1} B) - p(q B)}{(q^{-1}-q)B}\,-\, 
(1\!+\!q^2) ~ \frac{p(q^{-3} B) - p(q B)}{(q^{-3} -q) B}~A\right]f(A)\\
&\qquad\hspace{1.4cm}\quad 
+q^{-1/2} p(qB)B^\ast~~\frac{f(A)-f(q^2A)}{(1-q^2)A}~,\\[-24pt]
\end{align*}
\begin{align*}
&E\ang f(A)p(B^\ast) =
q^{-1/2}  \frac{f(q^{-2} A) - f(A)}{(q^{-2} -1) A}~B^\ast p(q B^\ast),
\hspace{140pt}\\[-24pt]
\end{align*}
\begin{align*}
&F\ang p(B) f(A) =- q^{-3/2} Bp(qB) \frac{f(q^{-2} A) - f(A)}{(q^{-2}-1)A},
\hspace{148pt}\\[-24pt]
\end{align*}
\begin{align*}
&F\ang f(A) p(B^\ast) = 
-q^{-3/2} ~\frac{f(A)-f(q^2A)}{(1-q^2)A} Bp(qB^\ast)\\
&\quad\qquad -q^{-1/2} f(A)  \left[\frac{p(q^{-1} B^\ast) 
- p(q B^\ast)}{(q^{-1}-q)B^\ast}~-~ (1\!+\!q^2) A~ \frac{p(q^{-3} B^\ast) 
- p(q B^\ast)}{(q^{-3} -q) B^\ast} \right],\\[-24pt]
\end{align*}
\begin{align*}
&K\ang p(B) f(A) = p(q^{-1} B) f(A), \quad 
K\ang f(A) p (B^\ast) = f(A) p (q B^\ast),\hspace{76pt}
\end{align*}
where $f$ is a function in $\cF^\infty(\sigma(A))$
and $p$ is a polynomial. 
It can be shown that these formulas define indeed an action of the 
Hopf $\ast$-al\-ge\-bra $\cU_q(\mathrm{su}_2)$ such that 
$\cF^\infty(\mathrm{S}^2_{qr})$ is a left 
$\cU_q(\mathrm{su}_2)$-mo\-dule $\ast$-al\-ge\-bra. 
We omit the details of these lengthy, 
but straightforward computations.

From the defining relations, it is clear that the coordinate 
$\ast$-al\-ge\-bra $\cO(\mathrm{S}^2_{qr})$ of the quantum sphere is a 
$\ast$-sub\-al\-ge\-bra of $\cF^\infty(\mathrm{S}^2_{qr})$. 
On $B$, $A$, $B^\ast$ 
considered as elements of $\cF^\infty(\mathrm{S}^2_{qr})$, 
the preceding formulas coincide with corresponding 
formulas for the actions on the generators $B$, $A$, $B^\ast$ 
of $\cO(\mathrm{S}^2_{qr})$. 
Hence $\cO(\mathrm{S}^2_{qr})$ is 
a left $\cU_q(\mathrm{su}_2)$-mo\-dule $\ast$-sub\-al\-ge\-bra 
of $\cF^\infty(\mathrm{S}^2_{qr})$.

Now we turn to the construction of 
a $\cU_q(\mathrm{su}_2)$-in\-var\-iant linear 
functional $h$. 
For $f\in\cF_{\mathrm{b}}(\sigma(A))$, we put
$$
h_0(f(A)):=\gamma_+ \sumop^\infty_{n=0} f(\lambda_+ q^{2n}) q^{2n} 
+\gamma_- \sumop^\infty_{n=0} f(\lambda_- q^{2n})q^{2n},
$$
where $\gamma_+ := (1\!-\!q^2) \lambda_+ (\lambda_+ \!-\! \lambda_-)^{-1}$ 
and   
$\gamma_- := (1-q^2)\lambda_- (\lambda_- \!-\!\lambda_+)^{-1}$.
Note that $\gamma_\pm = (1-q^2)(1/2 \pm (r+4)^{-1/2})$. 
(When $r=0$, the above equation simplifies to 
$h_0(f(A))=(1-q^2) \sum^\infty_{n=0} f(q^{2n}) q^{2n}$ since in this case 
$\gamma_+=\lambda_+=1$ and $\gamma_-=0$.) 
Define a functional $h$ on $\cF_{\mathrm{b}}(\mathrm{S}^2_{qr})$ by 
\begin{equation}\label{hadef}
h(p(B) f(A)) = p(0) h_0 (f(A)),\quad 
h(f(A) p(B^\ast)) = p(0) h_0 (f(A)).
\end{equation}
For $g\in\cF_{\mathrm{b}} (\sigma(A))$, we have 
\begin{equation}  \label{hqf1}
h(g(A)) =q^2 h(g(q^2 A)) + \gamma_+ g(\lambda_+) 
+ \gamma_- g(\lambda_-).
\end{equation}
First we show that $h$ is a faithful state on the 
$\ast$-al\-ge\-bra $\cF_{\mathrm{b}}(\mathrm{S}^2_{qr})$. 
We restrict ourselves to the case 
$r\in(0,\infty]$. 
Let $x=\sum_k(B^k f_k(A) + g_k(A) B^{\ast k}) \in 
\cF_{\mathrm{b}}(\mathrm{S}^2_{qr})$ 
with $f_k, g_k\in \cF_{\mathrm{b}}(\sigma(A))$. 
From (\ref{bfrel})--(\ref{barel2}), we obtain 
$B^{\ast k} B^k = \Pi^{k-1}_{i=0} p_r(q^{-2i} A)$ and 
$B^k B^{\ast k}= \Pi^k_{i=1} p_r(q^{2i} A)$, where  
$p_r(A) := (\lambda_+-A)(A-\lambda_-)$. By (\ref{hadef}), we have 
\begin{align}\label{hxx}
h(x^\ast x) &= \sumop_k h(B^k\overline{g_k} g_k B^{\ast k} 
+ \overline{f_k} B^{\ast k} B^k f_k)\nonumber\\
&=\sumop_k h_0 (|g_k( q^{2k} A)|^2 
\mathop{\mbox{$\prod$}}^k_{i=1} p_r(q^{2i} A)) 
+ h_0 (|f_k(A)|^2 \mathop{\mbox{$\prod$}}^{k-1}_{i=0} p_r(q^{-2i} A))
\end{align}
Note that $p_r(\lambda_\pm)=0$ and $p_r(q^{2j} \lambda_\pm)>0$ for $j\in\dN$. 
Hence $h(x^\ast x)\ge 0$ by (\ref{hxx}). Assume that $h(x^\ast x)=0$. 
From (\ref{hxx}) and the definition of $h_0$, it follows that 
$g_k(q^{2n} \lambda_\pm)=0$ for $n\in\dN_0$ and $f_k(q^{2n} \lambda_\pm)=0$ 
for $n\ge k$. The latter implies $g_k(A)=0$ and $B^k f_k(A)=0$, 
so $x=0$. That is, $h$ is a faithful state.

Now we prove that $h$ is $\cU_q(\mathrm{su}_2)$-in\-var\-iant 
on the $\ast$-sub\-al\-ge\-bra $\cF^\infty(\mathrm{S}^2_{qr})$, that is, 
$h(y\ang x)=\varepsilon (y) h(x)$ for $y\in \cU_q(\mathrm{su}_2)$ and 
$x\in\cF^\infty(\mathrm{S}^2_{qr})$. 
Clearly, it suffices to verify the latter condition 
for the generators $y=E,F,K$. For $y=K$, the assertion is obvious. 
Since the functional $h$ is hermitian and $\cF^\infty(\mathrm{S}^2_{qr})$ is a 
$\cU_q(\mathrm{su}_2)$-mo\-dule $\ast$-al\-ge\-bra, it is sufficient 
to check the 
invariance with respect to $E$. 
By the definitions of the action of $E$ and the functional $h$, it is 
enough to show that $h(E\ang Bf(A))=0$ for all 
$f(A)\in\cF^\infty(\mathrm{S}^2_{qr})$. 
We carry out the proof for $r\in(0,\infty)$. 
Since $f\in\cF^\infty(\mathrm{S}^2_{qr})$, the function 
$g(A) := f(A)-Af(A) + r(f(A)-f(0))/A$ is bounded on 
$\sigma(A)$. 
Hence $g\in\cF_{\mathrm{b}}(\sigma(A))$ and (\ref{hqf1}) 
applies to this function. 
Observe that $1-\lambda_\pm + r\lambda^{-1}_\pm=0$ and 
$\gamma_+\lambda^{-1}_+ +\gamma_- \lambda^{-1}_-=0$ by the definitions 
of these constants. Using these facts and the  
equation $BB^\ast=q^2 A-q^4A^2+r$,  
we compute
\begin{align*}
q^{-1/2}E\ang Bf(A)
&=f(A)-(1+q^2)Af(A)+ BB^\ast ~\frac{f(A)-f(q^2A)}{(1-q^2)A}\\
&= (1-q^2)^{-1} (g(A) - g(q^2 A) q^2)
\end{align*}
and hence
\begin{align*}
&q^{-1/2} (1{-}q^2) h(E\ang Bf(A))=\gamma_+g(\lambda_+) 
+ \gamma_- g(\lambda_-)=\\
&\gamma_+ f(\lambda_+) (1{-}\lambda_+ {+} r\lambda_+^{-1} ) 
+ \gamma_- f(\lambda_-)(1{-}\lambda_- {+} r \lambda_-^{-1}) - f(0) 
(\gamma_+ \lambda_+^{-1} {+} \gamma_- \lambda_-^{-1})
=0.
\end{align*}
Thus, $h$ is a $\cU_q(\mathrm{su}_2)$-in\-var\-iant state 
on $\cF^\infty(\mathrm{S}^2_{qr})$. 

For the coordinate algebra $\cO(\mathrm{S}^2_{qr})$,  
the preceding description of the invariant functional was obtained 
in \cite{NM}. 
However for our consideration  below it 
is crucial to have the invariant state 
on the larger $\ast$-al\-ge\-bra $\cF^\infty(\mathrm{S}^2_{qr})$. 

The preceding proof shows that, with the action of $E$ defined by the above
formula,  $h(E\ang Bf(A)) = 0$ 
for $f\in\cF_{\mathrm{b}}(\sigma(A))$ if $f(0):=\lim_{t\rightarrow 0}f(t)$ 
exists and 
the function $(f(A) {-} f(0))/A$ is bounded on $\sigma(A)$. 
For instance, $h(E\ang BA\chi_+(A))=0$, where $\chi_+$ is the
characteristic function of $[0,\infty)$. But we get 
$h(E\ang B\chi_+(A))=q^{1/2} (1-q^2)^{-1} 
\gamma_+ \lambda_-\ne 0$ for $r\in (0,\infty]$,
so the $\cU_q(\mathrm{su}_2)$-invariance 
of the functional $h$ does not hold 
on the larger $\ast$-al\-ge\-bra $\cF_{\mathrm{b}}(\mathrm{S}^2_{qr})$.

Now we develop a second operator-theoretic approach to 
the $\cU_q(\mathrm{su}_2)$-mo\-dule 
structure of $\cF^\infty(\mathrm{S}^2_{qr})$ and to the invariant state $h$. 
Suppose that $\pi$ is a $\ast$-re\-pre\-sen\-ta\-tion 
of the $\ast$-al\-ge\-bra $\cF_{\mathrm{b}}(\mathrm{S}^2_{qr})$ 
on a Hilbert space $\Hh$ such that $\ker \pi(A)=\{0\}$. Then all operators 
$\pi(x)$, $x\in\cF_{\mathrm{b}}(\mathrm{S}^2_{qr})$, 
are bounded and leave the dense domain 
$\cD := \cap^{\infty}_{n=1} \cD (\pi(A)^{-n})$ of $\Hh$ invariant. 
For notational 
simplicity, we write $x$ for the operators $\pi(x)$ and $\pi(x)\lceil\cD$. 
From the form of the representations of $\cO(\mathrm{S}^2_{qr})$ described 
in Subsection \ref{sec-4}, 
it is clear that $\sign A:=A|A|^{-1}$ commutes with 
all representation operators. 
For $T\in \cL^+(\cD)$, we define
\begin{align} 
E\ang T &=-q^{1/2} \lambda^{-1} A^{-1} 
\big[B^\ast, |A|^{1/2} T|A|^{-1/2} \big]
\nonumber  \\
&=-q^{-1/2} \lambda^{-1} \sign A \, |A|^{-1/2} \big[B^\ast,T\big] |A|^{-1/2},
                                                   \label{eactt1}\\[-24pt]
\nonumber
\end{align}
\begin{align}
\nonumber
F\ang T &=-q^{-5/2} \lambda^{-1} A^{-1} \big[B, |A|^{1/2} T|A|^{-1/2} \big] \\
&=-q^{-3/2} \lambda^{-1} \sign A \,|A|^{-1/2} \big[B,T\big] |A|^{-1/2},  
                                              \label{eactt2}\\[-24pt]
\nonumber
\end{align}
\begin{align}
                                                         \label{eactt3}
K\ang T &= |A|^{1/2} T|A|^{-1/2},\quad  K^{-1}\ang T = |A|^{-1/2} T |A|^{1/2}.
\end{align}
Then one verifies that the $\ast$-al\-ge\-bra $\cL^+(\cD)$ is a 
$\cU_q(\mathrm{su}_2)$-mo\-dule $\ast$-al\-ge\-bra 
with the actions of the generators 
$E$, $F$, $K$, $K^{-1}$ given by these formulas and
$\cF^\infty(\mathrm{S}^2_{qr})$ is 
a $\cU_q(\mathrm{su}_2)$-mo\-dule $\ast$-sub\-al\-ge\-bra of $\cL^+(\cD)$. 
Setting $T=p(B)f(A)$ or $T=f(A) p(B^\ast)$, one easily obtains the 
formulas for the actions of $E$, $F$, $K$, $K^{-1}$ 
on $\cF^\infty(\mathrm{S}^2_{qr})$ 
listed above.

In order to define the invariant state $h$,  
we specialize the representation $\pi$. 
For $r\in (0,\infty]$, let $\pi$ be the direct sum representation on 
$\Hh := \Hh^+ \oplus \Hh^-$, where the representations on $\Hh^\pm$ 
are given in Subsection \ref{sec-4}
with $\Hh^\pm_0=\dC$. For $r=0$, we take the 
representation on $\Hh := \Hh^+$ with $\Hh^+_0=\dC$. 
Since the operator $A$ has an orthonormal basis of eigenvectors with 
eigenvalues $q^{2n} \lambda_\pm$ and 
each eigenvalue has multiplicity one, $A$ is of trace 
class and so is $|A|x$ for all $x\in \cF_{\mathrm{b}}(\mathrm{S}^2_{qr})$. 
Obviously, 
$\tr_\Hh \,|A| =(1\!-\!q^2)^{-1}(\lambda_+\!-\!\lambda_-)$. Therefore,
\begin{equation}\label{htrace}
h(x) = (1\!-\!q^2)(\lambda_+\!-\!\lambda_-)^{-1}\,\tr_\Hh \,|A|x,\quad 
x\in\cF_{\mathrm{b}}(\mathrm{S}^2_{qr}),
\end{equation}
defines a state on the $\ast$-al\-ge\-bra 
$\cF_{\mathrm{b}}(\mathrm{S}^2_{qr})$.

Next we show that $h$ is $\cU_q(\mathrm{su}_2)$-in\-var\-iant 
on $\cF^\infty(\mathrm{S}^2_{qr})$. 
We carry out this verification in the case $r\in(0,\infty]$ for 
the generator $E$ and for an element 
$x=B^{\# n} f(A)$ of $\cF^\infty(\mathrm{S}^2_{qr})$, 
where $n\in\dZ$ and $f\in\cF^\infty(\sigma(A))$. The other cases are 
treated in a similar manner. Since $f\in\cF^\infty(\sigma(A))$, there is 
a bounded function $g$ on $\sigma(A)$ such that $f(A)\!-\!f(0)=Ag(A)$. 
Write $x=x_1+x_2$, where $x_1=B^{\# n} g(A)A$ and $x_2=f(0) B^{\# n}$. 
By (\ref{eactt1}) and (\ref{htrace}), we have
\begin{equation}\label{hey}
h(E\ang y)= {\rm const}\, \tr_\Hh \, |A|A^{-1} [B^\ast, y] = 
{\rm const}\,  (\tr_{\Hh^+}\, [B^\ast,y]-\tr_{\Hh^-}\, [B^\ast,y])
\end{equation}
for $y\in\cF^\infty(\mathrm{S}^2_{qr})$. Since $g(A)$ and $B^{\# n}$ are 
bounded operators and $A$ is of trace class, $x_1$ is of trace class. 
Therefore, $\tr_{\Hh^\pm}\, [B^\ast, x_1]=0$, so that $h(E\ang x_1)=0$ 
by (\ref{hey}). Since $[B^\ast,B^{\# n}]\in A{\cdot} \cO(\mathrm{S}^2_{qr})$ 
by Lemma \ref{1/1} below, $[B^\ast,B^{\# n}]$ is of trace class. 
If $n\ne 1$, then $\tr_{\Hh}\, |A|A^{-1} [B^\ast,B^{\# n}]=0$ because 
in this case $|A|A^{-1}[B^\ast,B^{\# n}]$ acts on $\Hh$ 
as a weighted shift operator. 
Hence 
$h(E\ang x_2)=0$. Suppose that $n=1$. 
Since $[B^\ast,B]=(1-q^2)(A-(1+q^2)A^2)$ and 
$$
(1\!-\!q^2) \tr_{\Hh^\pm} (A\!-\!(1+q^2)A^2)
=\lambda_\pm \!-\! \lambda^2_\pm = \lambda_+\lambda_-,
$$
it follows from (\ref{hey}) that $h(E\ang x_2)=0$. Thus, $h(E\ang x)=0$. 
This completes the proof of the $\cU_q(\mathrm{su}_2)$-invariance of $h$ on 
$\cF^\infty(\mathrm{S}^2_{qr})$.


In both approaches, we have proved the following theorem.
\begin{tht}                                                   \label{P1}
With the foregoing  definitions, $\cF^\infty(\mathrm{S}^2_{qr})$ is a left 
$\cU_q(\mathrm{su}_2)$-mo\-dule 
$\ast$-al\-ge\-bra which contains $\cO(\mathrm{S}^2_{qr})$ as a 
$\cU_q(\mathrm{su}_2)$-mo\-dule 
$\ast$-sub\-al\-ge\-bra. The functional $h$ is a faithful 
$\cU_q(\mathrm{su}_2)$-in\-var\-iant 
state on the $\cU_q(\mathrm{su}_2)$-mo\-dule $\ast$-al\-ge\-bra 
$\cF^\infty(\mathrm{S}^2_{qr})$.
\end{tht}

In the next subsection, we shall need the following lemma. 
It is the algebraic version of the operator-theoretic 
formulas (\ref{eactt1})--(\ref{eactt3}) stated above.
\begin{thl}\label{1/1} 
For any $x\in\cO(\mathrm{S}^2_{qr})$, we have
\begin{eqnarray*}
   &   Ax=(K^2\ang x)A,\qquad xA=A(K^{-2}\ang x),  &\\
&   [B^\ast,x] =-q^{1/2}\lambda A (K^{-1} E\ang x),\quad 
  [B,x] =-q^{3/2}\lambda A (K^{-1} F \ang x). &
\end{eqnarray*}
In particular, the commutators $[B,x]$ and $[B^\ast,x]$
are in $A{\cdot} \cO(\mathrm{S}^2_{qr})$.
\end{thl}
It suffices to prove the lemma for elements $x$ of the 
vector space basis 
$\{A^n B^{\# k}\,;$ $n\in \dN_0,$ $k\in \dZ\}$ of $\cO(\mathrm{S}^2_{qr})$. 
This can be done by a straightforward induction ar\-gu\-ment. 
We omit the details.
%
%
\subsection{Description of quantum line bundles by charts}
                                                      \label{sec-repcpa}
Throughout this subsection, we suppose that $j\in\frac{1}{2}\dZ$. 

From Proposition \ref{cor} and Theorems \ref{T1} and \ref{rep}, 
it follows that each irreducible integrable 
$\ast$-re\-pre\-sen\-ta\-tion of 
$\cO(\mathrm{S}^2_{qr})\rti \cU_q(\mathrm{su}_2)$ 
can be realized as a projective module 
$M_j\cong \cO(\mathrm{S}^2_{qr})^{2|j|+1}P_j$. 
It is known that the projective modules $M_j$ can be considered 
as line bundles over the quantum spheres $\mathrm{S}^2_{qr}$ 
\cite{BM,HM}. 
In this section, we describe the quantum line bundles $M_j$ 
by two ``charts''. The charts will be realized by algebras of functions 
on the positive and the negative part of the spectrum of the
self-adjoint operator $A$. 
These function algebras 
lead to $\ast$-re\-pre\-sen\-ta\-tions of the $\ast$-al\-ge\-bra
$\cF_{\mathrm{b}}(\mathrm{S}^2_{qr}) \otimes \cY_r$ 
by left and right multiplications. 
It will be shown that each chart is related to a tensor product of an 
irreducible $\ast$-re\-pre\-sen\-ta\-tion of $\cO(\mathrm{S}^2_{qr})$ 
from Subsection \ref{sec-4} and an irreducible 
$\ast$-re\-pre\-sen\-ta\-tion of the $\ast$-al\-ge\-bra $\cY_r$ from 
Subsection \ref{S4}.

Recall that the isomorphism 
$\Psi_j$ realizing $M_j\cong \cO(\mathrm{S}^2_{qr})^{2|j|+1}P_j$ 
is given by Equation (\ref{psi}).
In what follows, we consider $\cO(\mathrm{S}^2_{qr})^{2|j|+1}P_j$ 
as a subspace of  $\cF_{\mathrm{b}}(\mathrm{S}^2_{qr})^{2|j|+1}$. 
Our next aim is to equip $\cF_{\mathrm{b}}(\mathrm{S}^2_{qr})^{2|j|+1}$ 
with an inner product such that $\Psi_j$ becomes 
an isometry. We begin by proving an auxiliary lemma. 
\begin{thl}                                                  \label{2/2}
Let $x\in\cO(\mathrm{S}^2_{qr})$ and 
$k,l\in\{-|j|,-|j|+1,\ldots,-|j|\}$. Then 
\begin{equation}                                \label{hvxv}
h(v^{|j|\ast}_{kj}xv^{|j|}_{lj})=c_jq^{2|j|-2k}\hspace{1pt}
h(xv^{|j|}_{lj}v^{|j|\ast}_{kj}),
\end{equation}
where $c_j=h(v^{|j|}_{|j|,j}v^{|j|\ast}_{|j|,j})^{-1}$. 
For $y\in \cF_{\mathrm{b}}(\mathrm{S}^2_{qr})$ and 
$g\in \cF_{\mathrm{b}}(\sigma(A))$, we have 
\begin{equation}                                             \label{hxb}
h(y\hs g(A))=h(g(A)\hs y),\quad h(yB)=q^2\hspace{1pt}h(By),\quad 
h(yB^\ast) = q^{-2}\hspace{1pt} h(B^\ast y).
\end{equation}
\end{thl}
{\bf Proof.}
We first show that it suffices to prove (\ref{hvxv}) for $k=|j|$. 
Indeed, from (\ref{vlkj}), it follows that 
there is a non-zero real constant $\gamma^j_k$ such that 
$v^{|j|}_{kj}=\gamma^j_k\hspace{1pt} F^{|j|-k}\ang v^{|j|}_{|j|,j}$. 
Recall that $h((f\ang y)^\ast z)=h(y^\ast f\ang z)$ and
$(f\ang y)^\ast =S(f)^\ast \ang y^\ast$ for $f\in \cU_q(\mathrm{su}_2)$, 
$y,z\in\cO(\mathrm{SU}_q(2))$. 
Assuming that (\ref{hvxv}) holds for $k=|j|$, we get 
\begin{align*}
h(v^{|j|\ast}_{kj}xv^{|j|}_{lj})
&=\gamma^j_k\hspace{1pt} 
h((F^{|j|-k}\ang v^{|j|}_{|j|,j})^\ast xv^{|j|}_{lj})
=\gamma^j_k\hspace{1pt} 
h(v^{|j|\ast}_{|j|,j} E^{|j|-k}\ang (xv^{|j|}_{lj}))\\
&=c_j \gamma^j_k\hspace{1pt} 
h( (S(E^{|j|-k})^\ast \ang(xv^{|j|}_{lj})^\ast)^\ast v^{|j|\ast}_{|j|,j})
=c_j \gamma^j_k\hspace{1pt} 
h( xv^{|j|}_{lj} S(E^{|j|-k}) \ang v^{|j|\ast}_{|j|,j})\\
&=c_j \gamma^j_k\hspace{1pt} 
h( xv^{|j|}_{lj} (S^2(E^{|j|-k})^\ast \ang v^{|j|}_{|j|,j})^\ast)
=c_jq^{2(|j|-k)} \hspace{1pt} 
h( xv^{|j|}_{lj}v^{|j|\ast}_{kj}).
\end{align*}
Next we note that it suffices to prove that 
\begin{equation}                            \label{hvv}
h(v^{|j|\ast}_{|j|,j}v^{|j|+n}_{lj})=c_j\hspace{1pt}
h(v^{|j|+n}_{lj}v^{|j|\ast}_{|j|,j})
\end{equation}
for all $n\in\dN_0$ since $xv^{|j|}_{lj}\in M_j$ and 
the elements $v^{|j|+n}_{lj}$ span $M_j$ by Theorem \ref{T1}. 
As the elements $v^{|j|+n}_{lj}$ form an orthonormal 
set in $\cO(\mathrm{SU}_q(2))$ 
with  inner product defined by (\ref{ip}), we only have to show that 
the right-hand side of (\ref{hvv}) vanishes whenever 
$l\neq |j|$ or $n>0$. (Observe that Equation (\ref{hvv}) 
is trivially satisfied for $l=|j|$ and $n=0$.)  
If $l\neq |j|$, then it follows from 
$K\ang v^{|j|+n}_{lj}v^{|j|\ast}_{|j|,j}
= q^{|j|-l}v^{|j|+n}_{lj}v^{|j|\ast}_{|j|,j}$ and the 
$\cU_q(\mathrm{su}_2)$-invariance of $h$ that 
$h(v^{|j|+n}_{lj}v^{|j|\ast}_{|j|,j})=0$. If $n>0$, 
then $v^{|j|+n}_{lj}= \kappa^{|j|+n}_l\hspace{1pt}F\ang v^{|j|+n}_{l+1,j}$ 
with a non-zero real constant $\kappa^{|j|+n}_l$.  
Hence 
\begin{align*}
h(v^{|j|+n}_{lj}v^{|j|\ast}_{|j|,j}) 
&=\kappa^{|j|+n}_l\hspace{1pt}
h( (S(F)^\ast\ang v^{|j|+n\ast}_{l+1,j})^\ast v^{|j|\ast}_{|j|,j}) 
= \kappa^{|j|+n}_l\hspace{1pt}
h( v^{|j|+n\ast}_{l+1,j} S(F)\ang v^{|j|\ast}_{|j|,j}) \\
&=\kappa^{|j|+n}_l q^{-2}\hspace{1pt}
h( v^{|j|+n\ast}_{l+1,j} (E\ang v^{|j|}_{|j|,j})^\ast) = 0
\end{align*}
since $v^{|j|}_{|j|,j}$ is a highest weight vector. 

To prove (\ref{hxb}), we can assume that 
$y=f(A)B^{\# k}$, where $f\in \cF_{\mathrm{b}}(\sigma(A))$ 
and $k\in\dZ$, because these elements span 
$\cF_{\mathrm{b}}(\mathrm{S}^2_{qr})$. 
Since $h(y\hs g(A))=h(g(A)\hs y)=0$ for $k\ne 0$, 
the first equality of (\ref{hxb}) 
is obvious. The third equality of (\ref{hxb}) follows from the second one  
because the state $h$ is hermitian. As $h(yB)=h(By)=0$ for $k\ne -1$, 
it remains to treat the case $k=-1$. 
Using the relations $B^\ast B=(A\!-\!\lambda_-)(\lambda_+ \!-\! A)$, 
$BB^\ast =(q^2 A\!-\!\lambda_-) (\lambda_+ \!-\! q^2 A)$ 
and Equation (\ref{hqf1}), we get\\
\begin{parbox}[c]{0.4cm}{ 
\begin{align*}
&\\
&\\
& 
\end{align*} } \end{parbox} \hfill
\begin{parbox}[c]{12.0cm}{
\begin{align*}
h(yB)&=h(f(A)B^\ast B)=h(f(A)(A\!-\!\lambda_-)(\lambda_+ \!-\! A))\\
&= h(q^2 f(q^2 A)(q^2 A\!-\!\lambda_-)(\lambda_+ \!-\! q^2 A))
=h(q^2 f(q^2 A)BB^\ast)\\
&=h(q^2 Bf(A)B^\ast)=h(q^2 By).  
\end{align*}  } \end{parbox} \hfill
\begin{parbox}[c]{0.4cm}{ 
\begin{align*}
&\\
&\\
\Box
\end{align*} } \end{parbox} 
\begin{thp}                                                \label{6.4}
Let $\cL_2(\hspace{-0.3pt}\mathrm{S}^2_{qr}\hspace{-0.3pt})^{2|j|+1}$ 
be  the Hilbert space completion of 
$\cF_{\mathrm{b}}(\hspace{-0.3pt}\mathrm{S}^2_{qr}\hspace{-0.3pt})^{2|j|+1}$ 
with respect 
to the inner product given by 
\begin{equation}                                       \label{ipj}
 \langle (y_{-|j|},\ldots,y_{|j|}),(z_{-|j|},\ldots,z_{|j|})\rangle
=c_j q^{2|j|}\sumop_{k=-|j|}^{|j|}q^{-2k}h(z_k^\ast y_k),
\end{equation}
where $c_j=h(v^{|j|}_{|j|,j}v^{|j|\ast}_{|j|,j})^{-1}$. 
Then the right multiplication by $P_j$ on 
$\cF_{\mathrm{b}}(\mathrm{S}^2_{qr})^{2|j|+1}$
defines an orthogonal projection on the Hilbert space 
$\cL_2(\mathrm{S}^2_{qr})^{2|j|+1}$. 
The isomorphism $\Psi_j:\cO(\mathrm{S}^2_{qr})^{2|j|+1}P_j 
\stackrel{\cong}{\longrightarrow}M_j$ 
from Equation (\ref{psi}) is an isometry.
\end{thp}
{\bf Proof.} 
As $P_j^2=P_j$, we only have to show that $P_j$ is self-adjoint with respect 
to the inner product (\ref{ipj}). Since 
$K\ang v^{|j|}_{lj}v^{|j|\ast}_{kj}\!=\!q^{l-k}v^{|j|}_{lj}v^{|j|\ast}_{kj}$, 
there exists a polynomial $p^{j}_{lk}(A)$ such that 
$v^{|j|}_{lj}v^{|j|\ast}_{kj}\!=\!B^{\#k-l }p^{j}_{lk}(A)$. 
Hence, by (\ref{hxb}),
\begin{equation}                                              \label{hB}
h(xv^{|j|}_{lj}v^{|j|\ast}_{kj})=h(q^{2(k-l)}v^{|j|}_{lj}v^{|j|\ast}_{kj}x)
\end{equation}
for $x\in \cF_{\mathrm{b}}(\mathrm{S}^2_{qr})$. 
This gives 
\begin{align*}
&\langle  (y_{-|j|},\ldots,y_{|j|})P_j, (z_{-|j|},\ldots,z_{|j|})\rangle \\
&\qquad\qquad\qquad \qquad
=c_j q^{2|j|} \!\!\sumop_{k,l=-|j|}^{|j|}\!\! q^{-2k}
h( [2|j|\!+\!1]^{-1}q^{-(l+k)}z_k^\ast~ y_l v^{|j|}_{lj} v^{|j|\ast}_{kj})\\
&\qquad\qquad\qquad \qquad 
=c_j q^{2|j|} \!\!\sumop_{k,l=-|j|}^{|j|}\!\!\!q^{-2l}
h( [2|j|\!+\!1]^{-1}q^{-(l+k)}v^{|j|}_{lj} v^{|j|\ast}_{kj}z_k^\ast y_l )\\
&\qquad\qquad\qquad \qquad   
=\langle  (y_{-|j|},\ldots,y_{|j|}), (z_{-|j|},\ldots,z_{|j|})P_j\rangle,
\end{align*}
which proves the first assertion of the proposition. 

Let $(y_{-|j|},\hs.\hs.\hs.\hs,y_{|j|})
\hsp=\hsp(y_{-|j|},\hs.\hs.\hs.\hs,y_{|j|})P_j\in 
\cO(\mathrm{S}^2_{qr})^{2|j|+1}P_j$ and 
$(z_{-|j|},\hs.\hs.\hs.\hs,z_{|j|})=(z_{-|j|},\hs.\hs.\hs.\hs,z_{|j|})P_j\in 
\cO(\mathrm{S}^2_{qr})^{2|j|+1}P_j$. From the definition of $P_j$, 
it follows that $y_n=[2|j|\!+\!1]^{-1} 
\sum_{k=-|j|}^{|j|} q^{-(k+n)} y_k v^{|j|}_{kj} v^{|j|\ast}_{nj}$. Using 
Lemma \ref{2/2} and the last identity, we compute
\begin{align*}
&\langle \Psi_j(y_{-|j|},\ldots,y_{|j|}),
\Psi_j(z_{-|j|},\ldots,z_{|j|})\rangle\\
&\qquad\qquad= [2|j|\!+\!1]^{-1}h\big(
(\!\sumop^{|j|}_{n=-|j|} q^{-n} z_n v^{|j|}_{nj})^\ast
\!\sumop^{|j|}_{k=-|j|} q^{-k} y_k v^{|j|}_{kj} \big)\\
&\qquad\qquad=[2|j|\!+\!1]^{-1}\!\sumop^{|j|}_{n,k=-|j|}h\big( q^{-(n+k)}
v^{|j|\ast}_{nj}z_n^\ast y_k v^{|j|}_{kj} \big)\\
&\qquad\qquad=[2|j|\!+\!1]^{-1}\!\sumop^{|j|}_{n,k=-|j|}
c_j q^{2|j|-2n}h\big( q^{-(n+k)}
z_n^\ast y_k v^{|j|}_{kj}v^{|j|\ast}_{nj}\big)\\
&\qquad\qquad=c_j q^{2|j|}\!\sumop^{|j|}_{n=-|j|}q^{-2n}h(z_n^\ast y_n)
=\langle (y_{-|j|},\ldots,y_{|j|}),(z_{-|j|},\ldots,z_{|j|})\rangle,
\end{align*}
which shows the second assertion of the proposition.  \hfill$\Box$
\mn 

Since $\Psi_j$ is an isometric isomorphism, the $\ast$-re\-pre\-sen\-ta\-tion 
$\hat\pi_j$ of the cross product algebra 
$\cO(\mathrm{S}^2_{qr})\rti \cU_q(\mathrm{su}_2)$ on $M_j$ 
is unitarily equivalent to the $\ast$-re\-pre\-sen\-ta\-tion 
${\check{\pi}}_j:=\Psi_j^{-1}\circ\hat\pi_j\circ\Psi_j$ 
on $\cO(\mathrm{S}^2_{qr})^{2|j|+1}P_j$. 
The restriction of ${\check{\pi}}_j$ to 
$\cO(\mathrm{S}^2_{qr}) \otimes \cY_r$
can be described by left and right multiplications by matrices 
with entries from $\cO(\mathrm{S}^2_{qr})$. Clearly, 
for $z\in\cO(\mathrm{S}^2_{qr})$, 
$ {\check{\pi}}_j(z) (y_{-|j|},\ldots,y_{|j|})P_j
=z(y_{-|j|},\ldots,y_{|j|})P_j$.
Since $X=q^{3/2}\lambda A FK^{-1}  + q^{-1}B$ commutes with 
$y_k\in\cO(\mathrm{S}^2_{qr})$, we obtain
\begin{align*}
&[2|j|\!+\!1]^{1/2}{\check{\pi}}_j(X)(y_{-|j|},\ldots,y_{|j|})=\Psi_j^{-1}(
 \sumop^{|j|}_{k=-|j|} q^{-k} y_k X v^{|j|}_{kj})\\
&\qquad=\Psi_j^{-1}(
\sumop^{|j|}_{k=-|j|} q^{-k} y_k( 
q^{3/2-k}\lambda [|j|\!-\!k\!+\!1]^{1/2}[|j|\!+\!k]^{1/2} 
Av^{|j|}_{k-1,j}+q^{-1}Bv^{|j|}_{kj})). 
\end{align*}
This shows that ${\check{\pi}}_j(X)((y_{-|j|},\ldots,y_{|j|})P_j)
=((y_{-|j|},\ldots,y_{|j|})\fm_j)P_j$, where the
matrix 
$\fm_j=\big(m^j_{kl}\big)_{k,l=-|j|}^{|j|}
\in{\mathrm M}_{2|j|+1} (\cO(\mathrm{S}^2_{qr}) ) $ has the entries 
\begin{equation*}
m^j_{kk}=q^{-1}B,\ \,
m^j_{k,k-1}=q^{1/2-k}\lambda [|j|\!-\!k\!+\!1]^{1/2}[|j|\!+\!k]^{1/2}A,\ \,
m^j_{kl}=0,\ l\neq k,k-1.  
\end{equation*}
Similarly, we have ${\check{\pi}}_j(X^\ast)((y_{-|j|},\ldots,y_{|j|})P_j)
=((y_{-|j|},\ldots,y_{|j|})\fm^\dagger_j)P_j$ with matrix 
$\fm^\dagger_j=\big(m^{\dagger j}_{kl}\big)_{k,l=-|j|}^{|j|}$ given by 
\begin{equation*}
m^{\dagger j}_{kk}=qB^\ast,\ \,
m^{\dagger j}_{k-1,k}
=q^{5/2-k}\lambda [|j|\!-\!k\!+\!1]^{1/2}[|j|\!+\!k]^{1/2}A,\ \,
m^{\dagger j}_{kl}=0,\ l\neq k,k-1
\end{equation*}
and ${\check{\pi}}_j(Y)((y_{-|j|},\ldots,y_{|j|})P_j)
=((y_{-|j|},\ldots,y_{|j|})\fn_j)P_j$ with 
$\fn_j=\big(n^{j}_{kl}\big)_{k,l=-|j|}^{|j|}$ given by 
\begin{equation*}
n^j_{kk}=q^{-2k+1}A,\quad 
n^j_{kl}=0,\ \, l\neq k.  
\end{equation*}
It is easy to check that 
\begin{equation}                                \label{MNrel}
\fn_j\fm^\dagger_j = q^2 \fm^\dagger_j\fn_j,\ \, 
\fn_j\fm_j = q^{-2} \fm_j \fn_j,\ \,
\fm_j \fm^\dagger_j - q^2 \fm^\dagger_j\fm_j = (1-q^2) (\fn_j^2+r).
\end{equation}
Note that all operators are bounded. 
Hence the restriction of ${\check{\pi}}_j$ to 
$\cO(\mathrm{S}^2_{qr})\otimes \cY_r$ yields 
a bounded $\ast$-re\-pre\-sen\-ta\-tion 
${\check{\pi}}^\mathrm{b}_j$ 
on the Hilbert space 
$\cL_2(\mathrm{S}^2_{qr})^{2|j|+1}P_j$.
This representation can be extended to a bounded 
$\ast$-re\-pre\-sen\-ta\-tion,  
denoted again by ${\check{\pi}}^\mathrm{b}_j$, 
of $\cF_{\mathrm{b}}(\mathrm{S}^2_{qr})\otimes \cY_r$
on $\cL_2(\mathrm{S}^2_{qr})^{2|j|+1}P_j$ 
such that ${\check{\pi}}^\mathrm{b}_j(f)$, 
$f\in\cF_{\mathrm{b}}(\mathrm{S}^2_{qr})$, acts on 
$\cF_{\mathrm{b}}(\mathrm{S}^2_{qr})^{2|j|+1}P_j$ by left multiplication. 

Let $\chi_-$ and $\chi_+$ denote the characteristic functions of 
the intervals 
$(-\infty,0]$ and $[0,\infty)$, respectively. Set 
$\cF_{\mathrm{b}}(\mathrm{S}^2_{qr})_+
:=\cF_{\mathrm{b}}(\mathrm{S}^2_{qr})\chi_+(A)$ and, 
for $r>0$,  
$\cF_{\mathrm{b}}(\mathrm{S}^2_{qr})_-
:=\cF_{\mathrm{b}}(\mathrm{S}^2_{qr})\chi_-(A)$. Since 
$\chi_\pm(A)$ commutes with all elements of 
$\cF_{\mathrm{b}}(\mathrm{S}^2_{qr})$, 
$\cF_{\mathrm{b}}(\mathrm{S}^2_{qr})_\pm$ is a unital 
$\ast$-al\-ge\-bra with unit element $\chi_\pm(A)$ and an ideal 
of $\cF_{\mathrm{b}}(\mathrm{S}^2_{qr})$.  
Our next aim is to describe the ``charts''  
$\cF_{\mathrm{b}}(\mathrm{S}^2_{qr})^{2|j|+1}_\pm P_j$
of the quantum line bundle $\cF_{\mathrm{b}}(\mathrm{S}^2_{qr})^{2|j|+1}P_j$ 
by the function 
algebras $\cF_{\mathrm{b}}(\mathrm{S}^2_{qr})_-$ and 
$\cF_{\mathrm{b}}(\mathrm{S}^2_{qr})_+$ themselves. 
%
%
\begin{thl}                                       \label{vjjj}
Set $\xi(s):= -sab+(s^2-1)qbc+sq^2dc$, where the parameter $s$ is defined  
as in Section \ref{sec-decomp}. Then  
\begin{enumerate}
\item[(i)]
$(a-qsc)(d+sb)=1-\xi(s),\quad (d+q^{-1}sb)(a-sc)=1-q^{-2}\xi(s),$
\smallskip\\
$(b-qsd)(-qc-qsa)=q^2s^2+\xi(s),\quad (qc+sa)(-b+sd)=s^2+\xi(s).$

\item[(ii)]
$(a-qsc)\xi(qs)\!=\!q^2\xi(s)(a-qsc),$\ \ \smallskip 
$(d+q^{-1}sb)\xi(q^{-1}s)\!=\!q^{-2}\xi(s)(d+q^{-1}sb),$\\
$(b-qsd)\xi(qs)\!=\!\xi(s)(b-qsd),$\ \ 
$(c+q^{-1}sa)\xi(q^{-1}s)\!=\!\xi(s)(c+q^{-1}sa).$

\item[(iii)]
$v^j_{jj}v^{j\ast}_{jj}
=\gamma^+_j(\lambda_+-A)(\lambda_+-q^{-2}A)\dots(\lambda_+-q^{-4j+2}A),
\quad j>0,$\smallskip\\
$v^j_{-j,j}v^{j\ast}_{-j,j}
=\gamma^-_j(q^2A-\lambda_-)(q^4A-\lambda_-)\dots(q^{4j}A-\lambda_-),
\quad j>0,$\smallskip\\
$v^{|j|}_{jj}v^{|j|\ast}_{jj}
=\gamma^+_j(\lambda_+-q^2A)(\lambda_+-q^{4}A)\dots(\lambda_+-q^{4|j|}A),
\quad j<0,$\smallskip\\
$v^{|j|}_{|j|,j}v^{|j|\ast}_{|j|,j}
=\gamma^-_j(A-\lambda_-)(q^{-2}A-\lambda_-)\dots(q^{-4|j|+2}A-\lambda_-),
\quad j<0$,\smallskip\\
with non-zero constants $\gamma^\pm_j\in\dR$.

\item[(iv)]
The function $v^{|j|}_{jj}v^{|j|\ast}_{jj}\chi_-(A)$ is invertible 
in $\cF_{\mathrm{b}}(\mathrm{S}^2_{qr})_-$ for $r>0$.\smallskip\\
The function $v^{|j|}_{-j,j}v^{|j|\ast}_{-j,j}\chi_+(A)$ is invertible 
in $\cF_{\mathrm{b}}(\mathrm{S}^2_{qr})_+$.

\item[(v)]
$v^{|j|}_{kj}v^{|j|\ast}_{jj}
\big(v^{|j|}_{jj}v^{|j|\ast}_{jj}\chi_-(A)\big)^{-1}
v^{|j|}_{jj}v^{|j|\ast}_{lj}=v^{|j|}_{kj}v^{|j|\ast}_{lj}\chi_-(A)$\ \,for 
$r>0$, 
\medskip\\
$v^{|j|}_{kj}v^{|j|\ast}_{-j,j}
\big(v^{|j|}_{-j,j}v^{|j|\ast}_{-j,j}\chi_+(A)\big)^{-1}
v^{|j|}_{-j,j}v^{|j|\ast}_{lj}=v^{|j|}_{kj}v^{|j|\ast}_{lj}\chi_+(A).$
\end{enumerate}
\end{thl}
{\bf Proof.}
(i) and (ii) follow by straightforward computations. 

(iii) Let $j>0$. Then $v^j_{jj}=N^j_{jj}u_j$, where $u_j$ is defined 
by Equation (\ref{uj}). From (i) and (ii), it follows that 
$$
(d + q^{-1}s b)...(d+q^{-2j}s b)(a-q^{-2j+1} s c)... (a -  sc)
=(1-q^{-2}\xi(s))... (1-q^{-4j}\xi(s)). 
$$
Inserting $\xi(s)=sr^{-1/2}q^2A=\lambda_+^{-1}q^2A$ gives the result. 
The other cases are treated analogously.

(iv) As $A<0$ on the interval $(-\infty,0]$ and $\lambda_+>0$, 
each factor $(\lambda_+-q^{2k}A)$, $k\in\dZ$, is invertible
in $\cF_{\mathrm{b}}(\mathrm{S}^2_{qr})_-$. 
Likewise, for $r>0$, each factor  $(q^{2k}A-\lambda_-)$, $k\in\dZ$, is 
invertible in  $\cF_{\mathrm{b}}(\mathrm{S}^2_{qr})_+$ since $\lambda_-<0$. 

(v) Let $j>0$. In the proof of Proposition \ref{6.4}, we argued that 
$v^{j}_{kj}v^{j\ast}_{jj}\!=\!B^{j-k }p^{j}_{kj}(A)$. 
Note that $B\xi(s)=q^2\xi(s)B$. 
This together with (iii) gives 
\begin{align*}
&v^{j}_{kj}v^{j\ast}_{jj}
(v^{j}_{jj}v^{j\ast}_{jj}\chi_-(A))^{-1}
v^{j}_{jj}v^{j\ast}_{lj}\\
&\qquad\quad  =\big(N^{j2}_{jj}
(1\!-\!q^{2(j-k-1)}\xi(s))...(1\!-\!q^{-2(j+k)}\xi(s))\chi_-(A)\big)^{-1}
v^{j}_{kj}v^{j\ast}_{jj}v^{j}_{jj}v^{j\ast}_{lj}. 
\end{align*}
On the other hand, 
$v^{j\ast}_{jj}v^{j}_{jj}=N^{j2}_{jj}
(1-\xi(q^{-2j}s))...(1-q^{4j-2}\xi(q^{-2j}s))$ 
by (i) and (ii). 
The vector $v^{j}_{kj}$ is a linear combination of products consisting 
of factors $(d+q^{-n}sb)$ and $(c+q^{-m}sa)$, $1\leq n,m\leq 2j$, 
where the terms  $(d+q^{-n}sb)$ and $(c+q^{-m}sa)$ appear 
$j+k$-times and $j-k$-times, respectively. 
Applying (i) and (ii), we see that 
$v^{j}_{kj}v^{j\ast}_{jj}v^{j}_{jj}=N^{j2}_{jj}
(1-q^{-2(j+k)}\xi(s))...(1-q^{2(j-k-1)}\xi(s))v^{j}_{kj}$. 
Inserting this identity into above equation gives the result. 
The other cases are handled similarly.          \hfill $\Box$
\begin{thl}                               \label{charts}
With the inner product on $\cF_{\mathrm{b}}(\mathrm{S}^2_{qr})_\pm$ 
defined by $\langle f,g\rangle:=h(g^\ast f)$, 
there is an isometric isomorphism $\check{\Psi}_{j,\pm}$ from 
$\cF_{\mathrm{b}}(\mathrm{S}^2_{qr})_\pm$ onto 
$\cF_{\mathrm{b}}(\mathrm{S}^2_{qr})^{2|j|+1}_\pm P_j$
given by 
\begin{align*}
\check{\Psi}_{j,-}(f)
&=q^{2j}[2j\!+\!1]^{1/2}
(0,\ldots,0,f(v^{j}_{jj}v^{j\ast}_{jj}\chi_-(A))^{-1/2})P_j,\quad 
   j>0,\\
\check{\Psi}_{j,-}(f)
&=q^{2j}[2|j|\!+\!1]^{1/2}
(f(v^{|j|}_{jj}v^{|j|\ast}_{jj}\chi_-(A))^{-1/2},0,\ldots,0)P_j,\quad 
   j<0,\\
\check{\Psi}_{j,+}(g)
&=q^{-2j}[2j\!+\!1]^{1/2}
(g(v^{j}_{-j,j}v^{j\ast}_{-j,j}\chi_+(A))^{-1/2},0,\ldots,0)P_j,\quad
  j>0,\\
\check{\Psi}_{j,+}(g)
&=q^{-2j}[2|j|\!+\!1]^{1/2}
(0,\ldots,0,g(v^{|j|}_{|j|,j}v^{|j|\ast}_{|j|,j}\chi_+(A))^{-1/2})P_j,\quad 
  j<0,
\end{align*}
where $f\in    \cF_{\mathrm{b}}(\mathrm{S}^2_{qr})_-$ and  
$g\in    \cF_{\mathrm{b}}(\mathrm{S}^2_{qr})_+$.  
\end{thl}
{\bf Proof.}
We carry out the proof for $\check{\Psi}_{j,-}$ and  $j>0$. 
The other cases are 
treated in the same manner. 
By Lemma \ref{vjjj}, $v^{j}_{jj}v^{j\ast}_{jj}\chi_-(A)$ is invertible. 
Hence $\check{\Psi}_{j,-}$ is well defined. 
Fix
$(y_{-j},\ldots,y_{j})=(y_{-j},\ldots,y_{j})P_j\in
\cF_{\mathrm{b}}(\mathrm{S}^2_{qr})^{2j+1}_- P_j $. 
Let $z\!=\!\sum_{k=-j}^j q^{j-k}y_k
v^{j}_{kj}v^{j\ast}_{jj}
\big(v^{j}_{jj}v^{j\ast}_{jj}\chi_-(A)\big)^{-1}$ and 
$(z_{-j},\ldots,z_j)=(0,\ldots,0,z)P_j$. 
By Lemma \ref{vjjj}, 
\begin{align*}
z_l&=[2j\!+\!1]^{-1}\sumop_{k=-j}^j q^{-(k+l)} y_k 
v^{j}_{kj}v^{j\ast}_{jj}
\big(v^{j}_{jj}v^{j\ast}_{jj}\chi_-(A)\big)^{-1}v^{j}_{jj}v^{j\ast}_{lj}\\
&=[2j\!+\!1]^{-1}\sumop_{k=-j}^j 
q^{-(k+l)} y_k v^{j}_{kj}v^{j\ast}_{lj}\chi_-(A)=y_l.
\end{align*}
Thus 
$\check{\Psi}_{j,-}(q^{-2j}[2j\!+\!1]^{-1/2}z(v^{j}_{jj}
v^{j\ast}_{jj}\chi_-(A))^{1/2})
=(y_{-j},\ldots,y_{j})$, so $\check{\Psi}_{j,-}$ is surjective. 

Next we verify that $\check{\Psi}_{j,-}$ is isometric. 
Let $f\in\cF_{\mathrm{b}}(\mathrm{S}^2_{qr})_-$. 
Since $P_j$ is a projection, we have 
$[2j\!+\!1]^{-2}
\sum_{k=-j}^j q^{-2(j+k)}v^{j}_{jj}v^{j\ast}_{kj}v^{j}_{kj}v^{j\ast}_{jj}
=[2j\!+\!1]^{-1}q^{-2j}v^{j}_{jj}v^{j\ast}_{jj}$.
Using Equations (\ref{hxb}) and (\ref{hB}), we conclude that 
$$
h(x(v^{j}_{jj}v^{j\ast}_{jj}\chi_-(A))^{-1/2}v^{j}_{jj}v^{j\ast}_{kj})
=q^{2(k-j)}
h((v^{j}_{jj}v^{j\ast}_{jj}\chi_-(A))^{-1/2}v^{j}_{jj}v^{j\ast}_{kj}x).
$$
From these relations, it follows that, 
for $f\in\cF_{\mathrm{b}}(\mathrm{S}^2_{qr})_-$, 
\begin{align*}
\|\check{\Psi}_{j,-}(f)\|^2
&=[2j\!+\!1]^{-1}q^{4j}\sumop_{k=-j}^j q^{-2k}\\
&\quad\  
{\raisebox{0.4ex}{\scriptsize$\times$}}
h(q^{-2(j+k)}v^{j}_{kj}v^{j\ast}_{jj}
(v^{j}_{jj}v^{j\ast}_{jj}\chi_-(A))^{-1/2}f^\ast
f(v^{j}_{jj}v^{j\ast}_{jj}\chi_-(A))^{-1/2}v^{j}_{jj}v^{j\ast}_{kj})\\
&=h(f^\ast f)=\|f\|^2. \hspace{232pt}\Box                          
\end{align*}

As multiplication by elements of 
$\cF_{\mathrm{b}}(\mathrm{S}^2_{qr})$ leaves the 
decomposition $\cF_{\mathrm{b}}(\mathrm{S}^2_{qr})
=\cF_{\mathrm{b}}(\mathrm{S}^2_{qr})_-\oplus
\cF_{\mathrm{b}}(\mathrm{S}^2_{qr})_+$ invariant, we have 
$\cF_{\mathrm{b}}(\mathrm{S}^2_{qr})^{2|j|+1}P_j
=\cF_{\mathrm{b}}(\mathrm{S}^2_{qr})^{2|j|+1}_-P_j\oplus
\cF_{\mathrm{b}}(\mathrm{S}^2_{qr})^{2|j|+1}_+P_j$ 
and the 
$\ast$-re\-pre\-sen\-ta\-tion ${{\check \pi}}^\mathrm{b}_j$ of 
$\cF_{\mathrm{b}}(\mathrm{S}^2_{qr}) \otimes \cY_r$ decomposes 
into a direct sum
${{\check \pi}}^\mathrm{b}_j
={{\check \pi}}^\mathrm{b}_{j,-} \oplus {{\check \pi}}^\mathrm{b}_{j,+}$ 
of $\ast$-re\-pre\-sen\-ta\-tions 
${{\check \pi}}^\mathrm{b}_{j,\pm}$ on 
$\cF_{\mathrm{b}}(\mathrm{S}^2_{qr})^{2|j|+1}_\pm P_j$.  
Using the isometric isomorphism $\check{\Psi}_{j,\pm}$, 
the $\ast$- representation 
${{\check \pi}}^\mathrm{b}_{j,\pm}$ is unitarily equivalent to a 
$\ast$-re\-pre\-sen\-ta\-tion $\rho_{j,\pm}:= 
(\check{\Psi}_{j,\pm})^{-1}\circ 
{{\check \pi}}^\mathrm{b}_{j,\pm}\circ 
\check{\Psi}_{j,\pm}$ of 
$\cF_{\mathrm{b}}(\mathrm{S}^2_{qr}) \otimes \cY_r$ on
$\cF_{\mathrm{b}}(\mathrm{S}^2_{qr})_\pm$.
\begin{tht}                                               \label{reprho} 
The $\ast$-re\-pre\-sen\-ta\-tion
$\rho_{j,\pm}$ of $\cF_{\mathrm{b}}(\mathrm{S}^2_{qr}) \otimes \cY_r$ 
on $\cF_{\mathrm{b}}(\mathrm{S}^2_{qr})_\pm$  is given by 
\begin{align*}
&\rho_{j,-}(X)(f)=q^{-1}fB(\lambda_+ -q^{-4j}A)^{1/2}(\lambda_+ -A)^{-1/2},
\quad \rho_{j,-}(Y)(f)=q^{-2j+1}fA,\\
&\rho_{j,-}(X^\ast)(f)
=qf(\lambda_+-q^{-4j}A)^{1/2}(\lambda_+-A)^{-1/2}B^\ast,\\[-30pt]
\end{align*}
\begin{align*}
&\rho_{j,+}(X)(f)=q^{-1}fB(q^{4j}A-\lambda_-)^{1/2}(A-\lambda_-)^{-1/2},
\quad \rho_{j,+}(Y)(f)=q^{2j+1}fA,\hspace{13pt}\\
&\rho_{j,+}(X^\ast)(f)=qf(q^{4j}A-\lambda_-)^{1/2}(A-\lambda_-)^{-1/2}B^\ast,
\end{align*}
and $\rho_{j,\pm}(x)(f)=xf$ for $x\in\cF_{\mathrm{b}}(\mathrm{S}^2_{qr})$, 
$f\in \cF_{\mathrm{b}}(\mathrm{S}^2_{qr})_\pm$. 

In particular,
all representation operators $\rho_{j,\pm}(y)$, 
$y\in \cF_{\mathrm{b}}(\mathrm{S}^2_{qr}) \otimes \cY_r$, are bounded and 
$\rho_{j,\pm}$ extends to a $\ast$-re\-pre\-sen\-ta\-tion, 
denoted also by $\rho_{j,\pm}$, on the Hilbert space completion 
$\cL_2(\mathrm{S}^2_{qr})_\pm$
of $\cF_{\mathrm{b}}(\mathrm{S}^2_{qr})_\pm$. 
The restriction of this 
$\ast$-re\-pre\-sen\-ta\-tion $\rho_{j,^\pm}$ to 
$\cO(\mathrm{S}^2_{qr})\otimes \cY_r$ is unitarily equivalent 
to a tensor product representation $\sigma^\pm\otimes\sigma^\pm_j$ on 
$\Hh^\pm\otimes\cK^\pm$, where $\sigma^\pm$ denotes the irreducible 
$\ast$-re\-pre\-sen\-ta\-tion of $\cO(\mathrm{S}^2_{qr})$ on $\Hh^\pm$ 
from Subsection \ref{sec-4} and $\sigma^\pm_j$ denotes the irreducible 
$\ast$-re\-pre\-sen\-ta\-tion of $\cY_r$ on $\cK^\pm$ from 
Subsection \ref{S4}
with $Y_0=q^{\pm 2j+1}\lambda_\pm$. 
The restriction of the $\ast$-re\-pre\-sen\-ta\-tion 
${{\check \pi}}_j\cong \hat{\pi}_j\cong \pi_j$
of $\cO(\mathrm{S}^2_{qr})\rti \cU_q(\mathrm{su}_2)$ 
to $\cO(\mathrm{S}^2_{qr})\otimes \cY_r$ 
is unitarily equivalent to the direct sum 
representation $\rho_{j,-}\oplus\rho_{j,+}$. 
\end{tht}
{\bf Proof.}
Clearly, for $x\in\cF_{\mathrm{b}}(\mathrm{S}^2_{qr})$, 
$\rho_{j,- }(x)$ acts by left 
multiplication, that is, 
$\rho_{j,-}(x)(g)=xg$, $g\in\cF_{\mathrm{b}}(\mathrm{S}^2_{qr})_-$. 
Our next aim is to compute the action of the generators of $\cY_r$. 
Let $g\in\cF_{\mathrm{b}}(\mathrm{S}^2_{qr})_-$. 
By the definition the matrix $\fn_j$, 
\begin{align*}
&\rho_{j,-}(Y)(g)=(\check{\Psi}_{j,-})^{-1}(q^{2j}[2j\!+\!1]^{1/2}
(0,...,0,f(v^{j}_{jj}v^{j\ast}_{jj}\chi_-(A))^{-1/2})\fn_jP_j)\\
&=\hsp (\check{\Psi}_{j,-})^{-1}(q^{2j}[2j\!+\!1]^{1/2}
(0,...,0,q^{-2j+1}fA(v^{j}_{jj}v^{j\ast}_{jj}\chi_-(A))^{-1/2})P_j)
\hsp =\hsp q^{-2j+1}fA.
\end{align*}
From Lemma \ref{vjjj} (iii), it follows that 
$$
(v^{j}_{jj}v^{j\ast}_{jj}\chi_-(A))^{-1/2}B^\ast
(v^{j}_{jj}v^{j\ast}_{jj}\chi_-(A))^{1/2}\hspace{-0.2pt}=\hspace{-0.2pt}
\chi_-(A)(\lambda_+-q^{-4j}A)^{1/2}(\lambda_+-A)^{-1/2}B^\ast.
$$
Using this identity and the explicit form of the matrix $\fm^\dagger_j$, 
we see that 
\begin{align*}
\rho_{j,-}(X^\ast)(g)&=(\check{\Psi}_{j,-})^{-1}(q^{2j}[2j\!+\!1]^{1/2}
(0,...,0,f(v^{j}_{jj}v^{j\ast}_{jj}\chi_-(A))^{-1/2})\fm^\dagger_jP_j)\\
&=qg(\lambda
_+-q^{-4j}A)^{1/2}(\lambda_+-A)^{-1/2}B^\ast.
\end{align*}
The operator $\rho_{j,-}(X)$ is determined by the relation 
$\rho_{j,-}(X)=\rho_{j,-}(X^\ast)^\ast$.  Since 
$h(g^\ast f\varphi(A)B^\ast)=h(q^{-2}(gB\bar\varphi(A))^\ast f)$ for 
all $f,g,\varphi(A)\in\cF_{\mathrm{b}}(\mathrm{S}^2_{qr})_-$, we obtain 
$$
\rho_{j,-}(X)(g)
=q^{-1}gB(\lambda_+-q^{-4j}A)^{1/2}(\lambda_+-A)^{-1/2}.
$$
This proves the formulas of the theorem for $\rho_{j,-}$, $j>0$. 
The other cases are treated in the same way. From the preceding formulas 
it is clear that all representation operators
$\rho_{j,\pm}(y)$, 
$y\in \cF_{\mathrm{b}}(\mathrm{S}^2_{qr}) \otimes \cY_r$, are bounded.

For $n\in \dN_0$, let 
$\chi_{n,\pm}$ denote the characteristic 
function of the point $q^{2n} \lambda_\pm$. 
Define
\begin{equation}                                               \label{ob}
\vartheta^\pm_{nl}
:=  c_{nl}^\pm \chi_{n,\pm} (A) B^{\#l },\quad
n\in \dN_0,\ \ l\in \dZ,\  l\ge -n,
\end{equation}
where, with $\gamma_\pm$ defined in Subsection \ref{sec-algfun}, 
\begin{align*}
&  c_{n,0}^\pm=\gamma_\pm^{-1/2}q^{-n}, \\
 &  c_{nl}^\pm= \gamma_\pm^{-1/2} q^{-n-l}\big(
\mathop{\mbox{$\prod$}}_{m=0}^{|l|-1}(q^{2(n-m)}\lambda_\pm-\lambda_-)
(\lambda_+-q^{2(n-m)}\lambda_\pm)  \big)^{-1/2},\quad l<0,\\
&  c_{nl}^\pm= \gamma_\pm^{-1/2} q^{-n-l}\big(
\mathop{\mbox{$\prod$}}_{m=1}^{l}(q^{2(n+m)}\lambda_\pm-\lambda_-)
(\lambda_+-q^{2(n+m)}\lambda_\pm) \big)^{-1/2},\quad l>0.
\end{align*}
For $r=0$, only $\vartheta^+_{nl}$ is considered. 
The set 
$\{ \vartheta^{\pm}_{nl}\,;\, 
n\hsp\in\hsp \dN_0,\, l\hsp\in\hsp\dZ,\, l\ge -n,\}$ 
is an orthonormal basis of $\cL_2(\mathrm{S}^2_{qr})_\pm$.  
Note that $c_\pm (n)=\lambda_n (q^{2n} \lambda^2_\pm +r)^{1/2}
=((q^{2n} \lambda_\pm -\lambda_-)
(\lambda_+ - q^{2n}\lambda_\pm))^{1/2}$. 
Using the commutation rules in $\cF_{\mathrm{b}}(\mathrm{S}^2_{qr})$ 
and applying 
$f(A)\chi_{n,\pm}(A)=f(q^{2n}\lambda_\pm)\chi_{n,\pm}(A)$ 
for all $f\in \cF_{\mathrm{b}}(\sigma(A))$, 
we get
\begin{align*}
&\rho_{j,\pm}  (A)\vartheta^\pm_{nl} =  
\lambda_\pm  q^{2n} \vartheta^\pm_{nl},  \\
&\rho_{j,\pm}  (B)\vartheta^\pm_{nl} = 
c_{\pm}(n) \vartheta^\pm_{n-1,l+1},   \qquad  
\rho_{j,\pm} (B^\ast) \vartheta^\pm_{nl} = 
c_{\pm}(n\!+\!1) \vartheta^\pm_{n+1,l-1},\\
&\rho_{j,\pm}  (Y)\vartheta^\pm_{nl} =  
\lambda_{\pm} q^{\pm 2j+1}q^{2(n+l)} \vartheta^{\pm}_{nl},\\
&\rho_{j,\pm}  (X)\vartheta^\pm_{nl} = 
\lambda_{n+l+1}\big(q^{2(n+l)}
(\lambda_{\pm}q^{\pm 2j+1})^2+r\big)^{1/2} 
\vartheta^\pm_{n,l+1},\\
&\rho_{j,\pm} (X^\ast) \vartheta^\pm_{nl} = 
\lambda_{n+l}\big(q^{2(n+l-1)}
(\lambda_{\pm}q^{\pm 2j+1})^2+r\big)^{1/2}
\vartheta^\pm_{n,l-1}.
\end{align*}
Renaming $\zeta^\pm_{nk} := \vartheta^\pm_{n,k-n}$, $k,n\in\dN_0$, 
we obtain
\begin{align*}
&\rho_{j,\pm} (A)\zeta^\pm_{nk} =  
\lambda_\pm  q^{2n} \zeta^\pm_{nk},\ \ 
\rho_{j,\pm} (B)\zeta^\pm_{nk} = 
c_{\pm}(n) \zeta^\pm_{n-1,k},\\
&\rho_{j,\pm}(B^\ast)\zeta^\pm_{nk} =  
c_\pm (n\!+\!1)\zeta^\pm_{n+1,k},\\
&\rho_{j,\pm} (Y)\zeta^\pm_{nk} \hsp=\hsp  
q^{2k} q^{\pm 2j+1}\lambda_\pm \zeta^{\pm}_{nk},\ \,\hspace{-1.2pt}
\rho_{j,\pm} (X)\zeta^\pm_{nk} \hsp=\hsp 
\lambda_{k+1}(q^{2k} (q^{\pm 2j+1}\lambda_\pm )^2\hsp+\hsp r)^{1/2} 
\zeta^\pm_{n,k+1},\\
&\rho_{j,\pm}(X^\ast) \zeta^\pm_{nk} = 
\lambda_k (q^{2(k-1)} (q^{\pm 2j+1}\lambda_\pm)^2+r)^{1/2} 
\zeta^\pm_{n,k-1}.
\end{align*}
Let $\Hh^{\pm}$ and $\cK^{\pm}$ 
be the closed linear spans of orthonormal systems 
$\{\eta^\pm_{n}\,;\,n\hspace{1.5pt}{\in}\hspace{1.5pt} \dN_0\}$ and 
$\{\xi^\pm_{k}\,;\, k\in \dN_0\}$, respectively. 
Setting $\zeta^\pm_{nk}:= \eta^\pm_{n}\otimes\xi^\pm_{k}$,  
we see that the restriction of $\rho_{j,\pm}$ 
to $\cO(\mathrm{S}^2_{qr})\otimes \cY_r$ 
is unitarily equivalent to the 
tensor product representation
$\sigma^\pm\otimes\sigma_j^\pm$.
The last assertion of the theorem follows immediately from the 
preceding since $\cO(\mathrm{S}^2_{qr})\otimes 
\cY_r\subset \cF_{\mathrm{b}}(\mathrm{S}^2_{qr}) \otimes \cY_r$.  
                                                         \hfill$\Box$
%
\section{Description\hspace{-0.2pt} of\hspace{-0.2pt} 
the\hspace{-0.2pt} irreducible\hspace{-0.2pt} integrable 
representations by the second approach} 
                                                        \label{last}
As we have seen in Subsection \ref{sec-repcpa}, 
the irreducible integrable representation $\check{\pi}_j$ 
of $\cO(\mathrm{S}^2_{qr})\rti \cU_q (\mathrm{su}_2)$ leads to 
bounded $\ast$-re\-pre\-sen\-ta\-tions  $\rho_{j,\pm}$ 
of $\cO(\mathrm{S}^2_{qr}) \otimes \cY_r$ 
on the charts $\cF_{\mathrm{b}}(\mathrm{S}^2_{qr})_\pm$.
In this section, we recover the irreducible integrable representations 
$\check\pi_j$ of $\cO(\mathrm{S}^2_{qr})\rti \cU_q (\mathrm{su}_2)$ 
by taking the orthogonal sum of both charts and passing to another domain. 
Because this construction is based on 
$\ast$-re\-pre\-sen\-ta\-tions of 
$\cO(\mathrm{S}^2_{qr}) \otimes \cY_r$, 
we say that we have described the irreducible integrable 
representation $\check\pi_j$ by the second approach.

We begin by showing 
that the bounded $\ast$-re\-pre\-sen\-ta\-tion $\rho_{j,\pm}$ of 
$\cO(\mathrm{S}^2_{qr})\otimes \cY_r$ on the 
Hilbert space completion $\cL_2(\mathrm{S}^2_{qr})_\pm$ 
of $\cF_{\mathrm{b}}(\mathrm{S}^2_{qr})_\pm$
leads to a $\ast$-re\-pre\-sen\-ta\-tion
of the $\ast$-al\-ge\-bra 
$\hat{\cO}(\mathrm{S}^2_{qr})\rti \cU_q (\mathrm{su}_2)$. 
By a slight abuse of notation, 
we use the same symbol $\rho_{j,\pm}$ 
to denote the representation of the cross product algebra.
It is a $\ast$-re\-pre\-sen\-ta\-tion by unbounded operators acting 
on the invariant dense domain 
\begin{equation*}        
\cD_{j,\pm} 
:= \cap^\infty_{n,m=0} \cD(\rho_{j,\pm}(A)^{-n}\rho_{j,\pm}(Y)^{-m})
\subset \cL_2(\mathrm{S}^2_{qr})_\pm
\end{equation*} 
Note that $\pi(A)$ and $\pi(Y)$ 
are commuting
bounded self-adjoint operators  
but their inverses are unbounded.
For $\varphi \in \cD_{j,\pm}$, define
\begin{equation} \label{defka}
\rho_{j,\pm}(K)\varphi
:= q^{1/2} |\rho_{j,\pm}(Y)|^{-1/2}  |\rho_{j,\pm}(A)|^{1/2} \varphi ,\quad 
\rho_{j,\pm}(A^{-1})\varphi:=\rho_{j,\pm}(A)^{-1}\varphi
\end{equation}
and $\rho_{j,\pm}(K^{-1}):=\rho_{j,\pm}(K)^{-1}$. 
When $\varphi$ belongs to 
$\cF_{\mathrm{b}}(\mathrm{S}^2_{qr})_{\pm}\cap\cD_{j,\pm}$, 
we can write 
$\rho_{j,\pm}(K)\varphi=q^j |A|^{1/2} \varphi|A|^{-1/2}$ 
and $\rho_{j,\pm}(A^{-1})\varphi:=A^{-1}\varphi$. 
From the commutation rules in the algebra 
$\hat{\cO}(\mathrm{S}^2_{qr})\otimes \cY_r$, 
it follows easily that (\ref{defka}) defines indeed a 
$\ast$-re\-pre\-sen\-ta\-tion of the larger 
$\ast$-al\-ge\-bra 
$\hat{\cO}(\mathrm{S}^2_{qr}) \rti \cU_q (\mathrm{su}_2)$ on $\cD_{j,\pm}$.

Using the isometric isomorphism $\check{\Psi}_{j,\pm}$, we obtain 
an irreducible $\ast$-re\-pre\-sen\-ta\-tion 
$\check{\rho}_{j,\pm}:= 
\check{\Psi}_{j,\pm}\circ \rho_{j,\pm}\circ 
(\check{\Psi}_{j,\pm})^{-1}$ 
of $\hat{\cO}(\mathrm{S}^2_{qr}) \rti \cU_q (\mathrm{su}_2)$
on
$$
\cD(\check{\rho}_{j,\pm})
:= \cap^\infty_{n,m=0} \cD(\check{\rho}_{j,\pm}(A)^{-n}
\check{\rho}_{j,\pm}(Y)^{-m})\subset \cL_2(\mathrm{S}^2_{qr})_\pm^{2|j|+1}P_j, 
$$
The restriction of the direct sum 
$\check{\rho}_{j}:=\check{\rho}_{j,-}\oplus\check{\rho}_{j,+}$ 
to $\cO(\mathrm{S}^2_{qr})\otimes \cY_r$ yields a bounded representation 
which can be extended to a representation on the Hilbert space 
$\cL_2(\mathrm{S}^2_{qr})^{2|j|+1}P_j
=\cL_2(\mathrm{S}^2_{qr})^{2|j|+1}_-P_j\oplus
\cL_2(\mathrm{S}^2_{qr})^{2|j|+1}_+P_j$. 
By the definitions of $\rho_{j,\pm}$ and $\check{\rho}_{j,\pm}$, 
it is obvious that this representation coincides with 
${{\check \pi}}^\mathrm{b}_j$ and so it coincides with 
the restriction of ${{\check \pi}}_j$ 
to $\cO(\mathrm{S}^2_{qr})\otimes \cY_r$ 
on its common domain.  
However, the restriction of $\check{\rho}_{j}$
to $\cO(\mathrm{S}^2_{qr}) \rti \cU_q (\mathrm{su}_2)$ is
not unitarily equivalent to $\check{\pi}_j$ since the latter 
is irreducible while the former is not. 
Theorem \ref{adrep} below shows that
we can reconstruct the irreducible integrable  
representations $\check{\pi}_j\cong \hat{\pi}_j\cong \pi_j$ 
from $\check{\rho}_{j}$. 
The proof of Theorem \ref{adrep} is based on the following lemma. 
\begin{thl}                                                  \label{EXKY}
Let $\eta\in\cO(\mathrm{S}^2_{qr})^{2|j|+1}P_j$. Then
\begin{align}                                         \label{EX}
& \hspace{2.5cm}{{\check \pi}}^\mathrm{b}_j(A){{\check \pi}}_j(K^{-1}E)\eta=
q^{-3/2}\lambda^{-1}{{\check \pi}}^\mathrm{b}_j(X^\ast-qB^\ast)\eta,\\
&\hspace{2.5cm}{{\check \pi}}^\mathrm{b}_j(A){{\check \pi}}_j(FK^{-1})\eta=
q^{-3/2}\lambda^{-1}{{\check \pi}}^\mathrm{b}_j(X-q^{-1}B)\eta, \label{FX}\\
&|{{\check \pi}}^\mathrm{b}_j(Y)|^{1/2}{{\check \pi}}_j(K)\eta\hsp=\hsp
q^{1/2}|{{\check \pi}}^\mathrm{b}_j(A)|^{1/2}\eta, \ \ 
|{{\check \pi}}^\mathrm{b}_j(A)|^{1/2}{{\check \pi}}_j(K^{-1})\eta\hsp=\hsp
q^{-1/2}|{{\check \pi}}^\mathrm{b}_j(Y)|^{1/2}\eta.            \label{KY}
\end{align}
\end{thl}
{\bf Proof.} 
Let $\eta=(y_{-|j|},\ldots,y_{|j|})P_j\in\cO(\mathrm{S}^2_{qr})^{2|j|+1}P_j$. 
From the uniqueness of the square root of a positive operator, 
it follows that
\begin{align*} 
&|{{\check \pi}}^\mathrm{b}_j(A)|^{1/2}(y_{-|j|},\ldots,y_{|j|})P_j
=(|A|^{1/2}y_{-|j|},\ldots,|A|^{1/2}y_{|j|})P_j,\\
&|{{\check \pi}}^\mathrm{b}_j(Y)|^{1/2}(y_{-|j|},\ldots,y_{|j|})P_j
=q^{1/2}(q^{|j|}y_{-|j|}|A|^{1/2},\ldots,q^{-|j|}y_{|j|}|A|^{1/2})P_j.
\end{align*}
Further, the commutation rules in $\cF_{\mathrm{b}}(\mathrm{S}^2_{qr})$ 
imply for all $y\in\cO(\mathrm{S}^2_{qr})$ 
$$
|A|^{1/2}y=(K \ang y) |A|^{1/2},\ \quad   
y|A|^{1/2}=|A|^{1/2} (K^{-1}\ang y). 
$$ 
On the other hand, since  
${\check{\pi}}_j:=\Psi_j^{-1}\circ\hat\pi_j\circ\Psi_j$, 
we have 
\begin{align*}
{{\check \pi}}_j(K^{\pm 1})(y_{-|j|},\ldots,y_{|j|})P_j
&=[2|j|\!+\!1]^{-1/2}\Psi_j^{-1}
(K^{\pm 1}\ang \sumop^{|j|}_{k=-|j|} q^{-k} y_k v^{|j|}_{kj})\\
&=(q^{\mp |j|}K^{\pm 1}\ang y_{-|j|},\ldots,
q^{\pm |j|}K^{\pm 1}\ang y_{|j|})P_j.
\end{align*}
Combining the preceding equations proves (\ref{KY}).

Computing the action of 
${{\check \pi}}_j(FK^{-1})=\Psi_j^{-1}\circ\hat\pi_j\circ\Psi_j(FK^{-1})$
on the element  $\eta = (y_{-|j|},\ldots,y_{|j|})P_j$ gives
\begin{align*}
{{\check \pi}}_j(FK^{-1})(y_{-|j|},\ldots,y_{|j|})P_j
&=[2|j|\!+\!1]^{-1/2}\Psi_j^{-1}
\sumop^{|j|}_{k=-|j|}q^{-k}\big( (FK^{-1}\ang y_k) v^{|j|}_{kj} \\
&\qquad
+q^{-k}[|j|\!-\!k\!+\!1]^{1/2}[|j|\!+\!k]^{1/2}(K^{-2}\ang y_k)v^{|j|}_{k-1,j}
\big)
\end{align*}
so that, by using Lemma \ref{1/1},
\begin{align*}
&{{\check \pi}}^\mathrm{b}_j(A){{\check \pi}}_j(FK^{-1})
(y_{-|j|},\ldots,y_{|j|})P_j = 
(A(FK^{-1}\ang y_{-|j|}),\ldots,A(FK^{-1}\ang y_{|j|}))P_j \\
&+\!q^{-2}(q^{|j|}([2|j|][1])^{1/2}A(K^{-2}\ang y_{-|j|+1}),\ldots,
q^{-|j|+1}([1][2|j|])^{1/2}A(K^{-2}\ang y_{|j|}),0)P_j\\
&=-q^{-5/2}\lambda^{-1}( [B,y_{-|j|}],\ldots,[B,y_{|j|}])P_j \\
&\hspace{2.3cm}+q^{-2}(q^{|j|}([2|j|][1])^{1/2}y_{-|j|+1}A,\ldots,
q^{-|j|+1}([1][2|j|])^{1/2}y_{|j|}A,0)P_j.
\end{align*}
Comparing the last identity with the action of 
${{\check \pi}}^\mathrm{b}_j(X)$ and ${{\check \pi}}^\mathrm{b}_j(B)$ 
on $\cO(\mathrm{S}^2_{qr})^{2|j|+1}P_j$ from Subsection \ref{sec-repcpa} 
(see the discussion preceding Proposition \ref{6.4}) 
shows (\ref{FX}). Equation (\ref{EX}) is proved similarly. 
                                              \hfill$\Box$
\mn

Let us recall the notion of the adjoint of a $\ast$-re\-pre\-sen\-ta\-tion 
$\pi$ of a $\ast$-al\-ge\-bra $\cX$ (see e.g.\  \cite[Section 8.1]{S1} ).  
It is a representation $\pi^\ast$  of $\cX$ acting on the domain 
$\cD(\pi^\ast):= \cap_{x\in\cX} \cD (\pi (x)^\ast)$  
by $\pi^\ast (x)\varphi = \pi (x^\ast)^\ast \varphi$, where  
$\varphi\in\cD(\pi^\ast)$. 
In general, $\pi^\ast$ is not a $\ast$-re\-pre\-sen\-ta\-tion.
\begin{tht}                                              \label{adrep}
The irreducible integrable $\ast$-re\-pre\-sen\-ta\-tion ${\check \pi}_j$
of $\cO(\mathrm{S}^2_{qr})\rti \cU_q(\mathrm{su}_2)$ 
is the restriction of the adjoint $\check{\rho}_{j}^\ast$ of the 
$\ast$-re\-pre\-sen\-ta\-tion 
$\check{\rho}_{j}=\check{\rho}_{j,-}\oplus\check{\rho}_{j,+}$ of 
$\hat{\cO}(\mathrm{S}^2_{qr})\rti \cU_q(\mathrm{su}_2)$ 
to the domain $\cO(\mathrm{S}^2_{qr})^{2|j|+1}P_j$.
There exists a Hilbert space basis 
$\{\eta^\pm_{nm}\,;\,n,m\in\dN_0\}$ of 
$\cL_2(\mathrm{S}^2_{qr})^{2|j|+1}_\pm P_j$ such that 
$\check{\rho}_{j,\pm}$ is determined by the formulas 
$(I)_{\pm,q^{\mp j}}$ from Subsection \ref{sec-4}. 
\end{tht}
{\bf Proof.} 
First note that, by definition, 
$\check{\rho}_{j}(Z)\varphi ={{\check \pi}}^\mathrm{b}_j(Z)\varphi$ 
for all $Z\in\cO(\mathrm{S}^2_{qr})\otimes \cY_r$ and 
$\varphi\in\cD(\check{\rho}_{j,\pm})$. In particular, 
$\langle \check{\rho}_{j}(A^{-1})\varphi,
{{\check \pi}}^\mathrm{b}_j(A)\eta\rangle=\langle \varphi,\eta\rangle$ 
for all  $\eta$ from the Hilbert space 
$\cL_2(\mathrm{S}^2_{qr})^{2|j|+1}P_j$. 
Now let $\eta\in\cO(\mathrm{S}^2_{qr})^{2|j|+1}P_j$ and 
$\varphi\in\cD(\check{\rho}_{j,\pm})$. 
Then, by Lemma \ref{EXKY}, 
\begin{align*}                                        
\langle  \check{\rho}_{j}(K^{-1} E)\varphi,\eta\rangle 
&= q^{-3/2} \lambda^{-1} 
\langle \check{\rho}_{j}(X^\ast \!-\!q^{-1} B^\ast)
\check{\rho}_{j} (A^{-1}) \varphi, \eta\rangle   \\
&= q^{-3/2} \lambda^{-1} \langle \check{\rho}_{j} (A^{-1}) \varphi, 
{{\check \pi}}^\mathrm{b}_j (X \!-\! q^{-1} B) \eta\rangle   \\
&= \langle  \check{\rho}_{j}(A^{-1}) \varphi, 
  {{\check \pi}}^\mathrm{b}_j(A){{\check \pi}}_j(FK^{-1})\eta\rangle   
= \langle \varphi, {{\check \pi}}_j(FK^{-1})\eta\rangle.
\end{align*}
Similarly one shows 
$\langle  \check{\rho}_{j}(FK^{-1})\varphi,\eta\rangle 
= \langle \varphi, {{\check \pi}}_j(K^{-1} E)\eta\rangle$. 
As a above, we have 
$\langle |\check{\rho}_{j}(Y)|^{-1/2}\varphi,
|{\check{\pi}}^\mathrm{b}_j(Y)|^{1/2}\eta\rangle=\langle\varphi,\eta\rangle$
and thus, again by Lemma \ref{EXKY}, 
\begin{align*}                                        
\langle  \check{\rho}_{j}(K)\varphi,\eta\rangle 
&\hspace{-1.3pt}=\hspace{-1.3pt}q^{1/2}\langle |\check{\rho}_{j}(A)|^{1/2}
|\check{\rho}_{j}(Y)|^{-1/2}\varphi,\eta\rangle 
\hspace{-1.3pt}=\hspace{-1.3pt}
q^{1/2}\langle|\check{\rho}_{j}(Y)|^{-1/2}\varphi, 
|{\check{\pi}}^\mathrm{b}_j(A)|^{1/2}\eta\rangle\\ 
&\hspace{-1.3pt}=\hspace{-1.3pt}\langle|\check{\rho}_{j}(Y)|^{-1/2}\varphi, 
|{{\check \pi}}^\mathrm{b}_j(Y)|^{1/2}{{\check \pi}}_j(K)\eta\rangle
\hspace{-1.3pt}=\hspace{-1.3pt}\langle\varphi, {{\check \pi}}_j(K)\eta\rangle.
\end{align*}
Likewise, $\langle  \check{\rho}_{j}(K^{-1})\varphi,\eta\rangle=
\langle\varphi, {{\check \pi}}_j(K^{-1})\eta\rangle$. 
Clearly, $\langle  \check{\rho}_{j}(x^\ast)\varphi,\eta\rangle=
\langle\varphi, {{\check \pi}}_j(x)\eta\rangle$ 
for all $x\in\cO(\mathrm{S}^2_{qr})$. 
As the elements $K^{\pm 1}$, $K^{-1} E$, $FK^{-1}$ and 
$x\in\cO(\mathrm{S}^2_{qr})$ generate the algebra 
$\cO(\mathrm{S}^2_{qr})\rti \cU_q(\mathrm{su}_2)$, it follows that 
$\cO(\mathrm{S}^2_{qr})^{2|j|+1}P_j$ is contained in 
$\cD(\check{\rho}_{j}^\ast)$ and 
$\check{\rho}_{j}^\ast(Z)\eta={{\check \pi}}_j(Z)\eta$ for 
all $Z\in\cO(\mathrm{S}^2_{qr})\rti \cU_q(\mathrm{su}_2)$
and $\eta\in\cO(\mathrm{S}^2_{qr})^{2|j|+1}P_j$.

Computing the action of $\rho_{j,\pm}(K)$, 
$\rho_{j,\pm}(E)\hspace{-1.1pt}=\hspace{-1.1pt}q^{-3/2} \lambda^{-1} 
\rho_{j,\pm}(A^{-1} K(X^\ast-q B^\ast))$, 
$\rho_{j,\pm}(F)=\rho_{j,\pm}(E)^\ast$
on the basis vectors $\zeta^\pm_{nm}$ from the proof 
of Theorem \ref{reprho} shows that $\rho_{j,\pm}$ 
is determined by the formulas 
$(I)_{\pm,q^{\mp j}}$ from 
Subsection \ref{sec-4}. Setting 
$\eta^\pm_{nm}=\check{\Psi}_{j,\pm}(\zeta^\pm_{nm})$ 
establishes the second assertion of the theorem. 
                                                             \hfill$\Box$
\mn

Let us make the case $j=0$  more explicit. 
Then $P_0=1$, $M_0\cong \cO(\mathrm{S}^2_{qr})$, 
and $\pi_0\cong {{\check \pi}}_0$ 
is just the Heisenberg representation of 
$\cO(\mathrm{S}^2_{qr})\rti \cU_q(\mathrm{su}_2)$. 
Accordingly,
${{\check \pi}}_0(f)\eta=f\ang \eta$ for $f\in\cU_q(\mathrm{su}_2)$, 
$\eta\in\cO(\mathrm{S}^2_{qr})$. 
The $1\times 1$-matrices $\fm_0$, $\fm_0^\dagger$ and $\fn_0$ 
from Subsection \ref{sec-repcpa}
have the entries 
$q^{-1}B$, $qB^\ast$ and $qA$, respectively. 
Hence we obtain for 
$z,\varphi \in \cF_{\mathrm{b}}(\mathrm{S}^2_{qr})$
\begin{equation}                                             \label{pidef}
{{\check \pi}}^\mathrm{b}_0(z)\varphi =z\varphi, \quad 
{{\check \pi}}^\mathrm{b}_0(X)\varphi = q^{-1} \varphi B,\quad  
{{\check \pi}}^\mathrm{b}_0 (X^\ast)\varphi = q \varphi B^\ast,\quad  
{{\check \pi}}^\mathrm{b}_0 (Y) \varphi = q \varphi A. 
\end{equation}
Inserting these formulas into (\ref{EX})--(\ref{KY}), 
we recover Equations (\ref{eactt1})--(\ref{eactt3})
which we used to define a $\cU_q(\mathrm{su}_2)$-action on 
the operator algebra $\cL^+(\cD)$.
In particular, Equations (\ref{eactt1})--(\ref{eactt3}) and (\ref{pidef})
give a new description of the Heisenberg representation 
of $\cO(\mathrm{S}^2_{qr})\rti \cU_q(\mathrm{su}_2)$ on 
$ \cO(\mathrm{S}^2_{qr})\subset \cF_{\mathrm{b}}(\mathrm{S}^2_{qr})$ 
by left and right multiplications.


\begin{thebibliography}{999} 

\bibitem[1]{B}
Brzezi\'nski, T., {\it Quantum homogeneous spaces as quantum quotient spaces.} 
J. Math. Phys. {\bf 37} (1996), 2388--2399.

\bibitem[2]{BM}
Brzezi\'nski, T. and S. Majid, 
{\it Quantum geometry of algebra factorisations and coalgebra bundles.} 
Commun. Math. Phys. {\bf 213} (2000), 491--521.

\bibitem[3]{DK}
Dijkhuizen, M.  and T. K. Koornwinder, 
{\it Quantum homogeneous spaces, quantum duality and quantum 2-spheres.}
Geometriae Dedicata {\bf 52} (1994), 291--315.

\bibitem[4]{F}
Fiore, G., {\it On the decoupling of the homogeneous and inhomogeneous parts 
in inhomogeneous quantum groups.} J. Phys. A {\bf 35} (2002), 657--678.

\bibitem[5]{H}
Hajac, P. M., 
{\it Bundles over quantum sphere and noncommutative index theorem.} 
K-Theory {\bf 21} (1996), 141--150.

\bibitem[6]{HM}
Hajac, P. M. and S. Majid, 
{\it Projective module description of the $q$-mono\-pole.} 
Commun. Math. Phys. {\bf 206} (1999), 247--264.

\bibitem[7]{KS}
Klimyk, K. A. and K. Schm\"udgen, 
{\it Quantum Groups and Their Representations.} Springer, Heidelberg, 1997.


\bibitem[8]{KR}
Kulish, P. P. and N. Yu. Reshetikhin, 
{\it Quantum linear problem for the sine-Gordon equation 
and higher representations.} 
Zap. Nauchn. Sem. L0MI {\bf 101} (1981), 101--110.

\bibitem[9]{LR}
Lunts, V. A. and A. L. Rosenberg, 
{\it Differential operators on noncommutative rings.} 
Sel. math. {\bf 3} (1997), 335--359.

\bibitem[10]{NM}
Mimachi, K.  and M. Noumi, 
{\it Quantum 2-spheres and big $q$-Jacobi polynomials.}
Commun. Math. Phys. {\bf 128} (1990), 521--531.

\bibitem[11]{M}
Montgomery, S., {\it Hopf algebras and their actions on rings.} 
Amer. Math. Soc., Providence, R.I., 1993.

\bibitem[12]{MS}
M\"uller, E. F. and H.-J. Schneider, 
{\it Quantum homogeneous spaces with faithfully flat module structures.} 
Israel J. Math. {\bf 111} (1999), 157--190.


\bibitem[13]{P}
Podle\'s, P., {\it Quantum spheres.} 
Lett. Math. Phys. {\bf 14} (1987), 193--202.

\bibitem[14]{S1}
Schm\"udgen, K., 
{\it Unbounded Operator Algebras and Representation Theory.} 
Birkh\"auser, Basel, 1990.

\bibitem[15]{SW}
Schm\"udgen, K. and E. Wagner, 
{\it Hilbert space representations of cross product algebras.} 
J. Funct. Anal. {\bf 200} (2003), 451--493. 

\bibitem[16]{S2}
Schneider, H.-J., 
{\it Principal homogeneous spaces for arbitrary Hopf algebras.} 
Israel J. Math. {\bf 72} (1990), 167--195.

\bibitem[17]{VS}
Soibelman, Ya. S.  and L. L. Vaksman,  
{\it Algebra of functions on the quantum SU(2).}
Funct. Anal. Appl. {\bf 22} (1988), 170--181.

\bibitem[18]{SF}
Sz.-Nagy, B. and C. Foias, 
{\it Analyse harmonique des operateurs de l'espace de Hilbert.} 
Academiai Kiado, Budapest, 1979.

\bibitem[19]{W}
Woronowicz, S. L., 
{\it Twisted SU(2) group. 
An example of a non-commutative differential calculus.} 
Publ. RIMS Kyoto Univ. {\bf 23} (1987), 117--181.

\end{thebibliography}
\end{document}